\documentclass[11pt,reqno]{amsart}

\usepackage{amsmath,amsfonts,amssymb,mathrsfs}
\usepackage{amssymb,mathrsfs,graphicx,extpfeil}
\usepackage{amsmath,amsfonts,amssymb,amscd,amsthm,bbm}
\usepackage{epsfig}

\usepackage{graphicx,colortbl,here}
\usepackage{extpfeil}
\usepackage{subcaption}
\usepackage{caption}
\usepackage{hyperref}
\usepackage{etoolbox}

\newcommand{\bbr}{\mathbb R}

\newcommand{\E}{\mathbb E}
\newcommand{\R}{\mathbb R}
\newcommand{\bx}{\mathbf x}
\newcommand{\bv}{\mathbf v}

\newcommand{\by}{\mathbf y}
\newcommand{\bz}{\mathbf z}
\newcommand{\bh}{\mathbf h}
\newcommand{\bF}{\mathbf F}
\newcommand{\bG}{\mathbf G}

\newcommand{\bphi}{\boldsymbol{\phi}}
\newcommand{\bpsi}{\boldsymbol{\psi}}
\newcommand{\bDelta}{\boldsymbol{\Delta}}
\newcommand{\bPhi}{\boldsymbol{\Phi}}
\newcommand{\bPsi}{\boldsymbol{\Psi}}

\newcommand{\bX}{\mathbf{X}}
\newcommand{\bY}{\mathbf{Y}}
\newcommand{\bZ}{\mathbf{Z}}
\newcommand{\bV}{\mathbf{V}}

 \newcommand{\one}{\mathbf 1}
 \newcommand{\dd}{\,\mathrm d}
\newcommand{\norm}[1]{\left\|#1\right\|}

\newcommand{\Diag}{\operatorname{Diag}}

\newcommand*\di{\mathop{}\!\mathrm{d}}

\makeatletter
\def\@setkeywords{%
	{\itshape \keywordsname.}\enspace
	\@keywords\@addpunct.\par\vspace{0cm}}
\makeatother

\newcommand{\opnorm}[1]{{ \vert\kern-0.25ex \vert\kern-0.25ex \vert #1 
   \vert\kern-0.25ex \vert\kern-0.25ex \vert}}

\newtheorem{theorem}{Theorem}[section]
\newtheorem{lemma}{Lemma}[section]

\newtheorem{proposition}{Proposition}[section]
\newtheorem{remark}{Remark}[section]
\newtheorem{definition}{Definition}[section]

\renewcommand{\theequation}{\thesection.\arabic{equation}}

\allowdisplaybreaks[4]

\makeatletter

\patchcmd{\@setauthors}
{\MakeUppercase{\authors}}
{\authors}
{}{\PackageError{frontaff}{Failed to patch author capitalization}{}}

\newcommand{\frontaffiliations}{%
	\vspace{0.35em}
	\begin{center}
		\footnotesize\itshape
		\begin{minipage}{0.88\textwidth}
			\raggedright
			$^{a}$ Department of Mathematical Sciences and Research Institute of Mathematics,
			Seoul National University,
			Seoul 08826, Republic of Korea\\[1pt]
			
			$^{b}$ School of Mathematics, Harbin Institute of Technology,
			Harbin 150001, People's Republic of China\\[1pt]
			
			$^{c}$ School of Mathematical Sciences and Big Data,
			Anhui University of Science and Technology,
		Huainan 232001, People's Republic of China\\[1pt]
			
			$^{d}$ College of Mathematical Sciences,
			Harbin Engineering University,
			Harbin 150001, People's Republic of China
		
		\end{minipage}
	\end{center}
	\vspace{0.25em}
}

\patchcmd{\@maketitle}
{\ifx\@empty\authors \else \@setauthors \fi}
{\ifx\@empty\authors \else \@setauthors \fi
	\frontaffiliations}
{}{\PackageError{frontaff}{Failed to patch \string\@maketitle}{}}

\makeatother

\begin{document}

\title[Quantitative flocking for the DSCS model]{Quantitative flocking of Cucker-Smale flocks under time-delay and noise}

\author{
	Seung-Yeal Ha$^{a}$,
	Xinyu Wang$^{b}$,
	Chen Wu$^{c}$,
	and Fanqin Zeng$^{d,*}$
}

\thanks{$^{*}$ Corresponding author.
	E-mail addresses:
	syha@snu.ac.kr (S.-Y. Ha),
	wangxinyumath@hit.edu.cn (X. Wang),
	wuchen@aust.edu.cn (C. Wu),
	fqzeng@hrbeu.edu.cn (F. Zeng).}

\vspace{0.2cm}

\begin{abstract}
	We study the asymptotic flocking dynamics of the delayed stochastic Cucker-Smale (in short, DSCS)  model driven by agent-wise independent, component-wise multiplicative Brownian noise. The simultaneous interplay of stochastic perturbations and time delay yields fundamental challenges to asymptotic stochastic flocking analysis, as the communication weight may lose a uniform positive lower bound, while the velocity moment is no longer monotone. To overcome these difficulties, we present a sufficient framework and an analytical tool combining the
	method of time-varying effective region and Lyapunov functional approach with
	time-delay corrections. This approach avoids any reliance on uniformly positive communication weights or monotone velocity moments. The proposed approach recovers sufficiently strong effective velocity alignment despite stochastic fluctuations and delay effects. As a consequence, we establish the mean-square stochastic flocking with an explicit polynomial convergence rate. The admissible decay exponent is determined by the communication kernel exponent and the orders of the spatial and velocity history moments. These results establish a quantitative theory for long-time stochastic flocking under simultaneous delay and multiplicative noise, and we provide a positive answer to an open problem raised by~[R. Erban, J. Haskovec, Y. Sun, SIAM J. Appl. Math., 76, 1535--1557 (2016)]. We also develop robust analytical tools applicable to delayed stochastic collective dynamics with degenerate communications.
\end{abstract}

\keywords{Cucker-Smale model, multiplicative noise, stochastic flocking, time-delay}

\subjclass[2020]{60H10, 34K50, 60K35}

\maketitle

\section{Introduction} \label{sec:1}
\setcounter{equation}{0}
Collective behaviors are ubiquitous in interacting particle systems, where simple communications between agents may lead to coherent macroscopic patterns such as alignment, aggregation, synchronization, and flocking. Such emergent behaviors have often been observed in biological swarms, animal groups, robotic networks, opinion dynamics, and other multi-agent systems \cite{C-S,H-L,T-B,V-C,W2}. Over the past several decades,  mathematical models for such collective dynamics have attracted considerable attention, leading to significant progress on the local and global well-posedness, asymptotic behavior, and mean-field description of various alignment models \cite{C-F-R-T,C-S,H-W-X,M-T,S-T1,S-T2,V-C}. Among them, the Cucker-Smale (CS) model has become one of the most influential paradigms for self-organized flocking through pairwise velocity alignment. In this paper, we focus on the quantitative flocking analysis of a delayed stochastic Cucker-Smale (DSCS) model. To set the stage, we begin with the deterministic CS model \cite{C-S,C-S2}. 

Let $(\bx_i(t),\bv_i(t)) \in\R^{2d}$ be a position-velocity pair of the $i$-th CS particle at time $t$. Then, the dynamics of mechanical observables $(x_i, v_i)$ is governed by the CS model:
\begin{equation} \label{A-0}
	\begin{cases}
		\displaystyle \frac{\di \bx_i(t)}{\di t} = \bv_i(t),\quad t > 0, \quad  i\in[N]:=\{1,\ldots,N\},\vspace{6pt}\\
		\displaystyle \frac{\di\bv_i(t)}{\di t} =\frac{\kappa}{N}\sum_{j=1}^N
		\omega( \|\bx_j(t) -\bx_i(t)  \|)(\bv_j(t) -\bv_i(t)).
	\end{cases}
\end{equation}
Here, $\|\cdot \|$ is the standard $\ell^2$-norm in $\mathbb{R}^d$. The constant \(\kappa>0\) is the coupling strength. Since its precise value plays no role in our analysis, we set $\kappa=1$ throughout the paper. The map
\[
r\in[0,\infty)\quad\mapsto\quad\omega(r)\in(0,\infty)
\]
is the communication weight function taking the following ansatz:
\begin{equation}
	\omega(r)=\frac{1}{(1+r^2)^{\beta/2}},
	\quad r\geq0, \quad \beta \geq 0.
	\label{A-1}
\end{equation}
 When $0\le \beta\le1$, the communication weight function is not integrable at infinity and is therefore called a long-range kernel. In this regime, mono-cluster flocking emerges for arbitrary initial configurations \cite{C-S,C-S2,H-L,H-W2,W-X}. When $\beta>1$, the communication weight function is integrable at infinity and we call it a short-range kernel. In this case, multi-cluster flocking can emerge depending on initial configuration (see \cite{H-K-P-Z}). The CS model \eqref{A-0} has been extensively investigated from particle, kinetic, and hydrodynamic perspectives, including well-posedness, asymptotic flocking, mean-field limits, and various extensions of the interaction mechanism, e.g. \cite{C-C-R,C-C-H,C-F-R-T,H-T,H-W-Z,M-T,W} and the references therein. 

In this work, we are interested in the effect on the flocking due to the interplay between time delay and stochastic noise. In fact, the delay and stochastic forcing have both attracted considerable attention in the literature. Delay CS models have been investigated for various interaction mechanisms, while stochastic CS models have also been studied in connection with flocking, singular communication, control, phase transitions, and mean-field limits
\cite{C-H,C-L,C-O-T,C-P,C-S1,C,C-CF,C-M,D-Z,H2,H,H-M-D,H-L2,H-Z,L-Z-H,T-L-Y,Z-L-W-X}, etc.  Among others, we focus on the DSCS model proposed by Erban, Haskovec, and Sun in \cite{E-H-S}. To be consistent with notations in stochastic analysis, we use superscript and subscript to denote the particle number and time, i.e., $\bx^i_t$ and $\bv^i_t$ represent the position and velocity processes of the $i$-th agent at time $t$, respectively. Let
\((\Omega,\mathcal F,(\mathcal F_t)_{t\geq0},\mathbb P)\) be a filtered
probability space, and let $\{ {\bf B}_t^1,\ldots,{\bf B}_t^N \}$ be a family of independent standard \(\R^d\)-valued Brownian
processes adapted to the filtration \((\mathcal F_t)_{t\geq0}\),  where
\(
{\bf B}_t^i
=
\bigl(B_t^{i,1},\ldots,B_t^{i,d}\bigr),
\) with zero mean and covariance:
\[
\E\left[B_t^{i,k}B_s^{j,\ell}\right]
=
\delta_{ij}\delta_{k\ell}\min\{t,s\},
\qquad
1\leq i,j\leq N,\quad
1\leq k,\ell\leq d.
\]
Here \(\delta_{ij}\) is the Kronecker symbol:
\[
\delta_{i j} = \begin{cases}
	1, \quad & i  =j, \\
	0, \quad & i \neq j.
\end{cases}
\]
For \(\bz = (z^1, \cdots, z^d) \in\R^d\), the $d \times d$ matrix \(\Diag(\bz)\) denotes  the diagonal matrix with diagonal \(\bz  \). Then, the stochastic dynamics of $(\bx^i_t, \bv^i_t)$ is governed by the Cauchy problem for the coupled system of stochastic ordinary differential equations with the uniform time-delay \(\tau\):
\begin{equation}
	\begin{cases} \label{A-2}
		\dd \bx^i_t=\bv^i_t\dd t,
		\qquad t\geq0,\quad i\in[N],\\[0.1cm]
		\displaystyle
		\dd \bv^i_t
		=
		\frac 1 N\sum_{j=1}^N
		\omega\bigl( \|\bx^j_{t-\tau}-\bx^i_{t-\tau} \| \bigr)
		\bigl( \bv^j_{t-\tau}-\bv^i_{t-\tau}\bigr)\dd t\\[0.1cm]
		\displaystyle
		\hspace{0.4cm}
		+\sqrt{\sigma^i}\,
		\Diag\!\left(
		\frac1N\sum_{j=1}^N
		\omega\bigl( \|\bx^i_{t-\tau}-\bx^j_{t-\tau} \|\bigr)
		\bigl( \bv^j_{t-\tau}- \bv^i_{t-\tau}\bigr)
		\right)\dd {\bf B}^i_t, \vspace{6pt}\\
		\displaystyle(\bx^i_s, \bv^i_s)
		=
		(\bphi^i_s, \bpsi^i_s),
		\quad
		s\in [-\tau,0],
	\end{cases}
\end{equation}
where \(\sigma^i\geq0\) is the
noise intensity of the \(i\)-th particle, $\bPhi$ and $\bPsi$ are continuous initial histories:
\[
(\bPhi_s,\bPsi_s):= \bigl((\bphi^1_s, \ldots, \bphi^N_s), (\bpsi^1_s,\ldots, \bpsi^N_s )\bigr)
\in
L^2\!\left(
\Omega,\mathcal F_0;
C\!\left([-\tau,0];(\mathbb R^{2d})^N\right)
\right).
\]
That is, the initial history is
\(\mathcal F_0\)-measurable, has continuous sample paths, and satisfies bounded second moment estimates:
\[
\E\left[
\sup_{-\tau\le s\le0}
\Big(
\| \bPhi_s \|^2
+  \|\bPsi_s \|^2
\Big)
\right]
<\infty.
\]

In \cite{E-H-S}, the authors imposed a uniformly positive Fiedler number for the Laplacian matrix associated with communication weights to establish the asymptotic stochastic flocking of system \eqref{A-2}.  For the original configuration-dependent
CS kernel, however, this condition cannot in general be verified
a priori, since spatial spreading may cause the spectral gap to degenerate. As pointed out in \cite{E-H-S}, extending the analysis to the original
CS communication weight function \eqref{A-1} remains an open problem. Motivated by this open question, we are interested in the following questions:
\vspace{0.1cm} 

\begin{quote}
	\begin{itemize}
		\item (Q1):~Can we establish stochastic flocking for the DSCS model \eqref{A-2} with the original CS communication weight \eqref{A-1}?
		\vspace{0.1cm}
		\item (Q2):~If so, under what conditions on system parameters and initial histories can we establish stochastic flocking rigorously?
	\end{itemize}
\end{quote}
The purpose of this paper is to answer the above questions in an affirmative manner. For later use, we set the centroids of space and velocity configurations:
\[
\bar{\bx}_t:=\frac1N\sum_{i=1}^N \bx^i_t, \quad \bar{\bv}_t:=\frac1N\sum_{i=1}^N \bv^i_t.
\]
For (Q1), we introduce two types of stochastic flocking as follows:
\begin{definition}[Stochastic flocking]
	\label{D1.1}
	Let $(\bx_t^i, \bv_t^i)$  be a global solution process to the Cauchy problem \eqref{A-2}.
	\begin{enumerate}
		\item
		The Cauchy problem \eqref{A-2} exhibits mean-square stochastic flocking if the following relations hold.
		\begin{align*}
			\lim_{t\to\infty}
			\E\left[
			\frac1N\sum_{i=1}^N
			\norm{\bv^i_t- \bar{\bv}_t}^2
			\right]=0,\quad
			\sup_{0 \leq t <\infty}
			\E\left[
			\frac1N\sum_{i=1}^N
			\norm{\bx^i_t- \bar{\bx}_t}^2
			\right]<\infty.
		\end{align*}
		\item
		The Cauchy problem \eqref{A-2} exhibits almost-sure stochastic flocking if the following relations hold.
		\begin{align}\label{A-5}
			\lim_{t\to\infty}
			\max_{1\leq i,j\leq N} \| \bv^i_t-\bv^j_t \|=0,\quad
			\sup_{t\geq0}
			\max_{1\leq i,j\leq N} \| \bx^i_t-\bx^j_t \|<\infty,\quad\text{a.s.}
		\end{align}
	\end{enumerate}	
\end{definition}
\noindent To establish the above stochastic flocking estimates, we need to overcome three main difficulties. First, velocity alignment depends on
spatial configuration whereas spatial control itself requires the temporal decay of
the velocity fluctuation, creating a circular dependence. Second, the
delayed alignment field produces a sign-indefinite defect in the velocity moment estimate. Third, the lower bound of communication weight among particles may become arbitrarily small due to stochastic effects, i.e., the communication
weights \eqref{A-1} can become arbitrarily
small with positive probability. The method of time-varying effective region was developed for the kinetic model
with non-compact support \cite{H-W1,H-W,H-W-X}. Motivated by this method, we adopt a similar approach to avoid the
degeneration of the communication weight among particles. However, because of the time-delay effect, the velocity moment is no longer monotone. For this, we develop a new method which combines the method of time-varying effective region and time-delay Lyapunov argument to establish stochastic flocking in the sense of Definition \ref{D1.1}.\newline

For (Q2), we introduce a sufficient framework in terms of system parameters and initial history configurations:
\begin{align*}
	\begin{aligned}
		& 0 \leq \beta < 1,
		\quad
	p\ge 4,\quad p> \frac{2}{1-\beta},\quad
		0 \leq \tau \ll 1,
		\quad
		0 \leq \sigma^i \ll 1,
 \\
		& \E\left[
		\sup_{s\in[-\tau,0]}
		\frac1N\sum_{i=1}^N
		\norm{\bphi^i_s - \bar{\bphi}_s}^p
		\right] < \infty,\quad
		\E\left[
		\sup_{s\in[-\tau,0]}
		\frac1N\sum_{i=1}^N
		\norm{\bpsi^i_s - \bar{\bpsi}_s}^p
		\right] < \infty.
	\end{aligned}
\end{align*}
 Under the above set-up, we establish the mean-square stochastic flocking (see Theorem \ref{T2.3}): 
\[
\E\left[ \frac1N\sum_{i=1}^N
\| \bv^i_t- \bar{\bv}_t \|^2 \right] \lesssim (1+t)^{-\gamma}, \quad \forall~t \geq 0, \quad \sup_{ 0\leq t < \infty} \E\left[ \frac1N\sum_{i=1}^N
\| \bx^i_t- \bar{\bx}_t \|^2 \right]  < \infty,
\]
for any $\gamma>0$ satisfying $\gamma\beta <(p-2)(1-\beta)$. 
These estimates also imply almost-sure stochastic flocking \eqref{A-5} (see Theorem \ref{T2.4}).\newline

The rest of this paper is organized as follows. In Section~\ref{sec:2}, we study the basic properties of the DSCS model \eqref{A-2}, and then present the previous results on the DSCS model, and our main results.
In Section~\ref{sec:3}--Section \ref{sec:5}, we derive several preparatory estimates to be used in the proofs of Theorem \ref{T2.3} and Theorem \ref{T2.4}. In Section \ref{sec:6} and Section \ref{sec:7}, we provide proofs of mean-square stochastic flocking and almost-sure stochastic flocking, respectively. In Section \ref{sec:8}, we provide numerical simulations and compare them with analytical results. Finally, Section \ref{sec:9} is devoted to a brief summary of our main results and discussion of some remaining issues for future work. 

\vspace{0.5cm} 

\noindent\textbf{Gallery of Notations:} Throughout the paper, we use different notations for deterministic function and stochastic process. For example, for the map $t~\mapsto~f(t)$, we use subscript and superscript to denote particle number and coordinate component, i.e., $\bx_i^j(t)$ denotes the $j$-th coordinate of the $i$-th particle at time $t$. In contrast, we use superscripts to denote the particle number, component of coordinate, and double superscript index to denote the pair of particle number and component of its coordinate, respectively, and subscript to denote time-dependence, i.e., $\bx_t^i = (x_t^{i,1}, \cdots, x_t^{i,d}) \in {\mathbb R}^d$. Let $\bx$ and $\by$ be vectors in ${\mathbb R}^d$. Then, we set 
\[ \bx = (x^1, \cdots, x^d), \quad \by = (y^1, \cdots, y^d), \]
and we define inner product and its associated norm as follows.
\[ \bx \cdot \by := \sum_{i=1}^{d} x^i y^i, \quad  \| \bx \| := \sqrt{\bx \cdot \bx} = \Big( \sum_{i=1}^{d} |x^i|^2 \Big)^{\frac{1}{2}}. \]
For $N$-particle configuration vectors $\bX$ and $\bY$ in $({\mathbb R}^d)^N = {\mathbb R}^{dN}$, we set 
\[  \bX = (\bx^1, \cdots, \bx^N), \quad  \bY = (\by^1, \cdots, \by^N), \]
and we define the averaged inner product and norm by
\[ \langle \bX, \bY \rangle_N
:=\frac1N\sum_{i=1}^N  \bx^i \cdot \by^i,
\quad \| \bX \|_N^2
:=\frac1N\sum_{i=1}^N \| \bx^i \|^2.
\]
Moreover, for the matrix \(A=(a^{ij})\in\R^{m\times n}\), we denote by
\(\|\cdot\|_{\mathrm{F}}\) the Frobenius (Hilbert--Schmidt) norm:
\[
\|A\|_{\mathrm{F}}
:=
\left(
\sum_{i=1}^{m}\sum_{j=1}^{n}|a^{ij}|^2
\right)^{1/2}.
\]
Throughout the paper, the constant \(C\) denotes a positive constant independent
of \(N\) and time, whose value may change from line to line. We write
\(C_p\) and \(C_\gamma\) when the dependence on \(p\) or \(\gamma\)
is relevant. We also set
\[ \sigma^\infty :=\max_{1\leq i\leq N}\sigma^i<\infty.
\]
We set the spatial, velocity and phase configurations as 
\[ \bX_t := (\bx_t^1, \cdots, \bx_t^N), \quad  \bV_t := (\bv_t^1, \cdots, \bv_t^N), \quad  \bZ_t := (\bX_t, \bV_t). \]
	
\section{Preliminaries} \label{sec:2}
\setcounter{equation}{0}
In this section, we discuss the global well-posedness of system \eqref{A-2}, and previous results on stochastic flocking in Erban--Haskovec--Sun's work \cite{E-H-S}. Then, we present our main results on the stochastic flocking for the CS flock with the communication weight function \eqref{A-1}.

\subsection{The DSCS model}\label{sec:2.1}
In this subsection, we review the DSCS model with time-delay and recall the previous results on the asymptotic flocking dynamics. Let \((\Omega,\mathcal F,(\mathcal F_t)_{t\geq0},\mathbb P)\) be a filtered probability space. Recall the Cauchy problem for the DSCS model with time-delay:
\begin{equation} \label{B-1}
	\begin{cases}
		\dd \bx^i_t=\bv^i_t\dd t,
		\quad t\geq0,~~i\in[N],\\[0.1cm]
		\displaystyle
		\dd \bv^i_t
		=
		\frac 1 N\sum_{j=1}^N
		\omega\bigl( \|\bx^i_{t-\tau}-\bx^j_{t-\tau} \| \bigr)
		\bigl( \bv^j_{t-\tau}-\bv^i_{t-\tau}\bigr)\dd t\\[0.1cm]
		\displaystyle
		\hspace{0.8cm}
		+\sqrt{\sigma^i}\,
		\Diag\!\left(
		\frac1N\sum_{j=1}^N
		\omega\bigl( \|\bx^i_{t-\tau}-\bx^j_{t-\tau} \|\bigr)
		\bigl( \bv^j_{t-\tau}- \bv^i_{t-\tau}\bigr)
		\right)\dd {\bf B}^i_t, \vspace{4pt}\\
		\displaystyle (\bx^i_s, \bv^i_s) = (\bphi^{i}_s, \bpsi^{i}_s),
	\quad s\in[-\tau,0],
	\end{cases}
\end{equation}
%where initial history $ (\bx^i_s, \bv^i_s)$ is deterministic and is continuous.
where the initial history $ (\bphi^i_s, \bpsi^i_s)$ is $\mathcal{F}_0
$-measurable and has continuous sample paths on $[-\tau,0]$. Next, we recall the concept of strong solution as follows.
\begin{definition}
Let $(\mathbf X,\mathbf V)$ be a strong solution to the Cauchy problem \eqref{B-1} in the time interval $[0, T)$ if and only if the following relations hold.
\begin{enumerate}
\item
For almost every \(\theta\in\Omega\), the solution path $t~\mapsto~(\bX_t(\theta), \bV_t(\theta)) \in {\mathbb R}^{2Nd}$ is continuous and $(\bX_t, \bV_t)$ is ${\mathcal F}_t$-adapted. 
\vspace{0.1cm}
\item
The stochastic processes $(\bX_t, \bV_t)$ satisfy $\eqref{B-1}_1 - \eqref{B-1}_2$ with probability one.
\vspace{0.1cm}
\item
The initial history conditions hold pointwise for all $s \in [-\tau, 0]$. 
\end{enumerate}
\end{definition}
\begin{theorem}[Global well-posedness]
	\label{T2.1}
Assume that time-delay and initial history data satisfy the following conditions:
\begin{enumerate}
\item
The time-delay is nonnegative and finite:
\[ 0 \leq \tau < \infty. \]
\item
The initial history data $(\bPhi,\bPsi) = \bigl((\bphi^1,\ldots, \bphi^N),(\bpsi^1,\ldots, \bpsi^N)\bigr)$ are \(\mathcal F_0\)-measurable and have continuous sample paths on
	\([-\tau,0]\), and satisfy
	\begin{align}\label{B-2}
	\E\left[
	\sup_{s\in[-\tau,0]}
	\left(
	\|\bPhi_s \|^2+ \|\bPsi_s \|^2
	\right)
	\right]
	<\infty.
	\end{align}
\end{enumerate}	
	Then the Cauchy problem \eqref{B-1} admits a
	unique global continuous adapted strong solution process $(\bX_t, \bV_t)$ in the time interval $[-\tau, \infty)$ satisfying the following bound:
	\begin{align}\label{B-3}
	\E\left[
	\sup_{-\tau\leq t\leq T}
	\left(
	\| \bX_t \|^2+ \| \bV_t \|^2
	\right)
	\right]
	<\infty, \quad \mbox{for every $T > 0$}.
	\end{align}
\end{theorem}
\begin{proof}
	When \(\tau=0\), the assertion follows from the standard
	well-posedness theory for SDEs with locally Lipschitz coefficients
	and linear growth. Hence, below we consider \(\tau>0\). The proof is similar to that of Theorem~2.1 in \cite{E-H-S}. One difference is that the initial data need not be constant. We therefore impose \eqref{B-2} to guarantee the well-posedness of system \eqref{B-1}. For suitable $\bF$ and $\bG$, we can rewrite 
the velocity equation $\eqref{B-1}_2$ as follows:
		\[
		\dd \bV_t
		=
		\bF\!\left(\bX_{t-\tau}, \bV_{t-\tau}
		\right)\dd t
		+
		\bG\!\left( \bX_{t-\tau}, \bV_{t-\tau}
		\right)\dd {\bf B}_t.
		\] 
		On the interval \(0\le t\le\tau\), the delayed processes
		are completely determined by the prescribed history, and hence we have
		\[
		\begin{cases}
		\displaystyle \bV_t
				=
				\bPsi_0
				+
				\int_0^t
				\bF\!\left(\bX_{s-\tau},
				\bV_{s-\tau}
				\right)\dd s  
				+
				\int_0^t
				\bG\!\left(
				\bX_{s-\tau},
				\bV_{s-\tau}
				\right) \dd {\bf B}_s, \vspace{6pt}\\
		\displaystyle \bX_t
		=
		\bPhi_0
		+
		\int_0^t \bV_s \dd s.
		\end{cases}
		\]
		Thus the solution is uniquely determined on
		\([-\tau,\tau]\). Then, by induction,
		the equations again become explicit stochastic integral equations,
		which uniquely determine
		\((\bX_t, \bV_t)\)
		on every
		\([k\tau,(k+1)\tau]\).
		Repeating this argument over successive intervals, we derive a unique
		adapted continuous solution on every finite time interval. Similar to Theorem 3.1 of \cite{M} and the subsequent remark on p.~157, we use linear growth of the vector field of system \eqref{B-1} to establish the desired second moment estimate \eqref{B-3}.
\end{proof}

Throughout the paper, we will use bar-notation above a quantity to denote empirical average. For example, we set the average momentum as follows:
\begin{align*}
\begin{aligned}
\bar{\bx}_t &:= \frac{1}{N} \sum_{i=1}^{N} \bx^i_t, \quad \bar{\bv}_t := \frac{1}{N} \sum_{i=1}^{N} \bv^i_t, \quad  t > -\tau, \\
\bar{\bphi}_t &:= \frac{1}{N} \sum_{i=1}^{N} \bphi^i_t, \quad \bar{\bpsi}_t := \frac{1}{N} \sum_{i=1}^{N} \bpsi^i_t, \quad  t \in [-\tau, 0].
\end{aligned}
\end{align*}
\begin{lemma} \label{L2.1}
Let $(\bX_t, \bV_t)$ be a strong solution process to \eqref{B-1}. Then, the expectation of the average velocity process is conserved:
\[
\E[\bar{\bv}_t] = \E[\bar{\bpsi}_0],  \quad t > 0.
\]
\end{lemma}
\begin{proof}
We sum $\eqref{B-1}_2$ over all $i \in [N]$, and then we divide the resulting relation by $N$ to get 
\begin{align}
\begin{aligned} \label{B-1-1}
\dd \Big(\frac{1}{N} \sum_{i=1}^{N}  \bv^i_t \Big)
		&=
		\frac{1}{N^2} \sum_{i,j=1}^N
		\omega\bigl( \|\bx^i_{t-\tau}-\bx^j_{t-\tau} \|\bigr)
		\bigl( \bv^j_{t-\tau}-\bv^i_{t-\tau}\bigr)\dd t \\
		&+ \frac{1}{N} \sum_{i=1}^{N} \sqrt{\sigma^i}\,
		\Diag\!\left(
		\frac1N\sum_{j=1}^N
		\omega\bigl( \|\bx^i_{t-\tau}-\bx^j_{t-\tau} \|\bigr)
		\bigl( \bv^j_{t-\tau}- \bv^i_{t-\tau}\bigr)
		\right)\dd {\bf B}^i_t \\
		&=: I_{11} + I_{12}.
\end{aligned}
\end{align}
Note that the first term $I_{11}$ in the right-hand side of \eqref{B-1-1} is identically zero using the index exchange symmetry $(i, j)~\longleftrightarrow~(j,i)$:
\begin{align*}
\begin{aligned} 
& \sum_{i,j=1}^N
		\omega\bigl( \|\bx^j_{t-\tau}-\bx^i_{t-\tau} \| \bigr)
		\bigl( \bv^j_{t-\tau}-\bv^i_{t-\tau}\bigr)   \\
& \hspace{0.5cm}  = \sum_{i,j=1}^N
		\omega\bigl( \|\bx^i_{t-\tau}-\bx^j_{t-\tau} \|\bigr)
		\bigl( \bv^i_{t-\tau}-\bv^j_{t-\tau}\bigr) = - 	 \sum_{i,j=1}^N
		\omega\bigl( \|\bx^j_{t-\tau}-\bx^i_{t-\tau} \|\bigr)
		\bigl( \bv^j_{t-\tau}-\bv^i_{t-\tau}\bigr).
\end{aligned}
\end{align*}
This yields
\begin{equation} \label{B-1-2}
I_{11} = \frac{1}{N^2} \sum_{i,j=1}^N
	 	\omega\bigl( \|\bx^i_{t-\tau}-\bx^j_{t-\tau} \|\bigr)
		\bigl( \bv^j_{t-\tau}-\bv^i_{t-\tau}\bigr) \dd t = 0.
\end{equation}
For notational simplicity, we set
\[ \bF^i_{t} := \frac{1}{N} \sum_{j=1}^N
		\omega\bigl( \|\bx^i_{t}-\bx^j_{t} \| \bigr)
		\bigl( \bv^j_{t}- \bv^i_{t})
\]
and use \eqref{B-1-1} and \eqref{B-1-2} to find 
\[
\dd \Big(\frac{1}{N} \sum_{i=1}^{N}  \bv^i_t \Big)
	= \frac1N\sum_{\ell=1}^N
	\sqrt{\sigma^\ell}\,
	\Diag\bigl( \bF^\ell_{t-\tau} \bigr)\dd {\bf B}^\ell_t.
\]
We integrate the above relation from $0$ to $t$ to get 
\[
\frac{1}{N} \sum_{i=1}^{N}  \bv^i_t = \frac{1}{N} \sum_{i=1}^{N}  \bpsi^i_0 +  \frac1N\sum_{\ell=1}^N \sqrt{\sigma^\ell}\ \int_0^t 
	\Diag\bigl(\bF^\ell_{s-\tau} \bigr)\dd {\bf B}^\ell_s.\]
By \eqref{B-3} and \(0<\omega\leq1\), the stochastic integral is
square-integrable and hence has zero expectation. Now, we take the expectation of the above relation to obtain the desired estimate:
\[ 
\E \Big[ \frac{1}{N} \sum_{i=1}^{N}  \bv^i_t  \Big] = \E \Big[  \frac{1}{N} \sum_{i=1}^{N}  \bpsi^i_0 \Big] +  \frac1N\sum_{\ell=1}^N  \sqrt{\sigma^\ell} \E \Big[   \int_0^t 
	\Diag\bigl(\bF^\ell_{s-\tau} \bigr)\dd {\bf B}^\ell_{s}       \Big] =  \E \Big[  \frac{1}{N} \sum_{i=1}^{N}  \bpsi^i_0 \Big].
\]
\end{proof}

\begin{lemma} \label{L2.2}
Let $(\bX_t, \bV_t)$ be a strong solution process to \eqref{B-1}. Then, we have
\begin{equation*}
	\frac1{N^2}\sum_{i,j=1}^N
	\| \bv^i_t-\bv^j_t \|^2
	=
	\frac2N\sum_{i=1}^N \| \bv^i_t- \bar{\bv}_t \|^2, \quad t \geq 0.
	\label{B-1-3}
\end{equation*}
\end{lemma}
\begin{proof} By direct calculation, one has 
\begin{align*}
\begin{aligned}
& \frac{1}{N^2} \sum_{i,j=1}^N
\| \bv^i_t-\bv^j_t \|^2 \\
&= \frac{1}{N^2} \sum_{i,j=1}^N \Big( \|\bv^i_t \|^2  +  \|  \bv^j_t \|^2  - 2 \bv^i_t \cdot \bv^j_t \Big) =\frac{2}{N}   \sum_{i=1}^N \|\bv^i_t \|^2 - 2 \Big (  \frac{1}{N} \sum_{i=1}^N \bv^i_t  \Big ) \cdot  \Big (   \frac{1}{N}  \sum_{j=1}^N \bv^j_t \Big ) \\
& = \frac{2}{N}   \sum_{i=1}^N \|\bv^i_t \|^2 - 2  \| \bar{\bv}_t \|^2 = \frac{2}{N} \Big( \sum_{i=1}^N \|\bv^i_t \|^2 -  N  \| \bar{\bv}_t \|^2 \Big) = \frac2N\sum_{i=1}^N \| \bv^i_t- \bar{\bv}_t \|^2.
\end{aligned}
\end{align*}
This yields the desired estimate.
\end{proof}
\begin{lemma}\label{L2.3}
	For \(p\geq2\), let \(\bx,\by\in\bbr^d\). Then, there exist positive
	constants \(c_p\) and \(C_p\), depending only on \(p\), such that
	\begin{align*}
		& (i)~
		\left(
		\norm{\bx}^{p-2}\bx-\norm{\by}^{p-2}\by
		\right)\cdot(\bx-\by)
		\geq
		c_p
		\left(
		\norm{\bx}+\norm{\by}
		\right)^{p-2}
		\norm{\bx-\by}^2,
		\\
		& (ii)~
		\left(
		\norm{\bx}^{p-2}\bx-\norm{\by}^{p-2}\by
		\right)\cdot(\bx-\by)
		\geq
		c_p\norm{\bx-\by}^p,
		\\
		& (iii)~
		\left\|
		\norm{\bx}^{p-2}\bx-\norm{\by}^{p-2}\by
		\right\|
		\leq
		C_p\left(
		\norm{\by}^{p-2}\norm{\bx-\by}
		+
		\norm{\bx-\by}^{p-1}
		\right).
	\end{align*}
\end{lemma}

\begin{proof}
	The proof is provided in Appendix~\ref{App-A-0}.
\end{proof}

\subsection{Previous result}\label{sec:2.2}
In this subsection, we briefly recall Theorem~3.1 of Erban, Haskovec, and Sun~\cite{E-H-S}. We set a handy notation for communication weight and  the associated unnormalized Laplacian $W_t$:
\[
\omega^{ij}_t := \omega(\|\bx^i_t-\bx^j_t \| ), \quad \forall~t \geq -\tau, \quad \forall~ i, j \in [N], \quad 
(W_t)^{ij}:=
\begin{cases}
-\omega^{ij}_t, \quad & i\neq j, \vspace{6pt} \\
\sum_{j\neq i} \omega^{ij}_t, \quad & i = j.
\end{cases}
\]
In the sequel, we denote $\lambda_2(W_t)$ by the second smallest eigenvalue of $W_t$.  Then, the Cauchy problem for the momentum equation in \cite{E-H-S} reads as follows:
\begin{equation} \label{B-3-1}
\begin{cases}
\displaystyle \dd \bv^i_t
    =-\frac 1 N
    \bigl(W_{t-\tau} \bV_{t-\tau} \bigr)^i\dd t
    -\frac{\sqrt{\sigma^i}}{N}
    \bigl(W_{t-\tau} \bV_{t-\tau} \bigr)^i\dd  { B}^i_t, \quad & t > 0,~~i \in [N],\vspace{6pt} \\
\displaystyle   \bv^i_t = \bpsi^i_t,
    \quad & t \in[-\tau,0],
 \end{cases}
\end{equation}
where $\{B_t^i\}_{i\in[N]}$ are independent one-dimensional standard Brownian motions.

\begin{theorem} \emph{\cite{E-H-S}} \label{T2.2}
Suppose that the system parameters satisfy 
\[  \sigma^{\infty}:=\max_{1\leq i\leq N}\sigma^i< 1,  \]
and let  $(\bX_t, \bV_t)$ be a strong solution process to \eqref{B-3-1} satisfying the following a priori assumption:~there exist positive constants $\ell$ and $L$ such that 
\[
    \lambda_2\bigl(W_t \bigr)\geq\ell>0,
    \qquad
    \norm{W_t-W_{t-\tau}}_{\mathrm F}
    \leq L\tau, \quad \mbox{a.s.}~~t \geq 0.
\]
Then, there exists $\tau_c=\tau_c(\sigma^\infty,L,\ell,1)>0$ independent of $N$ such that for $0\leq\tau<\tau_c$, 
\[
    \lim_{t\to\infty}
  \E\Big[ \| \bv^i_t- \bv^j_t \| \Big]=0,
    \quad i,j\in[N].
\]
\end{theorem}
\begin{remark} \label{R2.1}
Theorem~\ref{T2.2} cannot be applied directly to the
configuration-dependent CS communication \eqref{A-1}. Under stochastic forcing,
no deterministic uniform bound on the relative position distances is
available a priori. Along spreading configurations, the communication
between separated particles or clusters may therefore become arbitrarily
weak, so the uniform spectral-gap condition
\[
 \lambda_2(W_t )\geq\ell>0
\]
cannot be guaranteed uniformly in time along the realization. Verifying
this condition would already require the spatial cohesion that we seek
to verify. Our argument in this paper instead recovers velocity alignment through empirical spatial moment bounds and the method of time-varying effective region.  
\end{remark}

\subsection{Main result}\label{sec:2.3}
In this subsection, we delineate our main result on the quantitative stochastic flocking estimates. Next, we introduce a sufficient framework $({\mathcal F}_A)$ for stochastic flocking in terms of system parameters and initial data:
\vspace{0.1cm}
\begin{itemize}
\item
$({\mathcal F}_A1)$:~Parameters satisfy the following relations:
\[ 0 \leq \max_{i \in [N]} \sigma^i =  \sigma^{\infty} < 2,\quad 0\le\beta<1,\quad 0 \leq \tau < 1.
\]
\item
$({\mathcal F}_A2)$:~The communication weight function takes the following form:
\[
\omega(r)=\frac{1}{(1+r^2)^{\beta/2}},
	\quad r\geq0.
\]
\item
$({\mathcal F}_A3)$:~Initial history $(\bPhi_s, \bPsi_s)$ satisfies the following moment bounds:~for \[p\ge4,\quad p>  \frac{2}{1-\beta},\] there exist positive constants $M^0_{p,x},M^0_{p,v}<\infty$ independent of $N$ such that
\begin{align}
\begin{aligned} \label{B-4-1}
& (i)~\E\left[
	\sup_{s\in[-\tau,0]}
	\left(
	\|\bPhi_s \|^2+ \|\bPsi_s \|^2
	\right)
	\right]
	<\infty. \\
& (ii)~\E\left[
\sup_{s\in[-\tau,0]}
\frac1N\sum_{i=1}^N
\|\bphi^i_s- \bar{\bphi}_s \|^p
\right]
\leq M^0_{p,x}, \quad 
\E\left[
\sup_{s\in[-\tau,0]}
\frac1N\sum_{i=1}^N
\|\bpsi^i_s- \bar{\bpsi}_s \|^p
\right]
\leq M^0_{p,v}.
\end{aligned}
\end{align}
\end{itemize}
Next, we study the propagation of the $p$-th moment bounds for velocity and spatial fluctuations around the local averages. These higher-order moment bounds will be used in Section~\ref{sec:5}
to control the contribution of particles outside the effective region,
which cannot be obtained from the second moment estimate alone.
\begin{proposition} \label{P2.1}
Let $p\geq4$, and suppose that the initial history satisfies the
moment bounds in \eqref{B-4-1}. Let $(\bX_t, \bV_t)$ be a strong solution process to \eqref{B-1}. Then, there exist
$ \sigma_* = \sigma_{p,*}>0,~~\tau_* = \tau_{p,*}>0,$ 
such that if parameters satisfy 
\[
    0\leq\sigma^\infty\leq\sigma_*,
    \quad
    0\leq\tau\leq\tau_*,
\]
 the following estimates hold
	\begin{align*}
	\begin{aligned}
		& (i)~
		\sup_{t\geq-\tau}
		\E\left[
		\frac1N\sum_{i=1}^N
		\| \bv^i_t- {\bar \bv}_t \|^p
		\right]
		\leq
		C_p\bigl(1+M_{p,v}^0\bigr)
		<\infty. \\
		& (ii)~\sup_{s\in[-\tau,t]}
		\E\left[
		\frac1N\sum_{i=1}^N
		\norm{\bx^i_s- \bar{\bx}_s}^p
		\right]
		\leq
		C_p\bigl(
		1+M_{p,x}^0+M_{p,v}^0
		\bigr)(1+t)^p,
		\qquad t\geq0.
	\end{aligned}
\end{align*}
\end{proposition}
\begin{proof}
	Since the proof is very lengthy, we leave its rigorous proof in
	Section~\ref{sec:4}. However, for readers' convenience, we briefly
	outline a proof strategy as follows.\newline
	
	\noindent~$\bullet$~{\bf Step A}:~
	First, we derive
	\[
	\frac{\dd}{\dd t}\E\left[
	\frac1N\sum_{i=1}^N
	\norm{\bv_t^i-\bar{\bv}_0}^p
	\right]
	\leq
	-\frac{p}{2}{\mathcal D}_p(t-\tau)
	+
	B_p
	\int_{t-3\tau}^{t-\tau}
	{\mathcal D}_p(s)\,\dd s,
	\quad
	\mbox{a.e. }t\geq3\tau,
	\]
	where
	\[
	{\mathcal D}_p(t)
	:=
	\E\left[
	\frac1{N^2}\sum_{i,j=1}^N
	\omega_t^{ij}
	\norm{\bv_t^i-\bar{\bv}_0}^{p-2}
	(\bv_t^i-\bar{\bv}_0)
	\cdot
	(\bv_t^i-\bv_t^j)
	\right].
	\]
	\vspace{0.2cm}
	
	\noindent~$\bullet$~{\bf Step B}:~
	Second, to deal with the integral term in the above estimate, we introduce
	\[
	\mathcal L_p(t)
	:=
\E\left[
\frac1N\sum_{i=1}^N
\norm{\bv_t^i-\bar{\bv}_0}^p
\right]
	+
	B_p
	\int_{t-3\tau}^{t-\tau}
	(s-t+3\tau){\mathcal D}_p(s)\,\dd s.
	\]
	Then, for sufficiently small \(\tau\) and \(\sigma^\infty\), we obtain
	\[
	\frac{\di}{\di t}\mathcal L_p(t)
	\leq
	-\frac{p}{4}{\mathcal D}_p(t-\tau),
	\quad
	\mbox{a.e. }t\geq3\tau.
	\]
	Together with the finite-time estimate on \([-\tau,3\tau]\) and
	\[
	\frac1N\sum_{i=1}^N
	\norm{\bv_t^i-\bar{\bv}_t}^p
	\leq
	2^p
	\frac1N\sum_{i=1}^N
	\norm{\bv_t^i-\bar{\bv}_0}^p,
	\]
	this implies
	\[
	\sup_{-\tau\leq s<\infty}
	\E\left[
	\frac1N\sum_{i=1}^N
	\norm{\bv_s^i-\bar{\bv}_s}^p
	\right]\le C_p\left(1+M_{p,v}^0\right).
	\]
	
	\vspace{0.2cm}
	
	\noindent~$\bullet$~{\bf Step C}:~
	Finally, we use
	\[
	\bx_t^i-\bar{\bx}_t
	=
	\bx_0^i-\bar{\bx}_0
	+
	\int_0^t
	(\bv_s^i-\bar{\bv}_s)\,\dd s
	\]
	and the estimate in Step B to obtain
	\[
	\sup_{-\tau\leq s\leq t}
	\E\left[
	\frac1N\sum_{i=1}^N
	\norm{\bx_s^i-\bar{\bx}_s}^p
	\right]
	C_p\left(1+M_{p,x}^0+M_{p,v}^0\right)(1+t)^p.
	\]
	For details, we refer to Section~\ref{sec:4}.
\end{proof}
Now, we are ready to provide the first main result. 
\begin{theorem} \label{T2.3}
Suppose the framework $({\mathcal F}_A1) - ({\mathcal F}_A3)$ holds, and let  $(\bX_t, \bV_t)$ be a strong solution process to \eqref{B-1}. Then, there exist
$ \sigma_* = \sigma_*(p)>0,~~\tau_* = \tau_*(p)>0,$ 
such that the following assertions hold whenever
\[
    0\leq\sigma^\infty\leq\sigma_*,
    \quad
    0\leq\tau\leq\tau_*.
\]
For \(0<\beta<1\), let
\[
0<\gamma<
\frac{(p-2)(1-\beta)}{\beta},
\]
and when \(\beta=0\), let \(\gamma>0\) be arbitrary.
Then, there exist \(C_\gamma>0\) and $C_x>0$ independent of \(N\) such that
\begin{align}
	\begin{aligned}\label{B-7}
		& \E\left[
		\frac1N\sum_{i=1}^N
		\|\bv_t^i-\bar{\bv}_t\|^2
		\right]
		\leq
		C_\gamma(1+t)^{-\gamma},
		~~ t\geq0,\quad \sup_{0\leq t<\infty}
		\E\left[
		\frac1N\sum_{i=1}^N
		\|\bx_t^i-\bar{\bx}_t\|^2
		\right]
		\leq
		C_x.
	\end{aligned}
\end{align}
\end{theorem}
\begin{proof}
	Since the proof is very lengthy, we leave its rigorous proof in
	Sections~\ref{sec:3}--\ref{sec:6}. However, for readers' convenience,
	we briefly outline a proof strategy as follows.\newline
	
	\noindent~$\bullet$~{\bf Step A}:~
	First, to deal with the memory term in the time-evolution estimate of
	\[
	\E\left[
	\frac1N\sum_{i=1}^N
	\norm{\bv_t^i-\bar{\bv}_t}^2
	\right]
	\] 
	we introduce the compensated energy functional
	\[
	\mathcal L(t)
	:=
		\E\left[
	\frac1N\sum_{i=1}^N
	\norm{\bv_t^i-\bar{\bv}_t}^2
	\right]
	+
	B\int_{t-2\tau}^{t-\tau}
	(s-t+2\tau)\Lambda^\omega(s)\,\dd s,
	\qquad t\geq\tau,
	\]
	where
	\[
	\Lambda^\omega(t)
	:=
	\E\left[
	\frac1{N^2}\sum_{i,j=1}^N
	\omega_t^{ij}
	\norm{\bv_t^i-\bv_t^j}^2
	\right].
	\]
	Then, under suitable smallness conditions on \(\tau\) and
	\(\sigma^\infty\), we obtain
	\[
	\frac{\di}{\di t}\mathcal L(t)
	\le-\frac{2-\sigma^{\infty}}{12}
	\Lambda^\omega(t-\tau),
	\qquad
	\mbox{a.e. }t\geq\tau.
	\]
	See Section~\ref{sec:3}.
	
	\vspace{0.2cm}
	
	\noindent~$\bullet$~{\bf Step B}:~
	Second, we obtain the following velocity and spatial \(p\)-th ($p\ge4$) moment estimates:
	\[
	\sup_{-\tau\leq t<\infty}
	\E\left[
	\frac1N\sum_{i=1}^N
	\norm{\bv_t^i-\bar{\bv}_t}^p
	\right]
	<\infty,\quad
	\sup_{-\tau\leq s\leq t}
	\E\left[
	\frac1N\sum_{i=1}^N
	\norm{\bx_s^i-\bar{\bx}_s}^p
	\right]
	\lesssim
	(1+t)^p,
	\quad t\geq0.
	\]
	These \(p\)-th moment estimates will be used in Step C to control the
	tails outside the effective region, whereas the flocking quantities
	in \eqref{B-7} are second-order moments. See Section~\ref{sec:4}.
	
	\vspace{0.2cm}
	
	\noindent~$\bullet$~{\bf Step C}:~
	Now, we are ready to use the method of time-varying effective region to
	estimate \(\Lambda^\omega(t)\). For the gauge function
	\[
	R(t):=(1+t)^\alpha,
	\]
	we derive
	\[
	\Lambda^\omega(t)
	\geq
	2^{-2\beta}(1+t)^{-\beta\alpha}\E\left[
	\frac1N\sum_{i=1}^N
	\norm{\bv_t^i-\bar{\bv}_t}^2
	\right]
	-
	2^{-\beta}\tilde{C}_p(1+t)^{-\beta\alpha-(p-2)(\alpha-1)},
	\qquad t\geq0.
	\]
	See Section~\ref{sec:5}.
	
	\vspace{0.2cm}
	
	\noindent~$\bullet$~{\bf Step D}:~
	We apply the estimate in Step A at time $t+\tau$ and combine it with the estimate in Step C to obtain
	\[
	\begin{aligned}
		\frac{\dd}{\dd t}\mathcal L(t+\tau)
		&\leq
		-\frac{2-\sigma^\infty}{24}\Lambda^\omega(t)
		-\frac{2-\sigma^\infty}{24}\,2^{-2\beta}
		(1+t)^{-\beta\alpha}
		\E\left[
		\frac1N\sum_{i=1}^N
		\|\bv_t^i-\bar{\bv}_t\|^2
		\right]
		\\
		&
		+\frac{2-\sigma^\infty}{24}\,2^{-\beta}\tilde C_p
		(1+t)^{-\beta\alpha-(p-2)(\alpha-1)}.
	\end{aligned}
	\]
%	We apply the estimate in Step A at time \(t+\tau\) and combine it with the estimate in Step C to obtain
%\[
%\frac{\di}{\di t}\mathcal L(t+\tau)
%\leq
%-2^{-\beta}\Lambda^\omega(t)
%-2^{-2\beta}(1+t)^{-\beta\alpha}
%\E\left[
%\frac1N\sum_{i=1}^N
%\norm{\bv_t^i-\bar{\bv}_t}^2
%\right]
%+
%2^{-\beta}\tilde{C}_p(1+t)^{-\beta\alpha-(p-2)(\alpha-1)}.
%\]
	For \(0<\beta<1\), for the above \(\gamma\), we choose
	\(\alpha\) such that
	\[
	1<\alpha<\frac1\beta,
	\qquad
	\gamma<
	\beta\alpha+(p-2)(\alpha-1)-1.
	\]  
	For \(\beta=0\), we choose \(\alpha>1\) sufficiently large.
	Then, we obtain
	\[
	(1+t)^\gamma\mathcal L(t+\tau)
	\leq
	C_\gamma.
	\]
	Together with
	\[
	\E\left[
	\frac1N\sum_{i=1}^N
	\norm{\bv_{t+\tau}^i-\bar{\bv}_{t+\tau}}^2
	\right]\leq\mathcal L(t+\tau),
	\]
	we have
	\[
	\E\left[
	\frac1N\sum_{i=1}^N
	\norm{\bv_t^i-\bar{\bv}_t}^2
	\right]
	\leq
	C_\gamma(1+t)^{-\gamma},
	\qquad t\geq0,
	\]
	which gives the desired velocity alignment estimate
	\(\eqref{B-7}_1\). See Section~\ref{sec:6.1}.
	
	\vspace{0.2cm}
	
	\noindent~$\bullet$~{\bf Step E}:~
Finally, when \(0<\beta<1\), the condition
\[
p>\frac{2}{1-\beta}
\]
allows us to choose any
\[
2<\gamma_0<
\frac{(p-2)(1-\beta)}{\beta}.
\]
When \(\beta=0\), we choose any \(\gamma_0>2\).
Then, we have
\[
\E\left[
\frac1N\sum_{i=1}^N
\norm{\bv_t^i-\bar{\bv}_t}^2
\right]
\leq
C_{\gamma_0}(1+t)^{-\gamma_0},
\qquad t\geq0.
\]
We use this estimate and the relation
	\[
	\begin{aligned}
		&\left(
		\E\left[
		\frac1N\sum_{i=1}^N
		\norm{\bx_t^i-\bar{\bx}_t}^2
		\right]
		\right)^{1/2}\\
	&
	\hspace{1cm}\leq
		\left(
		\E\left[
		\frac1N\sum_{i=1}^N
		\norm{\bx_0^i-\bar{\bx}_0}^2
		\right]
		\right)^{1/2}
		+
		\int_0^t
		\left(\E\left[
		\frac1N\sum_{i=1}^N
		\norm{\bv_s^i-\bar{\bv}_s}^2
		\right]\right)^{1/2}\,\dd s
	\end{aligned}
	\]
	to derive the desired spatial cohesion estimate
	\(\eqref{B-7}_2\). See Section~\ref{sec:6.2}.
\end{proof}

\begin{remark} \label{R2.2}
In the sequel, we comment on the framework and the results of Theorem \ref{T2.3} as follows.
		\begin{enumerate}
			\item
		Our conditions on initial data are parallel to those for the kinetic CS model \cite{H-W} with non-compact support and stochastic CS model \cite{H-Z}. That is, we require the high moments of velocity and position to be finite. This is important when we study the large-time behavior of system \eqref{A-2}. However, our velocity moment is not monotonic compared to that in literature \cite{H-W,H-Z}. This is why we need very precise estimation.
		
			\vspace{0.2cm}
			
			\item
		The results in Theorem \ref{T2.3} demonstrate that as long as the
		initial velocity and spatial distributions satisfy the prescribed
		moment conditions, the corresponding moment control will persist
		along system \eqref{A-2}.
		
			\vspace{0.2cm}
			
			\item
			Compared to the previous work \cite{C-S1,D-Z,E-H-S}, our results extend their findings by relaxing the assumption of communication weight function to allow for original CS weight \eqref{A-1}. This provides a positive answer to the open problem raised in \cite{E-H-S}. In particular, we show that a  polynomial moment yields a polynomial convergence rate.
			
			\vspace{0.2cm}
			
			\item Our result can be extended to a general locally Lipschitz communication weight
			satisfying
			\[
			\frac{c_{\omega}}{(1+r)^{\beta}}\leq\omega(r)\leq1,
			\]
			since our proof only relies on the decay rate of the communication weight
			function as \(r\to\infty\).
			
		\vspace{0.2cm}
		
				\item The numerical experiments in Section~\ref{sec:8} illustrate that
				large time-delay or noise intensity may deteriorate
				flocking behavior.
				
				\vspace{0.2cm}
				
			\item Our results cannot be extended to the critical exponent $\beta=1$, since the spatial moment has polynomial growth and our gauge function cannot provide sufficient decay for the mass outside the
			effective domain when \(\beta=1\). Therefore, we leave this critical case as a remaining open problem.
		\end{enumerate}
	\end{remark}
	
\begin{theorem} \label{T2.4}
Suppose the framework $({\mathcal F}_A1) - ({\mathcal F}_A3)$ holds, and let  $(\bX_t, \bV_t)$ be a strong solution process to \eqref{B-1}. Then, there exist
$ \sigma_* = \sigma_*(p)>0,~~\tau_* = \tau_*(p)>0,$ 
such that if parameters satisfy 
\begin{equation*}
    0\leq\sigma^\infty\leq\sigma_*,
    \quad
    0\leq\tau\leq\tau_*,
\end{equation*}
for every fixed \(N\), we have almost-sure stochastic flocking:
\[ 
	\lim_{t\to\infty}
	\max_{1\leq i,j\leq N} \|\bv^i_t-\bv^j_t \|=0,\quad
	\sup_{0 \leq t < \infty}
	\max_{1\leq i,j\leq N} \|\bx^i_t-\bx^j_t \| <\infty,\quad\text{a.s.}
\]	
\end{theorem}
\begin{proof}
Since the proof is very lengthy, we leave its proof in Section \ref{sec:7}.
\end{proof}	
\begin{remark}
	For the deterministic delayed Cucker--Smale system \cite{H}, the
	average velocity is conserved. In contrast, from Lemma \ref{L2.1}, we see that the
	stochastic perturbation in \eqref{B-1} destroys this pathwise
	conservation law, although the expectation of the average velocity is preserved. Thus, Theorem~\ref{T2.4} should be understood as a pathwise
	velocity alignment result rather than convergence toward the
	deterministic initial average velocity.
\end{remark}

\vspace{0.2cm}

\section{Compensated energy and dissipation} \label{sec:3}
\setcounter{equation}{0}
In this section, we construct a compensated energy functional and derive
a dissipation estimate for the velocity fluctuations. As a consequence,
we obtain the time-integrability of the alignment production functional,
which will be used in the subsequent flocking analysis.
\subsection{Preparatory estimates} \label{sec:3.1}
In this subsection, we study several elementary estimates to be used later.  Let  $(\bX_t, \bV_t)$ be a global strong solution process to \eqref{B-1}. Then, we define energy and its production functionals as follows.
 \begin{align}
 \begin{aligned} \label{C-1}
 & \mathcal E(t):=\E\left[
\frac1N\sum_{i=1}^N \|\bv^i_t-\bar{\bv}_t \|^2
\right], \quad \Lambda^{\omega}(t) :=\E\left[
	\frac1{N^2}\sum_{i,j=1}^N
	\omega^{ij}_{t}  \|\bv^{i}_{t} - \bv^{j}_{t} \|^2
	\right],\\
&  \bF^i_{t} := \frac{1}{N} \sum_{j=1}^N
		\omega^{ij}_t \bigl( \bv^j_{t}- \bv^i_{t}), \quad \Lambda^{\sigma}(t)
    := \E\left[
        \frac1N\sum_{i=1}^N
        \sigma^i \| \bF^i_{t} \|^2
    \right], \\
&   \Lambda^{c}(t,t-\tau) 
    := \E\left[
        \frac1N\sum_{i=1}^N
        \left(
            \bigl(\bv^i_t-\bar{\bv}_t\bigr)
            -
            \bigl(\bv^i_{t-\tau} -\bar{\bv}_{t-\tau} \bigr)
        \right)\cdot  \bF^i_{t-\tau}
    \right].
 \end{aligned}
 \end{align}   
 Here ${\mathcal E}$ measures velocity fluctuations around the time-dependent velocity average, whereas $\Lambda^\omega,~\Lambda^{\sigma}$ and  $\Lambda^{c}$ denote the production functionals for velocity alignment, stochastic noise, and crossing term, respectively.

 \begin{lemma} \label{L3.1}
Let $(\bX_t, \bV_t)$ be a strong solution process to \eqref{B-1}. Then, one has 
\[ 
\Lambda^{\omega}(t) \leq 2{\mathcal E}(t), \quad  \mbox{for}~t \geq 0.
\]
\end{lemma}
\begin{proof}
Since $0 < \omega^{ij}_t \leq1$, we use Lemma \ref{L2.2} to see that 
 \[
 \Lambda^{\omega}(t) \leq  \E\left[\frac1{N^2} 
	\sum_{i,j=1}^N  \|\bv^{i}_{t} - \bv^{j}_{t} \|^2 
	\right]  = 2\E \left[ \frac{1}{N} \sum_{i=1}^N \| \bv^i_{t}- \bar{\bv}_{t} \|^2 \right] = 2{\mathcal E}(t).
\]
\end{proof}
Next, we study the temporal-evolution of ${\mathcal E}(t)$.
\begin{lemma} \label{L3.2}
Suppose the framework $({\mathcal F}_A1) - ({\mathcal F}_A3)$ holds, and let  $(\bX_t, \bV_t)$ be a strong solution process to \eqref{B-1}. Then, the energy functional ${\mathcal E}$ satisfies the following properties:
\begin{enumerate}
\item[(i)]
${\mathcal E}(t)$ is locally absolutely continuous on $[0,\infty)$.
\vspace{0.2cm}
\item[(ii)]
For a.e.~$t > \tau$, one has 
\begin{equation} \label{C-1-1}
 \frac{\di{\mathcal E}(t)}{\di t} 
    =
    - {\Lambda}^{\omega}(t-\tau)
    + 2 \Lambda^{c}(t,t-\tau)
    +  \left(1-\frac1N\right) \Lambda^{\sigma}(t-\tau).
 \end{equation}
 \end{enumerate}
\end{lemma}
\begin{proof}
\noindent (i)~For $ i \in [N]$, we set 
\[
    \bDelta^{\bv, i}_t:=\bv^i_t-\bar{\bv}_t, \quad   \bDelta^{\bv}_t := ( \bDelta^{\bv, 1}_t, \cdots,  \bDelta^{\bv, N}_t),
    \quad  \|\bDelta^{\bv}_t \|_N
    := \Big( \frac1N\sum_{i=1}^N  \|\bDelta^{\bv,i}_t \|^2 \Big)^{\frac{1}{2}}.
\]
We first perform the time-evolution of ${\mathcal E}(t)$ under a stopping time. For
\(R>0\), we define a stopping time $T^R$:
\[
    T^R
    :=
    \inf\left\{
        t \in [0, \infty):~
        \sup_{s\in[-\tau,t]}
        \bigl( \| \bX_s \|  + \| \bV_s \| \bigr)
        \geq R
    \right\}.
\]
By \eqref{B-1-1} and \eqref{B-1-2}, we have
\begin{equation}\label{C-1-2}
	\begin{aligned}
		\dd\bigl(\bv^i_t-\bar{\bv}_t \bigr)
		&=
		\bF^i_{t-\tau}\,\dd t
		+
		\sqrt{\sigma^i}\,
		\Diag(\bF^i_{t-\tau})\,\dd {\bf B}^i_t
		-\frac1N\sum_{\ell=1}^N
		\sqrt{\sigma^\ell}\,
		\Diag(\bF^\ell_{t-\tau})\,\dd {\bf B}^\ell_t, \quad t\ge0.
	\end{aligned}
\end{equation}
We apply It\^o's formula to \eqref{C-1-2} up to \(t\wedge T^R\) to get
\begin{align}
\begin{aligned} \label{C-2}
\| \bDelta^{\bv}_{t\wedge T^R} \|_N^2
    &= \|\bDelta^{\bv}_0 \|_N^2
    + \frac{2}{N} \sum_{i=1}^N \int_0^t
    \one_{\{s\leq T^R\}}
    \bDelta^{\bv,i}_s \cdot \bF^i_{s-\tau}\,\dd s\\
    &
    +  \left(1-\frac1N\right)
    \frac1N\sum_{i=1}^N
    \int_0^t
    \one_{\{s\leq T^R\}}
    \sigma^i  \| \bF^i_{s-\tau} \|^2\,\dd s
    + M^R_{t},
\end{aligned}
\end{align}
where \(M^R_{t}\) is a square-integrable martingale with the form:
\[
\begin{aligned}
M^R_{t}
:=
\frac{2}{N}
\sum_{i,\ell=1}^N
\int_0^t
\one_{\{s\leq T^R\}}
\bDelta^{\bv,i}_{s}\cdot\left[
\delta_{i\ell}\sqrt{\sigma^i}\Diag(\bF^i_{s-\tau})
-\frac1N\sqrt{\sigma^\ell}\Diag(\bF^\ell_{s-\tau})
\right]
\dd \bf{B}^\ell_s.
\end{aligned}
\]
Here we use
\begin{equation}\label{C-3}
  \frac1N\sum_{i,\ell=1}^N
    \left\|
        \delta_{i\ell}\sqrt{\sigma^i}\,\Diag(\bF^i_{t-\tau})
        -
        \frac1N\sqrt{\sigma^\ell}\,\Diag(\bF^\ell_{t-\tau})
    \right\|_{\mathrm{F}}^2 =
    \left(1-\frac1N\right)
    \frac1N\sum_{\ell=1}^N
    \sigma^\ell  \|\bF^\ell_{t-\tau} \|^2.
\end{equation}
By \eqref{B-3} and \(0<\omega\leq1\), for every \(T>0\), it holds
\[
\E\left[
\sup_{0\leq s\leq T}\|\bDelta^{\bv}_s\|_N^2
\right]
+
\E\left[
\int_0^T
\frac1N\sum_{i=1}^N
\|\bF^i_{s-\tau}\|^2\,\dd s
\right]
<\infty.
\]
In particular, we have \(T^R\to\infty\) a.s. as \(R\to\infty\).
Since \(M_t^R\) is a square-integrable martingale,
\(\E[M_t^R]=0\). Taking expectations in \eqref{C-2} and then
letting \(R\to\infty\), the dominated convergence theorem yields
\[
\begin{aligned}
    \mathcal E(t)
    =
    \mathcal E(0)
    +
    2\int_0^t
    \E\left[
        \frac1N\sum_{i=1}^N
        \bDelta^{\bv,i}_s \cdot \bF^i_{s-\tau}
    \right]\dd s+ \Big(1 - \frac{1}{N} \Big)
    \int_0^t \Lambda^{\sigma}(s-\tau) \dd s,\quad t\ge0.
\end{aligned}
\]
This implies the local absolute continuity of \(\mathcal E\) on $[0,\infty)$. \newline

\noindent (ii)~For almost every \(t\geq0\), one has 
\begin{equation} \label{C-3-0}
 \frac{\di \mathcal E(t)}{\di t}
    =
    2\E\left[
        \frac1N\sum_{i=1}^N
        \bDelta^{\bv,i}_t\cdot \bF^i_{t-\tau}
    \right]
    +\Big(1 - \frac{1}{N} \Big)
    \Lambda^{\sigma}(t-\tau).
\end{equation}
It remains to separate the delayed dissipation from the defect term. Since the communication weight is invariant under the index exchange transformation $(i, j)~\longleftrightarrow~(j,i)$:
\[ \omega^{ij}_t=\omega^{ji}_t, \quad i, j \in [N], \quad t > 0, \]
we have
\begin{equation} \label{C-3-1}
   \frac1N\sum_{i=1}^N
    \bDelta^{\bv,i}_{t-\tau} \cdot \bF^i_{t-\tau}=-
    \frac1{2N^2}\sum_{i,j=1}^N \omega^{ij}_{t-\tau}
     \| \bv^{i}_{t-\tau}-  \bv^{j}_{t-\tau} \|^2.
\end{equation}
Then, we use \eqref{C-3-1} to get 
\begin{align}
\begin{aligned} \label{C-3-2}
\frac1N\sum_{i=1}^N
    \bDelta^{\bv,i}_t\cdot \bF^i_{t-\tau} &=\frac1N\sum_{i=1}^N\left[ \bDelta^{\bv,i}_{t-\tau}
    +
    \bigl(\bDelta^{\bv,i}_t- \bDelta^{\bv,i}_{t-\tau} \bigr)\right]\cdot \bF^i_{t-\tau} \\
&=-\frac1{2N^2}\sum_{i,j=1}^N \omega^{ij}_{t-\tau} \|\bv^{i}_{t-\tau}- \bv^{j}_{t-\tau} \|^2 +
    \frac1N\sum_{i=1}^N
    \bigl(\bDelta^{\bv,i}_t- \bDelta^{\bv,i}_{t-\tau} \bigr)\cdot \bF^i_{t-\tau}.
\end{aligned}
\end{align}
We take the expectations of \eqref{C-3-2} to find 
\begin{align}
\begin{aligned} \label{C-3-3}
& \E\left[ \frac1N\sum_{i=1}^N
    \bDelta^{\bv,i}_t\cdot \bF^i_{t-\tau}     \right]  \\
 & \hspace{1cm} = - \frac{1}{2} \E \left[ \frac1{N^2}\sum_{i,j=1}^N \omega^{ij}_{t-\tau} \| \bv^{i}_{t-\tau}- \bv^{j}_{t-\tau} \|^2               \right] + \E \left[   \frac1N\sum_{i=1}^N
    \bigl(\bDelta^{\bv,i}_t- \bDelta^{\bv,i}_{t-\tau} \bigr)\cdot \bF^i_{t-\tau}      \right].
\end{aligned}
\end{align}
Finally, we combine \eqref{C-1}, \eqref{C-3-0} and \eqref{C-3-3} to obtain the desired estimate \eqref{C-1-1} for almost every \(t\geq\tau\).
\end{proof}
Next, we estimate $\Lambda^{\sigma}(t-\tau)$ in terms of ${\Lambda}^\omega(t-\tau)$.
\begin{lemma}  \label{L3.3}
The stochastic noise production functional is dominated by the velocity alignment production functional:
\[
    0\leq \Lambda^{\sigma}(t) \leq \frac{\sigma^\infty}{2} \Lambda^{\omega}(t), \quad t \geq -\tau.
\]
\end{lemma}
\begin{proof} We use \(\sigma^i\leq\sigma^{\infty} \) and the definition of \(\Lambda^\sigma(t)\) to see that 
\begin{equation} \label{C-3-4}
    \Lambda^{\sigma}(t)
    =
    \E\left[
        \frac1N\sum_{i=1}^N
        \sigma^i  \| \bF^i_{t} \|^2
    \right]\leq
    \sigma^{\infty}
    \E\left[
        \frac1N\sum_{i=1}^N
        \norm{\bF^i_{t}}^2
    \right].
\end{equation}
Now, we claim that 
\begin{equation} \label{C-3-5}
  \E\left[
        \frac1N\sum_{i=1}^N
        \norm{\bF^i_{t}}^2
    \right] \leq \frac{1}{2} \Lambda^{\omega}(t).
\end{equation} 
{\it Proof of \eqref{C-3-5}}:~We use \(\omega_t^{ij}=\omega_t^{ji}\) to see
\begin{equation} \label{C-4-1}
 \sum_{i=1}^N  \bF^i_{t} =
    \frac{1}{2N}\sum_{i,j=1}^N
    (\omega^{ij}_{t} -\omega^{ji}_{t})
    \bigl(
         \bv^j_{t}- \bv^i_{t}
    \bigr)
    = {\bf 0}.
\end{equation}
Then, we use \eqref{C-4-1},  Young's inequality and \(0\leq\omega^{ij}_t \leq1\)  to get 
\begin{align*}
\begin{aligned}
    \sum_{i=1}^N \| \bF^i_{t} \|^2
    &=
    \frac1N\sum_{i,j=1}^N
    \omega^{ij}_{t}
    \bigl(
         \bv^j_{t}- \bv^i_{t}
    \bigr)\cdot \bF^i_{t} =
    \frac1{2N}\sum_{i,j=1}^N
    \omega^{ij}_{t}
    \bigl(\bF^i_{t}-\bF^j_{t}\bigr)\cdot
    \bigl(
         \bv^j_{t}- \bv^i_{t}
    \bigr)\\
    &\leq
    \frac1{4N}\sum_{i,j=1}^N
    \omega^{ij}_{t} \| \bF^i_{t}-\bF^j_{t} \|^2 +
    \frac1{4N}\sum_{i,j=1}^N
    \omega^{ij}_{t}
    \| \bv^j_{t}- \bv^i_{t} \|^2\\
    &\leq
    \frac1{4N}\sum_{i,j=1}^N
    \| \bF^i_{t}-\bF^j_{t} \|^2+
    \frac1{4N}\sum_{i,j=1}^N
    \omega^{ij}_{t}
    \| \bv^j_{t}- \bv^i_{t} \|^2\\
    &\leq
\frac12\sum_{i=1}^N\norm{\bF^i_{t}}^2+
    \frac1{4N}\sum_{i,j=1}^N
    \omega^{ij}_{t}
    \| \bv^j_{t}- \bv^i_{t} \|^2.
\end{aligned}
\end{align*}
This yields the desired estimate:
\begin{equation}\label{C-4-2}
	\frac1N\sum_{i=1}^N \|\bF_t^i\|^2
	\leq
	\frac1{2N^2}\sum_{i,j=1}^N
	\omega_t^{ij}
	\|\bv_t^i-\bv_t^j\|^2,
	\qquad\mbox{a.s.}
\end{equation}
Taking expectations in \eqref{C-4-2}, we obtain \eqref{C-3-5}.
Finally, we combine \eqref{C-3-4} and \eqref{C-3-5}
to get the desired estimate.
\end{proof}
In the following lemma, we estimate $\Lambda^{c}(t, t-\tau)$ in terms of $\Lambda^\omega(t-\tau)$ and its integral. 
\begin{lemma} \label{L3.4}
For $t \geq \tau$ and $\varepsilon > 0$, the crossing production functional can be controlled by the velocity alignment production functional and its integral. 
\[ \left| \Lambda^{c}(t, t-\tau) \right|
    \leq  \varepsilon \Lambda^{\omega}(t-\tau)   +
    \frac{(\tau+\sigma^{\infty})}{2\varepsilon}
    \int_{t-2\tau}^{t-\tau}\Lambda^{\omega}(s) \,\dd s, \quad t \geq \tau.
\]
\end{lemma}
\begin{proof} We use the definition of the crossing interaction production functional and the Cauchy--Schwarz inequality, and \eqref{C-3-5} to get 
\begin{align}
\begin{aligned}  \label{C-5-0}
 & \left| \Lambda^{c}(t,t-\tau)  \right| \\
 & \hspace{0.5cm} \leq \E\left[
        \frac1N\sum_{i=1}^N
         \|
            \big (\bv^i_t-\bar{\bv}_t \big )
            -
            \big (\bv^i_{t-\tau} -\bar{\bv}_{t-\tau} \big) \|
      \| \bF^i_{t-\tau} \|
    \right] \\
    & \hspace{0.5cm} \le  
    \left(
        \E\left[
            \frac1N\sum_{i=1}^N
            \norm{
                \bigl(\bv^i_t- \bar{\bv}_t \bigr)
                -
                \bigl( \bv^i_{t-\tau}-\bar{\bv}_{t-\tau} \bigr)
            }^2
        \right]
    \right)^{1/2}
    \left(
        \E\left[
            \frac1N\sum_{i=1}^N\norm{\bF^i_{t-\tau}}^2
        \right]
    \right)^{1/2}\\
    & \hspace{0.5cm} \leq  \left(
        \E\left[
            \frac1N\sum_{i=1}^N
            \norm{
                \bigl(\bv^i_t- \bar{\bv}_t \bigr)
                -
                \bigl( \bv^i_{t-\tau}-\bar{\bv}_{t-\tau} \bigr)
            }^2
        \right]
    \right)^{1/2} \Big (  \frac{1}{2} \Lambda^{\omega}(t-\tau) \Big )^{1/2}. 
\end{aligned}
\end{align}
Below, we estimate the first factor on the right-hand side of \eqref{C-5-0}. More precisely, we claim that for $t \geq \tau$, 
\begin{equation} \label{C-5}
~\E\left[
        \frac1N\sum_{i=1}^N
        \|
            \bigl(\bv^i_t- \bar{\bv}_t \bigr)
            -
            \bigl( \bv^i_{t-\tau}-\bar{\bv}_{t-\tau} \bigr)
        \|^2
    \right]\leq
    2(\tau+\sigma^\infty)
    \int_{t-2\tau}^{t-\tau} \Lambda^{\omega}(s) \dd s.
 \end{equation}   
{\it Proof of \eqref{C-5}}:~We integrate \eqref{C-1-2} over \([t-\tau,t]\) to obtain
\begin{align*}
\begin{aligned}
 &\bigl(\bv^i_t- \bar{\bv}_t \bigr)
    -  \bigl( \bv^i_{t-\tau}-\bar{\bv}_{t-\tau}\bigr)\\
    &\quad=
    \int_{t-\tau}^{t} \bF^i_{s-\tau}\,\dd s
    +\sum_{\ell=1}^N
    \int_{t-\tau}^{t}
    \left[
        \delta_{i\ell}\sqrt{\sigma^i}\Diag( \bF^i_{s-\tau})
        -
        \frac1N\sqrt{\sigma^\ell}\Diag(\bF^\ell_{s-\tau})
    \right]\dd {\bf B}^\ell_s.
\end{aligned}
\end{align*}
This implies
\begin{align}
\begin{aligned} \label{C-5-1}
& \E\left[ \frac1N\sum_{i=1}^N \| \bigl(\bv^i_t- \bar{\bv}_t \bigr)
    -  \bigl( \bv^i_{t-\tau}-\bar{\bv}_{t-\tau}\bigr) \|^2 \right]  \\
& \hspace{1cm} \leq 2 \E\left[  \frac{1}{N} \sum_{i=1}^N  \Big \| \int_{t-\tau}^{t} \bF^i_{s-\tau}\,\dd s \Big \|^2 \right]  \\
&  \hspace{1cm} +  2 \E\left[  \frac1N\sum_{i=1}^N \Big \|   \sum_{\ell=1}^N
    \int_{t-\tau}^{t}
    \left[
        \delta_{i\ell}\sqrt{\sigma^i}\Diag( \bF^i_{s-\tau})
        -
        \frac1N\sqrt{\sigma^\ell}\Diag(\bF^\ell_{s-\tau})
    \right]\dd {\bf B}^\ell_s \Big \|^2 \right] \\
& \hspace{1cm} =:  2 \Big(  I_{21} +I_{22} \Big),    
\end{aligned}
\end{align}
where we used $|a + b|^2 \leq 2 (|a|^2 + |b|^2)$. \newline

\noindent $\bullet$~Case 1 (Estimate of $ I_{21}$):~We use the Cauchy--Schwarz inequality to get 
\begin{equation} \label{C-5-2} I_{21}  \leq
    \E\left[
        \frac1N\sum_{i=1}^N
        \left(
            \int_{t-\tau}^{t}1\,\dd s
        \right)
        \left(
            \int_{t-\tau}^{t}
            \| \bF^i_{s-\tau} \|^2\,\dd s
        \right)
    \right] =
    \tau
    \int_{t-\tau}^{t}
    \E\left[
        \frac1N\sum_{i=1}^N
        \| \bF^i_{s-\tau} \|^2
    \right]\dd s.
\end{equation}

\noindent $\bullet$~Case 2 (Estimate of $I_{22}$):~
We use It\^o's isometry and \eqref{C-3} to obtain
\begin{align}
\begin{aligned} \label{C-5-3}
I_{22}  &=
    \int_{t-\tau}^{t}
    \E\left[
        \frac1N\sum_{i,\ell=1}^N
      \Big  \|
            \delta_{i\ell}\sqrt{\sigma^i}\Diag(\bF^i_{s-\tau})
            -
            \frac1N\sqrt{\sigma^\ell}\Diag(\bF^\ell_{s-\tau})
        \Big \|_{\mathrm{F}}^2
    \right]\dd s\\
    &=
    \left(1-\frac1N\right)
    \int_{t-\tau}^{t}
    \E\left[
        \frac1N\sum_{\ell=1}^N
        \sigma^\ell \| \bF^\ell_{s-\tau} \|^2
    \right]\dd s\le \sigma^{\infty} \int_{t-\tau}^{t}
    \E\left[
        \frac1N\sum_{\ell=1}^N
        \| \bF^\ell_{s-\tau} \|^2
    \right]\dd s.
\end{aligned}
\end{align}
In \eqref{C-5-1}, we combine the estimates \eqref{C-5-2} and \eqref{C-5-3} to get 
\begin{equation} \label{C-5-4}
\E\left[ \frac1N\sum_{i=1}^N \| \bigl(\bv^i_t- \bar{\bv}_t \bigr)
    -  \bigl( \bv^i_{t-\tau}-\bar{\bv}_{t-\tau}\bigr) \|^2 \right] \leq 2 (\tau + \sigma^{\infty})  \int_{t-2\tau}^{t-\tau}
    \E\left[
        \frac1N\sum_{\ell=1}^N
        \| \bF^\ell_{s} \|^2
    \right]\dd s.
\end{equation}
Again, we use \eqref{C-5-0}, \eqref{C-5-4} and Young's 
inequality to see
\begin{align*}
\begin{aligned} \label{C-5-5}
    \left| \Lambda^{c}(t,t-\tau)  \right| &\leq
    \left(
        2(\tau+\sigma^{\infty})
        \int_{t-2\tau}^{t-\tau}\Lambda^\omega(s) \,\dd s
    \right)^{1/2}
    (\Lambda^{\omega}(t-\tau))^{1/2} \\
    & \leq \varepsilon \Lambda^{\omega}(t-\tau)
    +
    \frac{(\tau+\sigma^{\infty})}{2\varepsilon}
    \int_{t-2\tau}^{t-\tau}\Lambda^{\omega}(s)\,\dd s.
\end{aligned}
\end{align*}
\end{proof}
In what follows, we derive a dissipation estimate for ${\mathcal E}$ using a modified energy functional. 
\subsection{Compensated energy functional} \label{sec:3.2} 
In this subsection, we first derive a differential inequality for
${\mathcal E}$ and then introduce a compensated energy functional
to obtain a dissipation estimate.
\begin{lemma} \label{L3.5}
Suppose the framework $({\mathcal F}_A1) - ({\mathcal F}_A3)$ holds, and let  $(\bX_t, \bV_t)$ be a strong solution process to \eqref{B-1}. Then, for any $\varepsilon > 0$, we have
\begin{equation} \label{C-6}
 \frac{\di{\mathcal E}(t)}{\di t} \leq   - \Big (1 - \frac{\sigma^{\infty}}{2} - 2 \varepsilon \Big) {\Lambda}^{\omega}(t-\tau) +  \frac{(\tau+\sigma^{\infty})}{\varepsilon}
    \int_{t-2\tau}^{t-\tau}\Lambda^{\omega}(s) \,\dd s, \quad \mbox{a.e.}~~t \geq \tau. 
 \end{equation}
\end{lemma}
\begin{proof} 
We combine all the estimates in Lemma \ref{L3.2}, Lemma \ref{L3.3} and Lemma \ref{L3.4} to get 
\begin{align*}
\begin{aligned}
\frac{\di{\mathcal E}(t)}{\di t} 
 &\leq
   - {\Lambda}^{\omega}(t-\tau)
    + 2 \Lambda^{c}(t,t-\tau)
    +  \Lambda^{\sigma}(t-\tau) \\
    &\leq  - \Big (1 - \frac{\sigma^{\infty}}{2} - 2 \varepsilon \Big) {\Lambda}^{\omega}(t-\tau) +  \frac{(\tau+\sigma^{\infty})}{\varepsilon}
    \int_{t-2\tau}^{t-\tau}\Lambda^{\omega}(s) \,\dd s, \quad t \geq \tau.
\end{aligned}
\end{align*}
\end{proof}
To prepare the delay compensation in Lemma~\ref{L3.5}, for
$0\leq \sigma^\infty<2$, we set
\[
\delta:=\delta_{\sigma^\infty}
:=\frac{2-\sigma^\infty}{3}>0,
\qquad
B:=B_{\sigma^\infty}
:=\frac{2(1+\sigma^\infty)}{\delta_{\sigma^\infty}}
=\frac{6(1+\sigma^\infty)}{2-\sigma^\infty}>0.
\]
Motivated by the delay-compensation approach used in \cite{H},
we introduce the following compensated energy functional and the
associated admissible time-delay bound:
\begin{equation}\label{C-11}
	\begin{cases}
		\displaystyle
		\mathcal L(t)
		:=
		\mathcal E(t)
		+
		B\int_{t-2\tau}^{t-\tau}
		(s-t+2\tau)\Lambda^\omega(s)\,\dd s,
		\qquad t\geq\tau,\vspace{6pt}\\
		\displaystyle
		\tau_{*,\sigma^\infty}
		:=
		\frac{\delta_{\sigma^\infty}}
		{4B_{\sigma^\infty}}
		=
		\frac{(2-\sigma^\infty)^2}
		{72(1+\sigma^\infty)}.
	\end{cases}
\end{equation}
Then, it is easy to see that $\mathcal L(t)$ is absolutely continuous and
\begin{equation}\label{C-11-0}
	\mathcal L(t)\geq \mathcal E(t)\geq0.
\end{equation}
\begin{proposition} \label{P3.1}
Suppose the framework $({\mathcal F}_A1) - ({\mathcal F}_A3)$ holds, and we also assume that system parameters $\tau$ and $\sigma^{\infty}$ satisfy
\begin{equation} \label{C-11-1}
\quad 0 \le \sigma^{\infty} <2, \quad  0\le  \tau \le  \tau_{*,\sigma^\infty},  
\end{equation}
let  $(\bX_t, \bV_t)$ be a strong solution process to \eqref{B-1}. Then, the following assertions hold.
\begin{enumerate}
\item[(i)]
The compensated energy functional ${\mathcal L}(t)$ satisfies the dissipation estimate:
\begin{equation}\label{C-12}
\frac{\di}{\di t}\mathcal L(t)
\leq
-\frac{2-\sigma^\infty}{12}
\Lambda^\omega(t-\tau),
\qquad \mbox{a.e. } t\geq\tau.
\end{equation}

\item[(ii)]
The alignment production functional $\Lambda^{\omega}(t)$ is integrable:
\begin{equation} \label{C-13}
\int_{0}^{\infty} \Lambda^{\omega}(t) \dd t \le  \frac{12}{2-\sigma^\infty} {\mathcal L}(\tau). 
\end{equation}
\end{enumerate}
\end{proposition}
\begin{proof}
\noindent (i)~We split the proof into two steps. \newline

\noindent $\bullet$~Step A (Temporal evolution of ${\mathcal E}$):~We first refine the estimate \eqref{C-6} for some $\varepsilon$. In Lemma \ref{L3.5}, we choose $\varepsilon = \frac{1}{2} \delta$  to get 
\begin{equation} \label{C-14}
 \frac{\di{\mathcal E}(t)}{\di t} \leq - \frac{1}{2}\Big (2 - \sigma^{\infty} - 2\delta \Big ) {\Lambda}^{\omega}(t-\tau) +  \frac{2(\tau+\sigma^{\infty})}{ \delta}
    \int_{t-2\tau}^{t-\tau}\Lambda^{\omega}(s) \,\dd s, \quad \mbox{a.e.}~~t \geq \tau. 
\end{equation}
Note that the relations in \eqref{C-11-1} imply
\[
    2-\sigma^{\infty}=
     3 \delta, \quad \mbox{i.e.} \quad    2-\sigma^{\infty} -2\delta \geq\delta. \]
Then, we use
\begin{equation} \label{C-14-1}
    \tau\leq\tau_{*,\sigma^\infty}\leq1,
\end{equation}
the definition of $\tau_{*,\sigma^\infty}$ in \eqref{C-11} to see
\begin{equation} \label{C-14-2}
    \frac{2(\tau+\sigma^{\infty})}{\delta}
    \leq
    \frac{2(1+\sigma^\infty)}{\delta} = B.
\end{equation}
Now, we use \eqref{C-14}, \eqref{C-14-1} and \eqref{C-14-2} to get 
\begin{equation}
 \frac{\di{\mathcal E}(t)}{\di t} 
    \leq
    -\frac{\delta}{2} \Lambda^\omega(t-\tau)
    +
    B \int_{t-2\tau}^{t-\tau}\Lambda^{\omega}(s)\,\dd s
    \qquad
    \mbox{a.e.}~t\geq\tau.
    \label{C-15}
\end{equation}

\vspace{0.2cm}

\noindent $\bullet$~Step B (Derivation of dissipation estimate for ${\mathcal L}$):~By Leibniz' rule, one has 
\begin{equation} \label{C-15-1}
 \frac{\dd}{\dd t}
    \int_{t-2\tau}^{t-\tau}
    (s-t+2\tau)\Lambda^{\omega}(s) \dd s=
    \tau\Lambda^{\omega}(t-\tau)
    -
    \int_{t-2\tau}^{t-\tau}\Lambda^{\omega}(s)\ \dd s, \quad \mbox{a.e.}~t \geq \tau.
\end{equation}
Now, we combine \eqref{C-11}, \eqref{C-15},
and \eqref{C-15-1} to obtain the desired estimate:
\begin{align*}
\begin{aligned} 
\frac{\di{\mathcal L}(t)}{\di t}
    &= \frac{\di \mathcal E(t)}{\di t}   +
    B \frac{\di}{\di t} \int_{t-2\tau}^{t-\tau}
    (s-t+2\tau)\Lambda^{\omega}(s) \,\dd s \\
    &\leq   -\frac{\delta}{2} \Lambda^{\omega}(t-\tau)
    +
    B \int_{t-2\tau}^{t-\tau}\Lambda^{\omega}(s) \,\dd s + B \Big(  \tau\Lambda^{\omega}(t-\tau)
    -
    \int_{t-2\tau}^{t-\tau}\Lambda^{\omega}(s)\,\dd s \Big) \\
    &= -\left(\frac{\delta}{2} - B \tau\right) \Lambda^{\omega}(t-\tau)\le -\frac{\delta}{4}\Lambda^{\omega}(t-\tau) = - \frac{2-\sigma^\infty}{12} \Lambda^{\omega}(t-\tau),
\end{aligned}
\end{align*}    
i.e.,
\begin{equation} \label{C-15-2}
\frac{\di{\mathcal L}(t)}{\di t} \leq  -\frac{2-\sigma^\infty}{12} \Lambda^{\omega}(t-\tau), \quad \mbox{a.e.}~t \geq \tau. 
\end{equation}

\vspace{0.5cm}

\noindent (ii)~ Next we integrate \eqref{C-15-2} over \([\tau,t]\) to get 
\[
{\mathcal L}(t)
    + \frac{2-\sigma^\infty}{12} \int_{0}^{t-\tau} \Lambda^{\omega}(s) \dd s \leq {\mathcal L}(\tau), \quad t \geq \tau.
\]
Since $ {\mathcal L}(t) \geq 0$, we have
\begin{equation*}
\frac{2-\sigma^\infty}{12} \int_{0}^{t-\tau} \Lambda^{\omega}(s) \dd s \leq {\mathcal L}(\tau).
\end{equation*}
Letting \(t\to\infty\), we have 
\[ 
\int_0^{\infty} \Lambda^{\omega} (t) \dd t \leq  \frac{12}{2-\sigma^\infty} {\mathcal L}(\tau) < \infty. 
\]
\end{proof}

\vspace{0.2cm}

\section{Proof of Proposition \ref{P2.1}} \label{sec:4}
\setcounter{equation}{0}
We now prove Proposition~\ref{P2.1}. Recall that $p\geq4$, the initial
history satisfies the moment bounds in \eqref{B-4-1}, and
$(\bX_t,\bV_t)$ is the unique global strong solution process to
\eqref{B-1}. It remains to prove that, for sufficiently small
$\sigma^\infty$ and $\tau$,
\begin{equation}\label{D-0-0}
	\sup_{-\tau\leq s<\infty}
	\E\left[
	\frac1N\sum_{i=1}^N
	\norm{\bv_s^i-\bar{\bv}_s}^{p}
	\right]
	\leq
	C_p\bigl(1+M_{p,v}^0\bigr),
\end{equation}
and
\begin{equation}\label{D-0-1}
	\sup_{-\tau\leq s\leq t}
	\E\left[
	\frac1N\sum_{i=1}^N
	\norm{\bx_s^i-\bar{\bx}_s}^{p}
	\right]
	\leq
	C_p\bigl(1+M_{p,x}^0+M_{p,v}^0\bigr)(1+t)^p.
\end{equation}
\subsection{Preparatory estimates}\label{sec:4.1}
In this subsection, we provide several elementary estimates on empirical
velocity and position \(p\)-th moments. To avoid the common empirical martingale
term appearing in the centered velocity equation, we set
\begin{equation}\label{D-0}
	\begin{aligned}
		\hat{\bv}_t^i
		&:=\bv_t^i-\bar{\bv}_0,
		\quad i\in[N],\quad
		{\mathcal D}_p(t):=
		\E\left[
		\frac1{N^2}\sum_{i,j=1}^N
		\omega_t^{ij}
		\|\hat{\bv}_t^i\|^{p-2}
		\hat{\bv}_t^i\cdot
		(\hat{\bv}_t^i-\hat{\bv}_t^j)
		\right].
	\end{aligned}
\end{equation}
\begin{lemma}
	\label{L4.1}
	Let $p\geq 2$. Then, the following assertions hold.
	\begin{enumerate}
		\item[(i)] The relative velocity $\hat{\bv}^i_t = \bv^i_t - {\bar \bv}_0$ satisfies 
		\begin{equation} \label{D-2}
			\dd\hat{\bv}^i_t
			=
			\bF^i_{t-\tau}\dd t
			+
			\sqrt{\sigma^i}\,\Diag(\bF^i_{t-\tau})\dd {\bf B}^i_t,
		\end{equation}
		and
		 \begin{equation}\label{D-4-3}
		 \bF^i_{t-\tau} =
		\frac1N\sum_{j=1}^N
		\omega^{ij}_{t-\tau}
		\bigl(
		\hat{\bv}^j_{t-\tau}-\hat{\bv}^i_{t-\tau}
		\bigr).\end{equation}
		\item[(ii)] The average $p$-th~moment of velocity fluctuation satisfies \begin{equation}
			\frac1N\sum_{i=1}^N
			\| \bv^i_t- \bar{\bv}_t \|^p
			\leq
			2^p   \frac1N\sum_{i=1}^N\norm{\hat{\bv}^i_t}^p,
			\quad t\geq-\tau.
			\label{D-3}
		\end{equation}
	\end{enumerate}
\end{lemma}
\begin{proof}
	\noindent (i)~Recall the equation $\eqref{B-1}_2$ for $\bv^i_t$:
	\begin{align}
		\begin{aligned} \label{D-4-1}
			\dd \bv^i_t
			=
			\bF^i_{t-\tau} \dd t
			+ \sqrt{\sigma^i}\,\Diag(\bF^i_{t-\tau})\dd {\bf B}^i_t.
		\end{aligned}
	\end{align}
	Since $\bar{\bv}_0$ is independent of time $t$, we have
	\[
	\dd\hat{\bv}_t^i
	=
	\dd(\bv_t^i-\bar{\bv}_0)
	=
	\dd\bv_t^i.
	\]
	Hence, \eqref{D-2} follows from \eqref{D-4-1}. Note that 
	\[
		\hat{\bv}^j_t-\hat{\bv}^i_t = \bv^j_t-  \bar{\bv}_0 - (\bv^i_t-  \bar{\bv}_0) = \bv^j_t-\bv^i_t.
	\]
	This implies 
	\[
			\bF^i_{t-\tau} = \frac{1}{N} \sum_{j=1}^N
			\omega^{ij}_{t-\tau}\bigl( \bv^j_{t-\tau}- \bv^i_{t-\tau})
			=
			\frac1N\sum_{j=1}^N
			\omega^{ij}_{t-\tau}
			\bigl(
			\hat{\bv}^j_{t-\tau}-\hat{\bv}^i_{t-\tau}
			\bigr).
			\]
This verifies \eqref{D-4-3}.  \newline
	
	\noindent (ii)~Note that 
	\[
	\bar{\hat{\bv}}_t 
	:=
	\frac1N\sum_{i=1}^N\hat{\bv}^i_t
	=
	{\bar \bv}_t - {\bar \bv}_0.
	\]
	This implies 
	\begin{equation} \label{D-4-4}
		\bv^i_t- \bar{\bv}_t 
		=
		\hat{\bv}^i_t  + {\bar \bv}_0 -  \bar{\bv}_t  =  \hat{\bv}^i_t  - \bar{\hat{\bv}}_t.
	\end{equation}
	By Jensen's inequality, one has 
	\begin{equation} \label{D-4-5}
		\| \bar{\hat \bv}_t \|^p =  \Big \| \frac1N \sum_{i=1}^{N} \hat{\bv}^i_t  \Big \|^p
		\leq
		\frac1N\sum_{i=1}^N\norm{\hat{\bv}^i_t}^p.
	\end{equation}
	Then we use \eqref{D-4-4}, the inequality:
	\[ |a + b|^p \leq 2^{p-1} (|a|^p + |b|^p) \quad \mbox{for $p \geq 2$}, \]
	\eqref{D-0} and \eqref{D-4-5} to get the desired estimate  \eqref{D-3}:
		\begin{align*} 
			\frac1N\sum_{i=1}^N
			\| \bv^i_t- \bar{\bv}_t \|^p &= \frac1N\sum_{i=1}^N \|  \hat{\bv}^i_t  - \bar{\hat{\bv}}_t \|^p  \leq  \frac{2^{p-1}}{N} \sum_{i=1}^N \Big( \|  \hat{\bv}^i_t \|^p + \| \bar{\hat{\bv}}_t \|^p        \Big) \\
			&  \leq
			2^{p-1}
			\frac1N\sum_{i=1}^N\norm{\hat{\bv}^i_t}^p
			+   2^{p-1} \| \bar{ \hat \bv}_t \|^p \leq
			2^p   \frac1N\sum_{i=1}^N\norm{\hat{\bv}^i_t}^p.
		\end{align*}
\end{proof}
Next, we collect several auxiliary estimates.
Detailed proofs of Lemmas~\ref{L4.2} and~\ref{L4.4}
are given in Appendix~\ref{App-A}.
\begin{lemma}
	\label{L4.2}
	Let \(p\geq2\) be given, and let
	$
	0\leq \omega^{ij}=\omega^{ji}\leq1$ for all $i,j\in[N].$ Then, there exists a constant \(C_p>0\) such that, for every \(\bZ = (\bz^1, \cdots, \bz^N)  \in(\R^d)^N\),  with
	\[
	\bF^i
	:=
	\frac1N\sum_{j=1}^N
	\omega^{ij}(\bz^j-\bz^i),
	\qquad i\in[N],
	\]
	the following estimates hold:
	\begin{align}
		\begin{aligned}  \label{D-5}
			& (i)~ \frac1N\sum_{i=1}^N  \| \bF^i \|^p \leq \frac{C_p}{N} \sum_{i=1}^N \| \bz^i \|^p. \\
			& (ii)~ \frac1N\sum_{i=1}^N \| \bF^i \|^p \leq \frac{C_p}{N^2} \sum_{i, j =1}^N  \omega^{ij}
			\Big( \| \bz^i \|^{p-2} \bz^i \cdot ( \bz^i- \bz^j) \Big). \\
			& (iii)~\frac{1}{N} \sum_{i=1}^N \| \bz^i \|^{p-2}  \| \bF^i  \|^2 \leq
			\frac{C_p}{N^2} \sum_{i, j =1}^N  \omega^{ij}
			\Big( \| \bz^i \|^{p-2} \bz^i \cdot ( \bz^i- \bz^j) \Big). \\
			& (iv)~0\leq  \frac{1}{N^2} \sum_{i, j =1}^N  \omega^{ij} 
			\Big( \| \bz^i \|^{p-2} \bz^i \cdot ( \bz^i- \bz^j) \Big) \leq
			\frac{C_p}{N} \sum_{i=1}^N \| \bz^i \|^p.
		\end{aligned}
	\end{align}
\end{lemma}

\begin{proof}
	The estimates follow from Jensen's inequality, the symmetry of the
	weights, and the estimates in Lemma~\ref{L2.3} for $\| \bz \|^{p-2} \bz$.  Direct proofs will be given in Appendix~\ref{App-A.1}.
\end{proof}
\begin{lemma}
	\label{L4.3}
 For every \(T>0\), there exists a constant
 $C_{p,T,\tau,\sigma^\infty}>0$ independent of $N$ such that
	\[
	\E\left[
	\sup_{0\leq t\leq T}    \frac1N\sum_{i=1}^N\norm{\hat{\bv}^i_t}^p
	\right]\le
C_{p,T,\tau,\sigma^\infty}M_{p,v}^0.
	\]
\end{lemma}

\begin{proof}
	Fix $T>0$. We first consider the case $\tau>0$.
	As in Lemma~\ref{L3.2}, we work with the stopped process and
	remove the stopping at the end by Fatou's lemma. For $k\geq0$ and $t\in[k\tau,(k+1)\tau]$, we integrate
	\eqref{D-2} over $[k\tau,t]$, and use H\"older's inequality
	in time and the Burkholder--Davis--Gundy inequality
	\cite[Theorem~7.3, p.~40]{M} to obtain
		\begin{equation}\label{D-5-0}
		\begin{aligned}
			\E\left[
			\sup_{k\tau\leq t\leq(k+1)\tau}
			\frac1N\sum_{i=1}^N
			\norm{\hat{\bv}_t^i}^p
			\right]
			&\leq
			C_p
			\E\left[
			\frac1N\sum_{i=1}^N
			\norm{\hat{\bv}_{k\tau}^i}^p
			\right]
			\\
			&\hspace{-2cm}+
			C_p
			\left(
			\tau^{p-1}
			+
			(\sigma^\infty)^{p/2}\tau^{\frac p2-1}
			\right)
			\int_{k\tau}^{(k+1)\tau}
			\E\left[
			\frac1N\sum_{i=1}^N
			\norm{\bF_{s-\tau}^i}^p
			\right]\dd s.
		\end{aligned}
		\end{equation}
	 On the other hand, by Lemma~\ref{L4.2}(i) we have
	\begin{align*}
		&\int_{k\tau}^{(k+1)\tau}
		\E\left[
		\frac1N\sum_{i=1}^N
		\norm{\bF_{s-\tau}^i}^p
		\right]\dd s\leq
		C_p\tau
		\E\left[
		\sup_{(k-1)\tau\leq t\leq k\tau}
		\frac1N\sum_{i=1}^N
		\norm{\hat{\bv}_t^i}^p
		\right].
	\end{align*}
	Therefore, with \eqref{D-5-0}, we have
	\begin{equation}\label{D-5-1}
	\E\left[
	\sup_{k\tau\leq t\leq(k+1)\tau}
	\frac1N\sum_{i=1}^N
	\norm{\hat{\bv}_t^i}^p
	\right]\le C_{p,\tau,\sigma^\infty}\E\left[
	\sup_{(k-1)\tau\leq t\leq k\tau}
	\frac1N\sum_{i=1}^N
	\norm{\hat{\bv}_t^i}^p
	\right],\qquad k\geq1.
	\end{equation}
We next prove the desired result by induction with \eqref{D-5-1}. We consider the first step $k=0$. By Lemma~\ref{L4.2}(i) and \eqref{B-4-1}, we have
\[
\int_0^\tau
\E\left[
\frac1N\sum_{i=1}^N
\norm{\bF_{s-\tau}^i}^p
\right]\dd s\le C_p\int_0^\tau\E\left[\frac{1}{N} \sum_{i=1}^N \| \bv^i_{s-\tau}-\bar{\bv}_{s-\tau} \|^p\right]\dd s
\leq
C_p\tau M_{p,v}^0,\]
\[\E\left[
\frac1N\sum_{i=1}^N
\norm{\hat{\bv}_0^i}^p
\right]=\E\left[\frac1N\sum_{i=1}^N\|\bpsi_0^i-\bar{\bpsi}_0\|^p\right]
\leq M_{p,v}^0.
\]
Therefore, with \eqref{D-5-0}, we have
\[
\E\left[
\sup_{0\leq t\leq\tau}
\frac1N\sum_{i=1}^N
\norm{\hat{\bv}_t^i}^p
\right]
\le C_{p,\tau,\sigma^\infty}M_{p,v}^0.
\]
Hence, with \eqref{D-5-1}, by induction, we get
	\begin{align*}
		\E\left[
		\sup_{0\leq t\leq k\tau}
		\frac1N\sum_{i=1}^N
		\norm{\hat{\bv}_t^i}^p
		\right]\le C_{p,\tau,\sigma^\infty,k}M_{p,v}^0,
	\end{align*}
	for every integer $k\geq1$. Choosing $k$ sufficiently large such that $T\leq k\tau$ gives the desired estimate.
	For $\tau=0$, we use the same localized H\"older inequality and  Burkholder--Davis--Gundy inequality estimate,
	together with Lemma~\ref{L4.2}(i) to get
	\[
	\E\left[
	\sup_{0\leq s\leq t}
	\frac1N\sum_{i=1}^N
	\norm{\hat{\bv}_s^i}^p
	\right]
	\leq
	C_{p,T,\sigma^\infty}
	\left(
	\E\left[
	\frac1N\sum_{i=1}^N
	\norm{\hat{\bv}_0^i}^p
	\right]
	+
	\int_0^t
	\E\left[
	\sup_{0\leq r\leq s}
	\frac1N\sum_{i=1}^N
	\norm{\hat{\bv}_r^i}^p
	\right]\dd s
	\right),
	\]
	for $0\leq t\leq T$. The desired estimate follows from
	Gronwall's inequality. 
\end{proof}

\begin{lemma}
	\label{L4.4}
	Suppose that 
	\[ 0\leq\sigma^{\infty} \leq1 \quad \mbox{and} \quad 0\leq\tau\leq1. \]
	Then, there exists a positive constant \(C_p>0\), independent of \(N\) such that the following assertions hold.
	\begin{enumerate}
		\item[(i)]
		For almost
		every \(t\geq2\tau\), we have
		\begin{equation} \label{D-6}
			\E\left[
			\frac1N\sum_{i=1}^N
			\norm{\hat{\bv}^i_t-\hat{\bv}^i_{t-\tau}}^p
			\right]
			\leq
			C_p(\tau+\sigma^{\infty})
			\int_{t-2\tau}^{t-\tau}{\mathcal D}_p(s)\,\dd s.
		\end{equation}
		Here, ${\mathcal D}_p$ is defined in \eqref{D-0}.
		\vspace{0.2cm}
		\item[(ii)]
		For almost every \(t\geq3\tau\), we have
		\begin{equation} \label{D-7}
			\E\left[
			\frac1N\sum_{i=1}^N
			\| \hat{\bv}^i_{t-\tau} \|^{p-2}
			\| \hat{\bv}^i_t-\hat{\bv}^i_{t-\tau} \|^2
			\right]
			\leq
			C_p(\tau+\sigma^{\infty})
			\int_{t-3\tau}^{t-\tau}{\mathcal D}_p(s)\,\dd s.
		\end{equation}
	\end{enumerate}
\end{lemma}
\begin{proof}
	Since the proofs are very lengthy, we leave them in Appendix \ref{App-A.2}.
\end{proof}

\subsection{Temporal evolution of
	$\displaystyle
	\E\left[\frac1N\sum_{i=1}^N\|\hat{\bv}_t^i\|^p\right]$}

We apply It\^o's formula to
$\|\hat{\bv}_t^i\|^p$. Recall that
\[
\dd \hat{\bv}_t^i
=
\bF_{t-\tau}^i\,\dd t
+
\sqrt{\sigma^i}\,
\Diag(\bF_{t-\tau}^i)\,\dd {\bf B}_t^i .
\]
Since
\[
\nabla \|\bz\|^p
=
p\|\bz\|^{p-2}\bz,\quad
D^2\|\bz\|^p
=
p\|\bz\|^{p-2}I_d
+
p(p-2)\|\bz\|^{p-4}\bz\bz^{\mathsf T},
\]
by It\^o's formula, we have
\begin{align*}
	\dd\|\hat{\bv}_t^i\|^p
	&=
	p\|\hat{\bv}_t^i\|^{p-2}
	\hat{\bv}_t^i\cdot\bF_{t-\tau}^i\,\dd t
	+
	p\sqrt{\sigma^i}
	\|\hat{\bv}_t^i\|^{p-2}
	\hat{\bv}_t^i\cdot
	\Diag(\bF_{t-\tau}^i)\,\dd\mathbf{B}_t^i
	\\
	&+
	\frac{p\sigma^i}{2}
	\|\hat{\bv}_t^i\|^{p-2}
	\|\bF_{t-\tau}^i\|^2\,\dd t
	+
	\frac{p(p-2)\sigma^i}{2}
	\|\hat{\bv}_t^i\|^{p-4}
	\big\|
	\Diag(\bF_{t-\tau}^i)\hat{\bv}_t^i
	\big\|^2\,\dd t .
\end{align*}
By Appendix~\ref{App-A.3}, the stochastic integral has zero expectation. We sum the above equality over $i$, divide it by $N$, and take expectation, to obtain that for almost every
$t\geq 3\tau$,
\begin{equation}\label{D-12-0}
	\begin{aligned}	
	\frac{\dd}{\dd t}
	\E\left[
	\frac1N\sum_{i=1}^N
	\|\hat{\bv}_t^i\|^p
	\right]
	={}&
	p\E\left[
	\frac1N\sum_{i=1}^N
	\|\hat{\bv}_t^i\|^{p-2}
	\hat{\bv}_t^i\cdot\bF_{t-\tau}^i
	\right]\\
	&\hspace{-2cm}+
	\frac{p}{2N}\sum_{i=1}^N\sigma^i
	\E\Bigg[
	\|\hat{\bv}_t^i\|^{p-2}
	\|\bF_{t-\tau}^i\|^2
	+
	(p-2)
	\|\hat{\bv}_t^i\|^{p-4}
	\big\|
	\Diag(\bF_{t-\tau}^i)\hat{\bv}_t^i
	\big\|^2
	\Bigg].
\end{aligned}
\end{equation}
We add and subtract
$\|\hat{\bv}_{t-\tau}^i\|^{p-2}\hat{\bv}_{t-\tau}^i$ to decompose the first term as
\begin{align*}
	&\E\left[
	\frac1N\sum_{i=1}^N
	\|\hat{\bv}_t^i\|^{p-2}
	\hat{\bv}_t^i\cdot\bF_{t-\tau}^i
	\right]\\
&\hspace{0.5cm}=
	\E\left[
	\frac1N\sum_{i=1}^N
	\|\hat{\bv}_{t-\tau}^i\|^{p-2}
	\hat{\bv}_{t-\tau}^i
	\cdot\bF_{t-\tau}^i
	\right]
	+
	\E\left[
	\frac1N\sum_{i=1}^N
	\left(
	\|\hat{\bv}_t^i\|^{p-2}\hat{\bv}_t^i
	-
	\|\hat{\bv}_{t-\tau}^i\|^{p-2}
	\hat{\bv}_{t-\tau}^i
	\right)
	\cdot\bF_{t-\tau}^i
	\right].
\end{align*}
By \eqref{D-4-3}, we have
\begin{align*}
\E\left[
\frac1N\sum_{i=1}^N
\|\hat{\bv}_{t-\tau}^i\|^{p-2}
\hat{\bv}_{t-\tau}^i
\cdot\bF_{t-\tau}^i
\right]=
-\mathcal D_p(t-\tau).
\end{align*}
We now define
\begin{align*}
	\mathcal I_p(t,t-\tau)
	:=&
	\frac{p}{2N}\sum_{i=1}^N\sigma^i
	\E\Bigg[
	\|\hat{\bv}_t^i\|^{p-2}
	\|\bF_{t-\tau}^i\|^2
	+
	(p-2)
	\|\hat{\bv}_t^i\|^{p-4}
	\big\|
	\Diag(\bF_{t-\tau}^i)\hat{\bv}_t^i
	\big\|^2
	\Bigg],\\
	\mathcal R_p(t,t-\tau)
	:=&
	\E\left[
	\frac1N\sum_{i=1}^N
	\left(
	\|\hat{\bv}_t^i\|^{p-2}\hat{\bv}_t^i
	-
	\|\hat{\bv}_{t-\tau}^i\|^{p-2}
	\hat{\bv}_{t-\tau}^i
	\right)
	\cdot\bF_{t-\tau}^i
	\right].
\end{align*}
Then, the relation \eqref{D-12-0} becomes
\begin{equation}
	\frac{\dd}{\dd t}
	\E\left[
	\frac1N\sum_{i=1}^N
	\|\hat{\bv}_t^i\|^p
	\right]
	=
	-p\mathcal D_p(t-\tau)
	+
	p\mathcal R_p(t,t-\tau)
	+
	\mathcal I_p(t,t-\tau).
	\label{D-12}
\end{equation}

Note that the first term in the right-hand side of \eqref{D-12} has a definite sign and provides a dissipation term. Thus we focus on the last two terms $\mathcal{R}_p$ and $\mathcal{I}_p$ in the sequel. 

\subsubsection{Estimate of  $|\mathcal{R}_p|$} By the Cauchy--Schwarz inequality, one has 
\begin{equation} \label{D-12-0-0}
	|{\mathcal R}_p(t,t-\tau)|
	\leq
	\E\left[
	\frac1N\sum_{i=1}^N
	\Big \| 
	\| \hat{\bv}^i_t \|^{p-2} \hat{\bv}^i_t 
	-
	\| \hat{\bv}^i_{t-\tau} \|^{p-2}  \hat{\bv}^i_{t-\tau}
	\Big \|  \| \bF^i_{t-\tau} \| \right].
\end{equation}
On the other hand, it follows from Lemma~\ref{L2.3}(iii) that
\begin{align}
	\begin{aligned} \label{D-12-0-1}
		& \Big\|  \| \hat{\bv}^i_t \|^{p-2} \hat{\bv}^i_t 
		-
		\| \hat{\bv}^i_{t-\tau} \|^{p-2}  \hat{\bv}^i_{t-\tau} \Big \| \leq C_p \Big (   \| \hat{\bv}^i_{t-\tau} \|^{p-2}  \|  \hat{\bv}^i_t  -  \hat{\bv}^i_{t-\tau} \| +  \|  \hat{\bv}^i_t  -  \hat{\bv}^i_{t-\tau} \|^{p-1} \Big).
	\end{aligned}
\end{align}            
By combining  \eqref{D-12-0-0} and \eqref{D-12-0-1}, we have
\begin{align}
	\begin{aligned} \label{D-12-1-1}
		|\mathcal R_p(t,t-\tau)|
		& \leq
		C_p \Big( \E\left[
		\frac1N\sum_{i=1}^N
		\| \hat{\bv}^i_{t-\tau} \|^{p-2}
		\| \hat{\bv}^i_t-\hat{\bv}^i_{t-\tau} \|
		\| \bF^i_{t-\tau} \|
		\right] \\
		& +  \E\left[
		\frac1N\sum_{i=1}^N \| \hat{\bv}^i_t-\hat{\bv}^i_{t-\tau} \|^{p-1}
		\norm{\bF^i_{t-\tau}}
		\right] \Big) \\
		& = C_p \Big(I_{31}(t,t-\tau) + I_{32}(t,t-\tau) \Big).
	\end{aligned}
\end{align}    
In the following lemma, we estimate \(\mathcal R_p\). 
\begin{lemma} \label{L4.5}
	Suppose that
	\[  0\leq\sigma^{\infty} \leq1 \quad \mbox{and} \quad 0\leq\tau\leq1. \]
	 Let  $(\bX_t, \bV_t)$ be a strong solution process to \eqref{B-1}. Then, for any $\eta\in(0,1]$ and for almost every $t\geq3\tau$, one has
	\begin{align}\label{D-13}
		\begin{aligned}
			  |\mathcal R_p(t,t-\tau)|
			\leq
			C_p\eta {\mathcal D}_p(t-\tau)
			+
			C_{p,\eta}(\tau+\sigma^{\infty})
			\int_{t-3\tau}^{t-\tau}{\mathcal D}_p(s)\,\dd s.
		\end{aligned}
	\end{align}
\end{lemma}
\begin{proof}
	\noindent (i)~For $\eta\in(0,1]$, we use Young's inequality to get 
	\begin{align}
		\begin{aligned} \label{D-12-3}
			I_{31}(t,t-\tau)
			&\leq
			\eta
			\E\left[
			\frac1N\sum_{i=1}^N
			\| \hat{\bv}^i_{t-\tau} \|^{p-2}
			\norm{\bF^i_{t-\tau}}^2
			\right]  +
			C_{p,\eta}
			\E\left[
			\frac1N\sum_{i=1}^N
			\| \hat{\bv}^i_{t-\tau} \|^{p-2}
			\| \hat{\bv}^i_t-\hat{\bv}^i_{t-\tau} \|^2
			\right]  \\
			&=: I_{311}(t-\tau) + I_{312}(t,t-\tau).
		\end{aligned}
	\end{align}
	Below, we estimate the terms $I_{31i},~i= 1,2$ one by one. \newline
	
	\noindent $\bullet$ (Estimate of $I_{311}$):~ We use Lemma \ref{L4.2} (iii) to get 
	\begin{equation} \label{D-12-4}
		I_{311}(t-\tau) \leq \eta C_p {\mathcal D}_p(t-\tau).
	\end{equation}   
	\noindent $\bullet$ (Estimate of $I_{312}$):~We use the second assertion in Lemma \ref{L4.4} to get 
	\begin{equation} \label{D-12-5}
		I_{312}(t,t-\tau) \leq C_{p,\eta} (\tau + \sigma^{\infty}) \int_{t-3\tau}^{t-\tau} {\mathcal D}_p(s) \dd s.
	\end{equation}
	Finally, we combine \eqref{D-12-3}, \eqref{D-12-4} and \eqref{D-12-5} to get the desired estimate:
	\begin{equation} \label{D-12-6}
		I_{31}(t,t-\tau)
		\leq
		C_p\eta {\mathcal D}_p(t-\tau)
		+
		C_{p,\eta}(\tau+\sigma^{\infty})
		\int_{t-3\tau}^{t-\tau}{\mathcal D}_p(s)\,\dd s.
	\end{equation}
	\noindent (ii)~For \(I_{32}\),  we use Young's inequality with conjugate exponents
	\(p\) and \(p/(p-1)\) to find 
	\[
	\begin{aligned}
		I_{32}(t,t-\tau)
		&\leq
		\eta
		\E\left[
		\frac1N\sum_{i=1}^N \| \bF^i_{t-\tau}\|^p
		\right]+
		C_{p,\eta}
		\E\left[
		\frac1N\sum_{i=1}^N \| \hat{\bv}^i_t-\hat{\bv}^i_{t-\tau} \|^p
		\right].
	\end{aligned}
	\]
	Using \eqref{D-5}$_2$ and \eqref{D-6}, we obtain
	\begin{equation} \label{D-12-7}
		I_{32}(t,t-\tau)
		\leq
		C_p\eta {\mathcal D}_p(t-\tau)
		+
		C_{p,\eta}(\tau+\sigma^\infty)
		\int_{t-2\tau}^{t-\tau}{\mathcal D}_p(s)\,\dd s.
	\end{equation}
	In \eqref{D-12-1-1}, we combine estimates \eqref{D-12-6} and \eqref{D-12-7} to find the desired estimate \eqref{D-13}.
\end{proof}

\begin{lemma} \label{L4.6}
	Suppose that 
	\[  0\leq\sigma^{\infty} \leq1 \quad \mbox{and} \quad 0\leq\tau\leq1.\]
Then, for any $\eta\in(0,1]$ and for almost every
	$t\geq3\tau$, the functional $\mathcal{I}_p(t,t-\tau)$ satisfies 
	\[
	\mathcal I_p(t,t-\tau)
	\leq
	C_p(\sigma^{\infty}+\eta){\mathcal D}_p(t-\tau)
	+
	C_{p,\eta}(\tau+\sigma^{\infty})
	\int_{t-2\tau}^{t-\tau}
	{\mathcal D}_p(s)\,\dd s.
	\]
\end{lemma}
\begin{proof}
Since
\[
\big\|
\Diag(\bF_{t-\tau}^i)\hat{\bv}_t^i
\big\|^2
\leq
\|\bF_{t-\tau}^i\|^2
\|\hat{\bv}_t^i\|^2,
\]
we have
\begin{equation*}
	\mathcal I_p(t,t-\tau)
	\leq
	\frac{p(p-1)}{2}\sigma^\infty
	\E\left[
	\frac1N\sum_{i=1}^N
	\|\hat{\bv}_t^i\|^{p-2}
	\|\bF_{t-\tau}^i\|^2
	\right].
\end{equation*}
	On the other hand, we note that
	\begin{equation} \label{D-13-2}
		\|\hat{\bv}^i_t\|^{p-2}
		\leq
		C_p\|\hat{\bv}^i_{t-\tau}\|^{p-2}
		+
		C_p
		\|\hat{\bv}^i_t-\hat{\bv}^i_{t-\tau}\|^{p-2}.
	\end{equation}
	Hence, one has
	\begin{align*}
			&\E\left[
			\frac1N\sum_{i=1}^N
			\|\hat{\bv}^i_t\|^{p-2}
			\|\bF^i_{t-\tau}\|^2
			\right]\\
			&\hspace{0.4cm}\leq
			C_p
			\E\left[
			\frac1N\sum_{i=1}^N
			\|\hat{\bv}^i_{t-\tau}\|^{p-2}
			\|\bF^i_{t-\tau}\|^2
			\right]
			+
			C_p
			\E\left[
			\frac1N\sum_{i=1}^N
			\|\hat{\bv}^i_t-\hat{\bv}^i_{t-\tau}\|^{p-2}
			\|\bF^i_{t-\tau}\|^2
			\right].
	\end{align*}
	For any $\eta\in(0,1]$, we use Young's inequality to see
	\begin{equation} \label{D-13-3}
		\|\hat{\bv}^i_t-\hat{\bv}^i_{t-\tau}\|^{p-2}
		\|\bF^i_{t-\tau}\|^2
		\leq
		\eta\|\bF^i_{t-\tau}\|^p
		+
		C_{p,\eta}
		\|\hat{\bv}^i_t-\hat{\bv}^i_{t-\tau}\|^p.
	\end{equation}
	Now, we use \eqref{D-13-2}, \eqref{D-13-3},
	Lemma~\ref{L4.2} (ii), (iii), Lemma~\ref{L4.4},
	and the same arguments as in the proof of Lemma~\ref{L4.5}
	to obtain
	\begin{align*}
			& \E\left[
			\frac1N\sum_{i=1}^N
			\|\hat{\bv}^i_t\|^{p-2}
			\|\bF^i_{t-\tau}\|^2
			\right]
			\\
			& \quad\leq
			C_p\Bigg(
			\E\left[
			\frac1N\sum_{i=1}^N
			\|\hat{\bv}^i_{t-\tau}\|^{p-2}
			\|\bF^i_{t-\tau}\|^2
			\right]
			+
			\E\left[
			\frac1N\sum_{i=1}^N
			\|\hat{\bv}^i_t-\hat{\bv}^i_{t-\tau}\|^{p-2}
			\|\bF^i_{t-\tau}\|^2
			\right]
			\Bigg)
			\\
			& \quad\leq
			C_p(1+\eta)\mathcal D_p(t-\tau)
			+
			C_{p,\eta}(\tau+\sigma^{\infty})
			\int_{t-2\tau}^{t-\tau}
			\mathcal D_p(s)\,\dd s.
		\end{align*}
	This implies
	\begin{align*}
		\begin{aligned}
			\mathcal I_p(t,t-\tau)
			&\leq
			C_p\sigma^{\infty}(1+\eta)
			\mathcal D_p(t-\tau)+
			C_{p,\eta}\sigma^{\infty}
			(\tau+\sigma^{\infty})
			\int_{t-2\tau}^{t-\tau}
			\mathcal D_p(s)\,\dd s
			\\
			&\leq
			C_p(\sigma^{\infty}+\eta)
			\mathcal D_p(t-\tau)
			+
			C_{p,\eta}(\tau+\sigma^{\infty})
			\int_{t-2\tau}^{t-\tau}
			\mathcal D_p(s)\,\dd s,
		\end{aligned}
	\end{align*}
	where we used
	\(0\leq\sigma^{\infty}\leq1\) and \(0<\eta\leq1\)
	in the last inequality.
\end{proof}
Finally, we use \eqref{D-12}, Lemma \ref{L4.5} and Lemma \ref{L4.6} to derive a differential inequality for $\displaystyle \E \left [  \frac1N\sum_{i=1}^N\norm{\hat{\bv}^i_t}^p \right ]$ in the following lemma. 
\begin{lemma} \label{L4.7}
	Suppose that
\[  0\leq\tau\leq1. \]
 Then, there exists
	$\sigma_{p,*}\in(0,1)$, depending only on $p$, such that, whenever
	\[
	0\leq\sigma^\infty\leq\sigma_{p,*},
	\] 
	we have
	\[
	\frac{\dd}{\dd t} \E \left [  \frac1N\sum_{i=1}^N\norm{\hat{\bv}^i_t}^p \right ] \leq
	-\frac{p}{2}{\mathcal D}_p(t-\tau)
	+
	B_p\int_{t-3\tau}^{t-\tau}{\mathcal D}_p(s)\,\dd s,
	\qquad
	\text{a.e. } t\geq3\tau,
	\]
	where \(B_p\) is a nonnegative constant given by
	\[
	B_p:= (p + 1) C_{p,\eta_p}(\tau+\sigma^\infty),
	\]
	with \(\eta_p\in(0,1]\) chosen in the proof below.
\end{lemma}
\begin{proof} It follows from \eqref{D-12}, Lemma \ref{L4.5} and Lemma \ref{L4.6} that 
	\begin{align*}
		\begin{aligned} 
			& \frac{\dd}{\dd t} \E \left [  \frac1N\sum_{i=1}^N\norm{\hat{\bv}^i_t}^p \right ] \\
			& \hspace{1cm} =-p{\mathcal D}_p(t-\tau)+p\mathcal{R}_p(t,t-\tau)+\mathcal{I}_p(t,t-\tau) \\
			& \hspace{1cm} \leq -p{\mathcal D}_p(t-\tau) + p \Big(    C_p\eta {\mathcal D}_p(t-\tau)
			+
			C_{p,\eta}(\tau+\sigma^{\infty})
			\int_{t-3\tau}^{t-\tau}{\mathcal D}_p(s)\,\dd s \Big) \\
			& \hspace{1cm}  + C_p (\sigma^{\infty} + \eta) {\mathcal D}_p(t-\tau)
			+
			C_{p,\eta}(\tau+\sigma^{\infty})
			\int_{t-2\tau}^{t-\tau}{\mathcal D}_p(s)\,\dd s  \\
			& \hspace{1cm}  \le \Big(  -p +  p C_p\eta + C_p (\sigma^{\infty} + \eta)  \Big) {\mathcal D}_p(t-\tau)  + (p +1) C_{p,\eta}(\tau+\sigma^{\infty}) \int_{t-3\tau}^{t-\tau}{\mathcal D}_p(s)\,\dd s.
		\end{aligned}
	\end{align*}
	For a given \(p\geq4\), we choose
	\(\eta=\eta_p\in(0,1]\) and \(\sigma_{p,*}\in(0,1)\) sufficiently small such that
	\[
	(p+1)C_p\eta_p\leq\frac p4,
	\qquad
	C_p\sigma_{p,*}\leq\frac p4.
	\]
	This implies that, for \(0\leq\sigma^\infty\leq\sigma_{p,*}\), we have
	\[
	pC_p\eta_p+C_p(\sigma^\infty+\eta_p)
	\leq\frac p2.
	\]
	Therefore, we set
	\[
	B_p:=(p+1)C_{p,\eta_p}(\tau+\sigma^\infty)
	\]
	to get the desired estimate.
\end{proof}

\subsection{Moments for spatial and velocity fluctuations} \label{sec:4.3} 
In this subsection, we are ready to verify the estimates \eqref{D-0-0} and \eqref{D-0-1} using lemmas obtained in the previous subsection. \newline

By \eqref{D-3}, it suffices to derive a uniform bound for
\[ \E\left[ \frac1N\sum_{i=1}^N
\| \bv^i_t- \bar{\bv}_0 \|^p \right]. \]
All the applications of It\^o's formula below are first made to the stopped
solution, and then we can use Lemma~\ref{L4.3} and
Appendix~\ref{App-A.3} to justify the passage to the original process. For simplicity of notation, we suppress the stopping notation. First, we recall the following differential inequality:
\begin{equation*}
	\frac{\dd}{\dd t} \E\left[ \frac1N\sum_{i=1}^N
	\| \bv^i_t- \bar{\bv}_0 \|^p \right] 
	\leq
	-\frac{p}{2}{\mathcal D}_p(t-\tau)
	+
	B_p\int_{t-3\tau}^{t-\tau}{\mathcal D}_p(s)\,\dd s,
	\quad
	t\geq3\tau.
\end{equation*}
For \(t\geq3\tau\), we define the new compensated functional:
\[
\mathcal L_p(t)
:= \E\left[ \frac1N\sum_{i=1}^N
\| \bv^i_t- \bar{\bv}_0 \|^p \right] 
+
B_p
\int_{t-3\tau}^{t-\tau}
(s-t+3\tau){\mathcal D}_p(s)\,\dd s.
\]
Then, it is easy to see that 
\begin{equation}
	\E\left[ \frac1N\sum_{i=1}^N
	\| \bv^i_t- \bar{\bv}_0 \|^p \right] 
	\leq
	\mathcal L_p(t),
	\quad t\geq3\tau.
	\label{D-15-2}
\end{equation}
\begin{lemma}
	\label{L4.8}
Let $\sigma_{p,*}\in(0,1)$ be the constant given in
Lemma~\ref{L4.7}. Then, there exists
$\tau_{p,*}\in(0,1]$, depending only on $p$, such that,  whenever
	\[
	0\leq\sigma^\infty\leq\sigma_{p,*}<1,
	\qquad
	0\leq\tau\leq\tau_{p,*}\le1,
	\]
	one has the following estimate:
	\begin{equation} \label{D-14}
		\frac{\di}{\di t} \mathcal L_p(t)
		\leq
		-\frac{p}{4}{\mathcal D}_p(t-\tau),  \quad\text{ a.e.} ~~t \ge 3\tau.
	\end{equation}
\end{lemma}
\begin{proof} 
	\noindent By Lemma~\ref{L4.2}(iv) and Lemma~\ref{L4.3}, \({\mathcal D}_p\) is locally integrable. Hence, we have
	\[
	\frac{\dd}{\dd t}
	\int_{t-3\tau}^{t-\tau}
	(s-t+3\tau){\mathcal D}_p(s)\,\dd s
	=
	2\tau{\mathcal D}_p(t-\tau)
	-
	\int_{t-3\tau}^{t-\tau}
	{\mathcal D}_p(s)\,\dd s, \quad \mbox{a.e.} ~t\geq3\tau.
	\]
We combine this identity with Lemma~\ref{L4.7} to obtain
	\begin{align*}
			\frac{\di \mathcal L_p(t)}{\di t} &=  \frac{\di}{\di t}  \E\left[ \frac1N\sum_{i=1}^N
			\| \bv^i_t- \bar{\bv}_0 \|^p \right] 
			+
			B_p \frac{\di}{\di t}
			\int_{t-3\tau}^{t-\tau}
			(s-t+3\tau){\mathcal D}_p(s)\,\dd s \\
			& \leq -\frac{p}{2}{\mathcal D}_p(t-\tau)
			+ B_p\int_{t-3\tau}^{t-\tau}{\mathcal D}_p(s)\,\dd s +   B_p  \Big( 2\tau {\mathcal D}_p(t-\tau)
			- \int_{t-3\tau}^{t-\tau}{\mathcal D}_p(s)\,\dd s \Big) \\
			&= -\bigl(\frac{p}{2}-2B_p\tau\bigr){\mathcal D}_p(t -\tau).
	\end{align*}
	Note that since \(0\leq\tau,\,\sigma^\infty\leq1\), one has 
	\[
	2 B_p \tau  = 2 (p + 1) C_{p,\eta_p}(\tau+\sigma^\infty) \tau \leq  4 (p + 1) C_{p,\eta_p} \tau.
	\]
	Now, we choose \[
	\tau_{p,*}
	:=
	\min\left\{
	1,\,
	\frac{p}{16(p+1)C_{p,\eta_p}}
	\right\}.
	\]  
	Then, whenever \(0\leq\tau\leq\tau_{p,*}\), we have
	\[
	2B_p\tau
	\leq
	\frac{p}{4}.
	\]
	Hence, we have the desired estimate \eqref{D-14}.
\end{proof}
Now, we are ready to give the proof of Proposition \ref{P2.1}.

\begin{proof}
	\noindent (i)~For the first estimate, the history interval
	$[-\tau,0]$ follows directly from \eqref{B-4-1}.
	For nonnegative times, we consider two cases:
	\[ \mbox{Either} \quad t \geq 3 \tau \quad \mbox{or} \quad 0 \leq t < 3 \tau. \]
	\noindent $\bullet$~Step A $(t \geq 3 \tau)$:~It follows from \eqref{D-3},  \eqref{D-15-2} and Lemma \ref{L4.8} that 
	\begin{equation} \label{D-15-5}
		\E\left[
		\frac1N\sum_{i=1}^N
		\| \bv^i_t- {\bar \bv}_t \|^p
		\right]  \leq 2^p  \E\left[
		\frac1N\sum_{i=1}^N
		\| \bv^i_t- {\bar \bv}_0 \|^p
		\right] \leq 2^p {\mathcal L}_p(t) \leq 2^p {\mathcal L}_p(3\tau), \quad t \geq 3 \tau.  
	\end{equation}
	In the next step, we will show that ${\mathcal L}_p(3\tau)$ is also finite. \newline
	
	\noindent $\bullet$~Step B $(0 \leq t < 3 \tau)$:~We claim that 
	\begin{align}
		\begin{aligned} \label{D-16}
			& (i)~ \sup_{0 \leq t\leq 3 \tau} \E\left[
			\frac1N\sum_{i=1}^N
			\| \bv^i_t- {\bar \bv}_0 \|^p
			\right] \leq  C_p\bigl(1+M_{p,v}^0\bigr). \\
			& (ii)~ \mathcal L_p(3\tau) \leq  C_p\bigl(1+M_{p,v}^0\bigr).
		\end{aligned}
	\end{align}
	{\it Proof of \eqref{D-16}}:~We first prove (i).
	For \(k=0,1,2,3\), we set
\[
H_k
:=
\sup_{0\leq t\leq k\tau}
\E\left[
\frac1N\sum_{i=1}^N
\|\bv_t^i-\bar{\bv}_0\|^p
\right].
\]
	Since $\hat{\bv}^i_0 = \bv^i_0 - {\bar \bv}_0$, the bounds \eqref{B-4-1} yield
	\begin{equation}\label{D-16-0}
	H_0
	=  \E\left[
	\frac1N\sum_{i=1}^N
	\| \bv^i_0- {\bar \bv}_0 \|^p
	\right]
	\leq
	M_{p,v}^0.
	\end{equation}
	On the first delay interval, we have
	\[
	\bF_{s-\tau}^i
	=
	\frac1N\sum_{j=1}^N
	\omega_{s-\tau}^{ij}
	\Big[
	(\bv_{s-\tau}^j-\bar{\bv}_{s-\tau})
	-
	(\bv_{s-\tau}^i-\bar{\bv}_{s-\tau})
	\Big].
	\]
	Thus Lemma~\ref{L4.2}(i) and conditions \eqref{B-4-1} imply
	\begin{equation}\label{D-16-1}
	\begin{aligned}
		&\int_0^\tau
		\E\left[
		\frac1N\sum_{i=1}^N
		\| \bF^i_{s-\tau} \|^p
		\right]\dd s\leq
		C_p\tau M_{p,v}^0.
	\end{aligned}
	\end{equation}
	For any \(t\in[0,\tau]\), we integrate \eqref{D-2} over \([0,t]\), then use H\"older's inequality in time and the
	Burkholder--Davis--Gundy inequality to obtain
	\begin{align}
		\begin{aligned}
			\E\left[
			\frac1N\sum_{i=1}^N
			\|\hat{\bv}_t^i\|^p
			\right]
			&\leq
			C_p H_0
			+
			C_p t^{p-1}
			\int_0^t
			\E\left[
			\frac1N\sum_{i=1}^N
			\|\bF_{s-\tau}^i\|^p
			\right]\dd s
			\\
			&+
			\frac{C_p}{N}
			\sum_{i=1}^N
			(\sigma^i)^{p/2}
			\E\left[
			\left(
			\int_0^t
			\|\Diag(\bF_{s-\tau}^i)\|_{\mathrm F}^2
			\,\dd s
			\right)^{p/2}
			\right]
			\\
			&\leq
			C_p H_0
			+
			C_p t^{p-1}
			\int_0^t
			\E\left[
			\frac1N\sum_{i=1}^N
			\|\bF_{s-\tau}^i\|^p
			\right]\dd s
			\\
			&+
			C_p(\sigma^\infty)^{p/2}
			t^{\frac p2-1}
			\int_0^t
			\E\left[
			\frac1N\sum_{i=1}^N
			\|\bF_{s-\tau}^i\|^p
			\right]\dd s,
		\end{aligned}
		\label{D-16-2}
	\end{align}
	where the Burkholder--Davis--Gundy inequality gives
	\[
	\E\left[
	\left\|
	\int_0^t
	\Diag(\bF_{s-\tau}^i)\,\dd\mathbf{B}_s^i
	\right\|^p
	\right]
	\leq
	C_p
	\E\left[
	\left(
	\int_0^t
	\|\Diag(\bF_{s-\tau}^i)\|_{\mathrm F}^2
	\,\dd s
	\right)^{p/2}
	\right],
	\]
	and the last inequality is due to 
	\[
	\|\Diag(\bF_{s-\tau}^i)\|_{\mathrm F}
	=
	\|\bF_{s-\tau}^i\|,
	\]
	and H\"older's inequality in time:
	\[
	\left(
	\int_0^t
	\|\bF_{s-\tau}^i\|^2\,\dd s
	\right)^{p/2}
	\leq
	t^{\frac p2-1}
	\int_0^t
	\|\bF_{s-\tau}^i\|^p\,\dd s.
	\]
	Hence, by \eqref{D-16-0}, \eqref{D-16-1} and \eqref{D-16-2}, together with \(0\leq\tau,\sigma^\infty\leq1\), we see 
	\[
	\begin{aligned}
		H_1=\sup_{0\le t\le\tau}\E\left[
		\frac1N\sum_{i=1}^N
		\|\hat{\bv}_t^i\|^p
		\right]
		&\leq
		C_p H_0
		+
		C_p
		\left(
		\tau^p
		+
		(\sigma^\infty)^{p/2}\tau^{p/2}
		\right)
		M_{p,v}^0\leq
		C_p\bigl(1+M_{p,v}^0\bigr).
	\end{aligned}
	\]
	For \(k=1,2\), if
	\(
	s\in[k\tau,(k+1)\tau],
	\)
	then
	\[
	s-\tau\in[(k-1)\tau,k\tau].
	\]
	Consequently, Lemma~\ref{L4.2}(i) gives
	\[
	\begin{aligned}
		&\int_{k\tau}^{(k+1)\tau}
		\E\left[
		\frac1N\sum_{i=1}^N
		\| \bF^i_{s-\tau} \|^p
		\right]\dd s\leq
		C_p\tau H_k.
	\end{aligned}
	\]
	We apply the same integral estimate on
	\([k\tau,(k+1)\tau]\) to obtain
	\[
	H_{k+1}
	\leq
	C_p(1+H_k),
	\qquad
	k=1,2.
	\]
	This yields
	\begin{equation}
		\sup_{0\leq t\leq3\tau}
		\E\left[
		\frac1N\sum_{i=1}^N
		\| \bv^i_t- {\bar \bv}_0 \|^p
		\right]
		=
		H_3
		\leq
		C_p\bigl(1+M_{p,v}^0\bigr).
		\label{D-17}
	\end{equation}
	Moreover, it follows from Lemma~\ref{L4.2}(iv) that 
	\begin{equation}\label{D-17-0}
		\int_0^{2\tau}{\mathcal D}_p(s)\,\dd s
		\leq
		C_p
		\int_0^{2\tau} \E\left[
		\frac1N\sum_{i=1}^N
		\| \bv^i_{s}- {\bar \bv}_0 \|^p
		\right] \dd s \leq
		C_p\bigl(1+M_{p,v}^0\bigr).
	\end{equation}
	At \(t=3\tau\), the memory term becomes
	\[
	\int_0^{2\tau}s{\mathcal D}_p(s)\,\dd s.
	\]
	Since \({\mathcal D}_p(s)\geq0\) and \(0\leq s\leq2\tau\), we have
	\[
	\int_0^{2\tau}s{\mathcal D}_p(s)\,\dd s
	\leq
	2\tau\int_0^{2\tau}{\mathcal D}_p(s)\,\dd s
	\leq
	2\tau C_p\bigl(1+M_{p,v}^0\bigr).
	\]
Combining \eqref{D-17} and \eqref{D-17-0}, and using the boundedness of
\(B_p\tau\) under the chosen parameter restrictions, we obtain
\begin{equation*}
	\mathcal L_p(3\tau)
	=
	\E\left[
	\frac1N\sum_{i=1}^N
	\| \bv^i_{3\tau}- {\bar \bv}_0 \|^p
	\right]
	+
	B_p
	\int_0^{2\tau}
	s{\mathcal D}_p(s)\,\dd s
	\leq
	C_p\bigl(1+M_{p,v}^0\bigr).
\end{equation*}
	\noindent $\bullet$~Final Step:~For $t \geq 3 \tau$, with \eqref{D-15-5} and \eqref{D-16}, we have
	\[
	\E\left[
	\frac1N\sum_{i=1}^N
	\| \bv^i_t- {\bar \bv}_t \|^p
	\right]  \leq 2^p {\mathcal L}_p(3\tau) \leq C_p\bigl(1+M_{p,v}^0\bigr), \quad t \geq 3 \tau.  
	\]
	For $0 \leq t < 3\tau$, one has
	\[
	\E\left[
	\frac1N\sum_{i=1}^N
	\| \bv^i_t- {\bar \bv}_t \|^p
	\right]  \leq 2^p \E\left[
	\frac1N\sum_{i=1}^N
	\| \bv^i_t- {\bar \bv}_0 \|^p
	\right]  \leq
	C_p\bigl(1+M_{p,v}^0\bigr).
	\]
%	Finally, the relation \eqref{D-3} yields
%	\[
%	\begin{aligned}
%		\sup_{t\geq0}
%		\E\left[
%		\frac1N\sum_{i=1}^N
%		\| \bv^i_t- \bar{\bv}_t \|^p
%		\right]
%		\leq
%		2^p 
%		\sup_{t\geq0} \E\left[
%		\frac1N\sum_{i=1}^N
%		\| \bv^i_t- {\bar \bv}_0 \|^p
%		\right]        \leq
%		C_p\bigl(1+M_{p,v}^0\bigr).
%	\end{aligned}
%	\]
	The corresponding estimate on \([-\tau,0]\) follows directly from
	conditions \eqref{B-4-1}. This verifies \eqref{D-0-0}.
	
	\vspace{0.2cm}
	
	\noindent (ii)~For \(t\geq0\), we have
	\[
	\bx^i_t- \bar{\bx}_t
	= \bx^i_0 -\bar{\bx}_0 +
	\int_0^t (\bv^i_s -\bar{\bv}_s \bigr)\,\dd s.
	\]
	This implies 
	\[
	\frac1N\sum_{i=1}^N
	\| \bx^i_t - {\bar \bx}_t \|^p
	\leq \frac{C_p}{N} \sum_{i=1}^N
	\| \bx^i_0 - {\bar \bx}_0 \|^p
	+ C_p t^{p-1}
	\int_0^t  \Big(
	\frac1N\sum_{i=1}^N
	\| \bv^i_s - {\bar \bv}_s \|^p \Big) \dd s.
	\]
	Now, we take an expectation and use the first assertion to obtain 
	\begin{align*}
		\begin{aligned}
			\E\left[ \frac1N\sum_{i=1}^N
			\| \bx^i_t - {\bar \bx}_t \|^p \right]
			&\leq C_p \E\left[ \frac{1}{N} \sum_{i=1}^N
			\| \bx^i_0 - {\bar \bx}_0 \|^p \right]
			+ C_p t^{p-1}
			\int_0^t  \E\left[ 
			\frac1N\sum_{i=1}^N
			\| \bv^i_s - {\bar \bv}_s \|^p \right] \dd s \\
			&\leq C_p \E\left[ \frac{1}{N} \sum_{i=1}^N
			\| \bx^i_0 - {\bar \bx}_0 \|^p \right] + C_p \bigl(1+M_{p,v}^0\bigr) t^p \\
			&\leq C_p \bigl(
			1+M_{p,x}^0+M_{p,v}^0
			\bigr)(1+t)^p.
		\end{aligned}
	\end{align*}
	We take the supremum over \(s\in[0,t]\), and use conditions \eqref{B-4-1} on \([-\tau,0]\) to derive the desired estimate $\eqref{D-0-1}$.
\end{proof}
\begin{remark}
	Compared to the kinetic CS model with non-compact support \cite{H-W,H-W1} and stochastic CS model \cite{H-Z}, the estimates of velocity and position moments are more delicate. Therefore, we introduce a memory-corrected Lyapunov functional to derive the desired estimates.
\end{remark}
\vspace{0.2cm}
\section{Derivation of a time-dependent lower bound for $\Lambda^{\omega}$} \label{sec:5}
\setcounter{equation}{0}
In this section, we derive a time-dependent lower bound for $\Lambda^\omega(t)$ in terms of the
velocity fluctuation:
\[
\Lambda^\omega(t)
\geq
c(1+t)^{-\beta\alpha}\mathcal E(t)
-
C(1+t)^{-\beta\alpha-(p-2)(\alpha-1)},
\]
where the exponent $\alpha$ will be chosen in Section~6.1.
This estimate will be combined with the dissipation inequality obtained
in Proposition~\ref{P3.1} in the proof of the main flocking result.
A uniform positive lower bound is not available globally, since the
communication weight may be small when the particle configuration
spreads out. We therefore
use the method of time-varying effective region to compensate for the
degeneration of the communication weight.
\begin{lemma}
	\label{L5.1}
	Assume that the framework
	$({\mathcal F}_{A}1)$--$({\mathcal F}_{A}3)$ holds. Let  $(\bX_t, \bV_t)$ be a strong solution process to \eqref{B-1}. Then, we have for every \(R>0\) and \(t\geq 0\), 
	\begin{equation}\label{D-19}
		\Lambda^\omega(t)
		\geq
		\frac{2^{-\beta}}{(1+R)^{\beta}} \mathcal E(t)
		-\frac{2^{-\beta}}{(1+R)^{\beta}} \left ( T^R(t)
		+ \E\left[
		q^R_{t}
		\frac1N\sum_{i=1}^N
		\|
		\bv^i_{t}-\bar{\bv}_{t}
		\|^2 
		\right] \right).
	\end{equation}
	Here, $T^R(t)$ and $q_t^R$ are defined as follows. 
	\[
	T^R(t)
	:=
	\E\left[
	\frac1N\sum_{i=1}^N
	\norm{ \bv^i_{t}-\bar{\bv}_{t}}^2
	\one_{\{ \|\bx^i_{t}- \bar{\bx}_{t} \| >R\}}
	\right], \quad 
	q^R_{t}
	:=
	\frac1N\sum_{i=1}^N
	\one_{\{\|\bx^i_{t}- {\bar \bx}_{t} \| >R\}}.
	\]
\end{lemma}
\begin{proof}
	For notational simplicity, we first fix \(t\geq0\), \(R>0\), and a sample point, and we set
	\begin{align*}
		\begin{aligned}
			& {\tilde \bv}^i_{t}:= \bv^i_{t}- \bar{\bv}_{t}, \quad
			{\tilde \bx}^i_{t}:=\bx^i_{t}- \bar{\bx}_{t}, \quad i \in [N], \\
			& \mathcal G^R_{t}
			:=
			\left\{
			i \in [N]: \| \tilde{\bx}^i_{t} \| \leq R
			\right\},
			\quad
			\mathcal B^R_{t}
			:=
			\left\{
			i \in [N]: \| \tilde{\bx}^i_{t} \|>R
			\right\},
			\quad
			q_{t}^R:=\frac{|\mathcal B^R_{t}|}{N}.
		\end{aligned}
	\end{align*}
	For \(i,j\in\mathcal G^R_{t}\), we have
	\[
	\| \bx^i_{t}-\bx^j_{t} \|
	=
	\norm{\tilde{\bx}^{i}_{t}- \tilde{\bx}^{j}_{t}}
	\leq 2R.
	\]
	Therefore, the lower bound in \eqref{A-1} yields
	\[
	\begin{aligned}
		\omega^{ij}_{t} =  \omega\left(
		\| \bx^i_{t}-\bx^j_{t} \|
		\right)
		=
		\frac{1}{
			\bigl(1+\|\bx_t^i-\bx_t^j\|^2\bigr)^{\beta/2}}
		\geq
		(1+4R^2)^{-\beta/2}
		\geq
		2^{-\beta}(1+R)^{-\beta}.
	\end{aligned}
	\]  
	Consequently, we have
	\[
	\Lambda^\omega(t)
	\geq
	2^{-\beta}(1+R)^{-\beta}
	\E\left[
	\frac1{2N^2}
	\sum_{i,j\in\mathcal G^R_{t}}
	\| \tilde{\bv}^j_t - \tilde{\bv}^i_t \|^2 \right].
	\]
	Next, we estimate the pairwise sum over \(\mathcal G^R_{t}\). The
	pairwise-variance identity gives
	\[
	\begin{aligned}
		\frac1{2N^2}
		\sum_{i,j\in\mathcal G^R_{t}}
		\| \tilde{\bv}^j_t - \tilde{\bv}^i_t \|^2
		&=
		\frac{|\mathcal G^R_{t}|}{N}
		\frac1N\sum_{i\in\mathcal G^R_{t}} \| \tilde{\bv}^i_{t} \|^2
		-
		\frac1{N^2}  \Big \|
		\sum_{i\in\mathcal G^R_{t}} \tilde{\bv}^i_t  \Big \|^2.
	\end{aligned}
	\]
	Since
	\[
	\sum_{i\in\mathcal G^R_{t}} \tilde{\bv}^i_t 
	=
	-\sum_{i\in\mathcal B^R_{t}}  \tilde{\bv}^i_t,
	\]
	by the Cauchy--Schwarz inequality, we have
	\begin{align*}
		\begin{aligned}
			\frac1{N^2}
			\Big \|
			\sum_{i\in\mathcal G^R_{t}} \tilde{\bv}^i_t 
			\Big \|^2
			&=
			\frac1{N^2}
			\Big \|
			\sum_{i\in\mathcal B^R_{t}} {\tilde \bv}^i_t  \Big \|^2\leq
			\frac{|\mathcal B^R_{t}|}{N^2}
			\sum_{i\in\mathcal B^R_{t}} \| \tilde{\bv}^i_t \|^2=
			q^R_{t}
			\frac1N\sum_{i\in\mathcal B^R_{t}} \| \tilde{\bv}_t^i \|^2.
		\end{aligned}
	\end{align*}
	Since 
	\[ \frac{|\mathcal G^R_{t}|}{N} =1-q^R_{t}, \]
	we have
	\begin{align*}
		\begin{aligned}
			& \frac1{2N^2}
			\sum_{i,j\in\mathcal G^R_{t}} 
			\| \tilde{\bv}^i_t - \tilde{\bv}^j_t \|^2  \geq
			\frac{(1-q^R_{t})}{N}  \sum_{i\in\mathcal G^R_{t}} \| \tilde{\bv}^i_t \|^2 
			- \frac{q^R_{t}}{N} \sum_{i\in\mathcal B^R_{t}}  \| \tilde{\bv}^i_t \|^2  \\
			& \hspace{1cm}= \frac{(1-q^R_{t})}{N} 
			\sum_{i=1}^N \| \tilde{\bv}^i_t \|^2 
			-
			\frac1N\sum_{i\in\mathcal B^R_{t}} \| \tilde{\bv}^i_t \|^2 =
			\frac1N\sum_{i=1}^N \| \tilde{\bv}^i_t \|^2 
			-
			\frac1N\sum_{i\in\mathcal B^R_{t}} \| \tilde{\bv}^i_t \|^2 
			-
			\frac{q^R_{t}}{N} \sum_{i=1}^N \| \tilde{\bv}^i_t \|^2.
		\end{aligned}
	\end{align*}
	We substitute this estimate into the preceding lower bound, taking expectations, and recalling the definitions
	of \(T^R(t)\) and \(q^R_{t} \), we obtain \eqref{D-19}.
\end{proof}
Lemma~\ref{L5.1} shows that the loss of alignment force is entirely due to
particles lying outside the ball of radius \(R\). Indeed, \(q^R_{t}\)
measures the proportion of such particles, whereas \(T^R(t)\) measures
the centered velocity energy carried by them. 
\begin{lemma}
	\label{L5.2}
	Assume that the framework
	$({\mathcal F}_{A}1)$--$({\mathcal F}_{A}3)$ holds. Let  $(\bX_t, \bV_t)$ be a strong solution process to \eqref{B-1}. Then, for every \(R>0\) and \(t\geq0\),
	\begin{equation} \label{D-19-0}
		T^R(t)
		+
		\E\left[
		q^R_{t}
		\frac1N\sum_{i=1}^N
		\norm{ \bv^i_{t}-\bar{\bv}_{t}}^2
		\right]\leq\tilde{C}_p\frac{(1+t)^{p-2}}{R^{p-2}},
	\end{equation}
where 
\(\tilde{C}_p:=2C_p
\bigl(1+M_{p,v}^0\bigr)^{2/p}
\bigl(1+M_{p,x}^0+M_{p,v}^0\bigr)^{(p-2)/p}\)
with $C_p>0$ defined in Proposition~\ref{P2.1} which is independent of \(N, R\).
\end{lemma}
\begin{proof}
	On the set
	\(
	\left\{
	\| \bx^i_{t}- \bar{\bx}_{t} \| >R
	\right \},
	\)
	we have
	\[
	R^{-(p-2)}
	\| \bx^i_{t}- \bar{\bx}_{t} \|^{p-2}\ge1.
	\]
	Therefore, we have
	\[
	\begin{aligned}
		&\frac1N\sum_{i=1}^N
		\| \bv^i_{t}-\bar{\bv}_{t} \|^2
		\one_{\{\norm{\bx^i_{t}- \bar{\bx}_{t}}>R\}}
		\leq
		R^{-(p-2)}
		\frac1N\sum_{i=1}^N \| \bv^i_{t}-\bar{\bv}_{t} \|^2 \| \bx^i_{t}- \bar{\bx}_{t} \|^{p-2}.
	\end{aligned}
	\]
	We use H\"older's inequality for the empirical sum to find 
	\begin{align*} &\frac1N\sum_{i=1}^N
	\norm{\bv^i_{t}-\bar{\bv}_{t}}^2
	\one_{\{\norm{\bx^i_{t}- \bar{\bx}_{t}}>R\}} \\
	&\hspace{0.5cm}\leq
	R^{-(p-2)}
	\left(
	\frac1N\sum_{i=1}^N
	\norm{\bv^i_{t}-\bar{\bv}_{t}}^p
	\right)^{2/p}
	\left(
	\frac1N\sum_{i=1}^N
	\norm{\bx^i_{t}- \bar{\bx}_{t}}^p
	\right)^{(p-2)/p}.
	\end{align*}
	Now, we take expectations and apply H\"older's inequality in probability and Proposition~\ref{P2.1} to obtain
	\[
	\begin{aligned}
	T^R(t)
	&\leq
	\frac1{R^{p-2}}\!
	\left(
	\E\left[
	\frac1N\sum_{i=1}^N
	\norm{\bv^i_{t}-\bar{\bv}_{t}}^p
	\right]
	\right)^{2/p}\!\!
	\left(
	\E\left[
	\frac1N\sum_{i=1}^N
	\norm{\bx^i_{t}- \bar{\bx}_{t}}^p
	\right]
	\right)^{(p-2)/p}\\
	&\le C_p
	\bigl(1+M_{p,v}^0\bigr)^{2/p}
	\bigl(1+M_{p,x}^0+M_{p,v}^0\bigr)^{(p-2)/p}\frac{(1+t)^{p-2}}{R^{p-2}}.
	\end{aligned}
	\]
	Next, we estimate the second term in the left side of \eqref{D-19-0}. We use H\"older's inequality in probability to get
	\begin{equation}\label{D-19-0-0}
	\E \left [q^R_{t}  \frac1N\sum_{i=1}^N
	\norm{\bv^i_{t}-\bar{\bv}_{t}}^2 \right ] \leq
	\left(
	\E \left[ \Big(  \frac1N\sum_{i=1}^N
	\norm{\bv^i_{t}-\bar{\bv}_{t}}^2 \Big)^{p/2} \right ]
	\right)^{2/p}
	\left(
	\E[ (q^R_{t})^{p/(p-2)}]
	\right)^{(p-2)/p}.
	\end{equation}
	By Jensen's inequality, we have
	\[
	\Big(  \frac1N\sum_{i=1}^N
	\norm{\bv^i_{t}-\bar{\bv}_{t}}^2 \Big)^{p/2}
	\leq
	\frac1N\sum_{i=1}^N
	\norm{\bv^i_{t}-\bar{\bv}_{t}}^p.
	\]
	Moreover, since \(0\leq q^R_{t} \leq1\), we have
	\[
	( q^R_{t})^{p/(p-2)}
	\leq
	q^R_{t}.
	\]
	By the definition of \(q^R_{t}\), we have
	\[
	\begin{aligned}
		q^R_{t}
		&=
		\frac1N\sum_{i=1}^N
		\one_{\{\norm{\bx^i_{t}- \bar{\bx}_{t}}>R\}}\leq
		\frac1{R^p}
		\frac1N\sum_{i=1}^N
		\norm{\bx^i_{t}- \bar{\bx}_{t}}^p.
	\end{aligned}
	\]
Then, we have
	\[
	\left(\E\left[(q_t^R)^{p/(p-2)}\right]\right)^{(p-2)/p}
	\leq
	R^{-(p-2)}
	\left(
	\E\left[
	\frac1N\sum_{i=1}^N
	\|\bx_t^i-\bar{\bx}_t\|^p
	\right]
	\right)^{(p-2)/p}.
	\]
	We now go back to \eqref{D-19-0-0}. It follows from Proposition \ref{P2.1} that
	\[
	\begin{aligned}
		\E \left [q^R_{t} \Big(  \frac1N\sum_{i=1}^N
		\norm{\bv^i_{t}-\bar{\bv}_{t}}^2 \Big) \right ]
		&\leq
		C_p
		\bigl(1+M_{p,v}^0\bigr)^{2/p}
		\bigl(1+M_{p,x}^0+M_{p,v}^0\bigr)^{(p-2)/p}\frac{(1+t)^{p-2}}{R^{p-2}}.
	\end{aligned}
	\]
	Finally, we combine the two estimates to derive the desired estimate.
\end{proof}
As can be seen in the next proposition, moment estimates obtained in Lemma \ref{L5.2}
control the two terms in \eqref{D-19} and yield the following lower
bound for \(\Lambda^\omega\). For this, we set 
\begin{equation} \label{D-19-1}
	R(t):=(1+t)^\alpha \geq 1, \quad \alpha \geq 1,
\end{equation}
where the constant $\alpha$ will be chosen in Section~\ref{sec:6.1}. 
\begin{proposition}
	\label{P5.1}
	Assume that the framework
	$({\mathcal F}_{A}1)$--$({\mathcal F}_{A}3)$ holds. Let  $(\bX_t, \bV_t)$ be a strong solution process to \eqref{B-1}. Then, for $R(t)$ given by \eqref{D-19-1}, we have the following estimate:
	\begin{equation}
		\begin{aligned}
			\Lambda^\omega(t)
			&\geq
			2^{-2\beta}(1+t)^{-\beta\alpha}\mathcal E(t)-
			2^{-\beta} \tilde{C}_p(1+t)^{-\beta\alpha-(p-2)(\alpha-1)},
			\quad t\geq 0.
		\end{aligned}
		\label{D-20}
	\end{equation}
\end{proposition}
\begin{proof}
	Since $R(t) \geq 1$, we have
	\[
	1+R(t)\le2R(t).
	\]
	This and \eqref{D-19-1} imply
	\begin{align}
		\begin{aligned} \label{D-20-1}
			(1+R(t))^{-\beta}
			& \ge
			2^{-\beta}R(t)^{-\beta}
			=
			2^{-\beta}(1+t)^{-\alpha \beta}, \\
			(1+R(t))^{-\beta} &\le
			R(t)^{-\beta}
			=
			(1+t)^{-\alpha \beta}.
		\end{aligned}
	\end{align}
	We substitute \eqref{D-19-0} into \eqref{D-19}, and then use \eqref{D-19-1}, \eqref{D-20-1} to find the desired estimate:
	\begin{align*}
		\begin{aligned} 
			\Lambda^\omega(t)
			&\geq
			\frac{2^{-\beta}}{(1+R(t))^{\beta}} \mathcal E(t)
			-\frac{2^{-\beta}}{(1+R(t))^{\beta}} \left ( T^R(t)
			+ \E\left[
			q^R_{t}
			\frac1N\sum_{i=1}^N
			\|
			\bv^i_{t}-\bar{\bv}_{t}
			\|^2 
			\right] \right) \\
			& \geq \frac{2^{-\beta}}{(1+R(t))^{\beta}} \mathcal E(t)
			-\frac{2^{-\beta}\tilde{C}_p }{(1+R(t))^{\beta} R(t)^{p-2}} (1+t)^{p-2} \\
			& \geq 2^{-2\beta} (1 + t)^{-\alpha \beta} \mathcal E(t) -2^{-\beta} \tilde{C}_p (1+t)^{-\alpha \beta-(p-2)(\alpha-1)}.
		\end{aligned}
	\end{align*}
\end{proof}

\vspace{0.2cm}

\section{Mean-square stochastic flocking} \label{sec:6}
\setcounter{equation}{0}
In this section, we provide a proof of Theorem~\ref{T2.3}. First, we establish the polynomial decay of the mean-square velocity
alignment, and then using this quantitative velocity alignment estimate, we show the uniform empirical mean-square
confinement.
\subsection{Velocity alignment}\label{sec:6.1}
Suppose the framework $({\mathcal F}_A1) - ({\mathcal F}_A3)$ holds, and let  $(\bX_t, \bV_t)$ be a strong solution process to \eqref{B-1}. Let $p\geq 4$ be fixed, and 
\[
\sigma_*:=\sigma_{p,*},
\qquad
\tau_*:
=
\min\left\{
\tau_{p,*},
\frac1{144}
\right\},
\]
where \(\sigma_{p,*}\) and \(\tau_{p,*}\) are given by
Proposition~\ref{P2.1}, while \(1/144\) follows from
\eqref{C-11} by taking \(\sigma^\infty=1\). \newline

Suppose that
\[
0\leq\sigma^\infty\leq\sigma_*,
\qquad
0\leq\tau\leq\tau_*.
\]
We claim that, for \(0<\beta<1\) and
\[
0<\gamma<
\frac{(p-2)(1-\beta)}{\beta},
\]
or for \(\beta=0\) and any \(\gamma>0\) it holds
\begin{equation}\label{E-0}
	\E\left[
	\frac1N\sum_{i=1}^N
	\norm{\bv_t^i-\bar{\bv}_t}^2
	\right]
	\leq
	C_\gamma(1+t)^{-\gamma},
	\qquad t\geq0,
\end{equation}
where \(C_\gamma
\) is a positive constant independent of $N$.\newline

{\it Proof of \eqref{E-0}}:~The case \(\tau=0\) follows from the same argument without the memory
terms, so below we consider \(\tau>0\). We prove it in the following two steps.

\vspace{0.2cm}

\noindent$\bullet$~\textbf{Step A.} Recall \eqref{C-12} in Proposition~\ref{P3.1} at time $t+\tau$:
\begin{equation}
	\frac{\di}{\di t}\mathcal L(t+\tau)
	\leq
	-\frac{2-\sigma^\infty}{12}\Lambda^\omega(t),
	\qquad
	\text{for a.e. }t\geq0.
	\label{E-1}
\end{equation}
We combine \eqref{E-1} with \eqref{D-20} in Proposition~\ref{P5.1} to get
\begin{equation}\label{E-1-1}
	\begin{aligned}
		\frac{\dd}{\dd t}\mathcal L(t+\tau)
		&\leq-\frac{2-\sigma^\infty}{24}\Lambda^\omega(t)
		-
		\frac{2-\sigma^\infty}{24}\Lambda^\omega(t)\\
		&\le-\frac{2-\sigma^\infty}{24}\Lambda^\omega(t)
		-\frac{(2-\sigma^\infty)2^{-2\beta}}{24}
		(1+t)^{-\beta\alpha}\mathcal E(t)
		\\
		&
		+\frac{(2-\sigma^\infty)2^{-\beta}}{24}
		\tilde C_p
		(1+t)^{-\beta\alpha-(p-2)(\alpha-1)},
		\qquad \text{for a.e. }t\geq0.
	\end{aligned}
\end{equation}
Next, we compare \(\mathcal L(t+\tau)\) with the velocity
fluctuation $\mathcal{E}(t)$ from the definition of $\mathcal{L}$ in \eqref{C-11}. Note that the identity
\[
\begin{aligned}
	\bv^i_{t+\tau}- \bar{\bv}_{t+\tau}
	&=
	\bv^i_{t}-\bar{\bv}_{t}+
	\left[
	\bigl(\bv^i_{t+\tau}- \bar{\bv}_{t+\tau} \bigr)
	-
	\bigl( \bv^i_{t}-\bar{\bv}_{t}\bigr)
	\right]
\end{aligned}
\]
and the estimate \eqref{C-5} at time $t+\tau$ yield
\begin{equation*}
	\begin{aligned}
		\mathcal E(t+\tau)
		&\leq
		2\mathcal E(t)
		+
		2\E\left[
		\frac1N\sum_{i=1}^N
		\left\|
		\bigl(\bv^i_{t+\tau}-\bar{\bv}_{t+\tau}\bigr)
		-
		\bigl(\bv^i_t-\bar{\bv}_t\bigr)
		\right\|^2
		\right]
		\\
		&\leq
		2\mathcal E(t)
		+
		4(\tau+\sigma^\infty)
		\int_{t-\tau}^{t}\Lambda^\omega(s)\,\dd s.
	\end{aligned}
\end{equation*}
Moreover, by \eqref{C-11} we have
\[
0
\leq
B\int_{t-\tau}^{t}
(s-t+\tau)\Lambda^\omega(s)\,\dd s
\leq
B\tau
\int_{t-\tau}^{t}\Lambda^\omega(s)\,\dd s.
\]
Since \(\tau\leq1\) and \(\sigma^\infty\leq\sigma_*\), it follows from the definition of $\mathcal{L}$ that
\begin{equation}
	\mathcal L(t+\tau)
	\leq
	2\mathcal E(t)
	+
	\left[
	4(\tau+\sigma^\infty)
	+
	\frac{6\tau(1+\sigma^\infty)}
	{2-\sigma^\infty}
	\right]\int_{t-\tau}^{t}\Lambda^\omega(s)\,\dd s.
	\label{E-4}
\end{equation}
The algebraically decaying term in \eqref{E-1-1},
together with the comparison estimate \eqref{E-4}, suggests
considering a polynomially weighted Lyapunov functional.
We therefore consider
\[
(1+t)^\gamma\mathcal L(t+\tau),
\qquad \gamma>0,
\]
and differentiate it. By \eqref{E-1-1} and \eqref{E-4}, we obtain
\begin{equation}\label{E-4-1}
\begin{aligned}
	&\frac{\dd}{\dd t}
	\left[
	(1+t)^\gamma\mathcal L(t+\tau)
	\right]\\
	&\hspace{0.4cm}=
	\gamma(1+t)^{\gamma-1}\mathcal L(t+\tau)
	+
	(1+t)^\gamma
	\frac{\dd}{\dd t}\mathcal L(t+\tau)\\
	&\hspace{0.4cm}\leq
	\left[
	2\gamma(1+t)^{\gamma-1}
	-
	\frac{(2-\sigma^\infty)2^{-2\beta}}{24}
	(1+t)^{\gamma-\beta\alpha}
	\right]\mathcal E(t)
	-\frac{2-\sigma^\infty}{24}
	(1+t)^\gamma\Lambda^\omega(t)\\
	&\hspace{0.4cm}
	+\gamma
	\left[
	4(\tau+\sigma^\infty)
	+
	\frac{6\tau(1+\sigma^\infty)}
	{2-\sigma^\infty}
	\right]
	(1+t)^{\gamma-1}
	\int_{t-\tau}^{t}\Lambda^\omega(s)\,\dd s\\
	&\hspace{0.4cm}
	+\frac{(2-\sigma^\infty)2^{-\beta}}{24}
	\tilde{C}_p
	(1+t)^{\gamma-\beta\alpha-(p-2)(\alpha-1)}.
\end{aligned}
\end{equation}

\vspace{0.1cm}

\noindent$\bullet$~\textbf{Step B.} To ensure that the coefficient of $\mathcal E(t)$ in \eqref{E-4-1}
is negative for all sufficiently large $t$, and that the last term
in \eqref{E-4-1} is integrable on $[1,\infty)$, we require
\begin{equation}\label{E-5-0}
\beta\alpha<1,\qquad 
\gamma-\beta\alpha-(p-2)(\alpha-1)<-1
\end{equation}

\vspace{0.2cm}

\noindent $\bullet$ Case 1: \((0<\beta<1)\) If $\gamma$ satisfies
\[
0<\gamma<
\frac{(p-2)(1-\beta)}{\beta}=\lim_{\alpha\uparrow1/\beta}
\left[
\beta\alpha+(p-2)(\alpha-1)-1
\right],
\]
there exists an admissible exponent $\alpha_\gamma$ satisfying \eqref{E-5-0}:
\[\alpha_\gamma=
\frac12
\left(
\frac{\gamma+p-1}{p-2+\beta}
+
\frac1\beta
\right).\]

\vspace{0.2cm}

\noindent $\bullet$ Case 2: \((\beta=0)\) For any $\gamma>0$, we have an admissible exponent $\alpha_\gamma$ for \eqref{E-5-0}
\[
\alpha_\gamma=1+\frac{\gamma+2}{p-2}.
\]

\vspace{0.2cm}

\noindent 
Since $\beta\alpha_\gamma<1$, there exists a sufficiently large time $T_\gamma\ge\tau$ such that
\[
2\gamma(1+t)^{\gamma-1}
\leq
\frac{(2-\sigma^\infty)2^{-2\beta}}{24}
(1+t)^{\gamma-\beta\alpha_\gamma},\qquad t\geq T_\gamma.
\]
Specifically, we choose \[
T_\gamma
:=
\max\left\{
\tau,\,
\left(
\frac{48\,2^{2\beta}\gamma}{2-\sigma^\infty}
\right)^{\frac{1}{1-\beta\alpha_\gamma}}
-1,
\frac{48\gamma\,2^{(\gamma-1)_+}\tau}
{2-\sigma^\infty}
\left[
4(\tau+\sigma^\infty)
+
\frac{6\tau(1+\sigma^\infty)}
{2-\sigma^\infty}
\right]
-1
\right\}.
\]
Hence the coefficient of $\mathcal E(t)$ in \eqref{E-4-1} is nonpositive for $t\geq T_\gamma$:
\[
\left[
	2\gamma(1+t)^{\gamma-1}
	-
	\frac{(2-\sigma^\infty)2^{-2\beta}}{24}
	(1+t)^{\gamma-\beta\alpha_\gamma}
	\right]\mathcal E(t)
	\leq0.
\]
Therefore, we return to \eqref{E-4-1}. We integrate \eqref{E-4-1} over $[T_\gamma,t]$ and discard the nonpositive term involving $\mathcal E(t)$ to obtain
\begin{equation}	\label{E-5}
\begin{aligned}
	(1+t)^\gamma\mathcal L(t+\tau)
	&\leq
	(1+T_\gamma)^\gamma\mathcal L(T_\gamma+\tau)-\frac{2-\sigma^\infty}{24}
	\int_{T_\gamma}^{t}
	(1+s)^\gamma\Lambda^\omega(s)\,\dd s
\\
&+
	\gamma
	\left[
	4(\tau+\sigma^\infty)+
	\frac{6\tau(1+\sigma^\infty)}
	{2-\sigma^\infty}
	\right]
	\int_{T_\gamma}^{t}
	(1+s)^{\gamma-1}
	\int_{s-\tau}^{s}
	\Lambda^\omega(r)\,\dd r\,\dd s
\\
	&+
	\frac{(2-\sigma^\infty)2^{-\beta}}{24}
	\tilde{C}_p
	\int_{T_\gamma}^{t}
	(1+s)^{
		\gamma-\beta\alpha_\gamma
		-(p-2)(\alpha_\gamma-1)
	}\,\dd s, \qquad t\ge T_\gamma.
\end{aligned}
\end{equation}
By Fubini's theorem, we obtain
\[
\begin{aligned}
	&\int_{T_\gamma}^{t}
	(1+s)^{\gamma-1}
	\int_{s-\tau}^{s}\Lambda^\omega(r)\dd r\dd s\le
	2^{(\gamma-1)_+}\tau
	\int_{T_\gamma-\tau}^{t}
	(1+s)^{\gamma-1}\Lambda^\omega(s)\dd s\\
	&\hspace{3cm}=
	2^{(\gamma-1)
		_+}\tau
	\int_{T_\gamma-\tau}^{T_\gamma}
	(1+s)^{\gamma-1}\Lambda^\omega(s)\dd s+
	2^{(\gamma-1)_+}\tau
	\int_{T_\gamma}^{t}
	(1+s)^{\gamma-1}\Lambda^\omega(s)\dd s\\
	&\hspace{3cm}=:I_{41}+I_{42},
\end{aligned}
\]
where $(\gamma-1)_+:=\max\{\gamma-1,0\}$, and the first inequality uses
\[
(1+s)^{\gamma-1}
\leq
2^{(\gamma-1)_+}(1+r)^{\gamma-1}
\]
for $r\leq s\leq r+\tau$ and $\tau\leq1$.

\vspace{0.2cm}

\noindent $\bullet$~(Estimate of $I_{41}$) By Lemma~\ref{L3.1}, Proposition~\ref{P2.1}, Jensen's inequality, we have
\begin{equation}\label{E-5-1}
\sup_{s\geq0}\Lambda^\omega(s)
\leq
2\sup_{s\geq0}\mathcal E(s)
\leq
C_p\bigl(1+M_{p,v}^0\bigr)^{2/p}.
\end{equation}
Therefore, we obtain
\begin{align*}
I_{41}
&\leq
2^{(\gamma-1)_+}\tau
C_p\bigl(1+M_{p,v}^0\bigr)^{2/p}
\int_{T_\gamma-\tau}^{T_\gamma}
(1+s)^{\gamma-1}\,\dd s \\
&=
\frac{2^{(\gamma-1)_+}\tau C_p}{\gamma}
\Bigl[
(1+T_\gamma)^\gamma-(1+T_\gamma-\tau)^\gamma
\Bigr]
\bigl(1+M_{p,v}^0\bigr)^{2/p}.
\end{align*}

\vspace{0.2cm}

\noindent $\bullet$~(Estimate of $I_{42}$)
By the choice of $T_\gamma$, we have
\[
\gamma
\left[
4(\tau+\sigma^\infty)
+
\frac{6\tau(1+\sigma^\infty)}
{2-\sigma^\infty}
\right]
2^{(\gamma-1)_+}\tau
(1+s)^{\gamma-1}
\leq
\frac{2-\sigma^\infty}{48}(1+s)^\gamma, \qquad s\geq T_\gamma.
\]
Hence, we have
\[
\begin{aligned}
	&\gamma
	\left[
	4(\tau+\sigma^\infty)
	+
	\frac{6\tau(1+\sigma^\infty)}
	{2-\sigma^\infty}
	\right]
	I_{42} \leq
	\frac{2-\sigma^\infty}{48}
	\int_{T_\gamma}^{t}
	(1+s)^\gamma\Lambda^\omega(s)\,\dd s.
\end{aligned}
\]

\vspace{0.2cm}

\noindent We go back to \eqref{E-5}. The term involving $I_{42}$ can be absorbed into the
dissipation term. So with the estimate of $I_{41}$, we obtain
\begin{equation}\label{E-6}
	\begin{aligned}
		(1+t)^\gamma\mathcal L(t+\tau)
		&\leq
		(1+T_\gamma)^\gamma\mathcal L(T_\gamma+\tau)
		-\frac{2-\sigma^\infty}{48}
		\int_{T_\gamma}^{t}
		(1+s)^\gamma\Lambda^\omega(s)\,\dd s
		\\
		&+
		2^{(\gamma-1)_+}\tau C_p
		\left[
		4(\tau+\sigma^\infty)
		+
		\frac{6\tau(1+\sigma^\infty)}
		{2-\sigma^\infty}
		\right]
		\\
		&\times
		\Bigl[
		(1+T_\gamma)^\gamma
		-(1+T_\gamma-\tau)^\gamma
		\Bigr]
		\bigl(1+M_{p,v}^0\bigr)^{2/p}
		\\
		&+
		\frac{(2-\sigma^\infty)2^{-\beta}}{24}
		\tilde C_p
		\int_{T_\gamma}^{t}
		(1+s)^{
			\gamma-\beta\alpha_\gamma
			-(p-2)(\alpha_\gamma-1)
		}\,\dd s,
		\qquad t\geq T_\gamma.
	\end{aligned}
\end{equation}
By \eqref{E-5-0}, we have
\begin{equation}\label{E-7-0}
		\int_{T_\gamma}^{\infty}
		(1+s)^{\gamma-\beta\alpha_\gamma
			-(p-2)(\alpha_\gamma-1)}\,\dd s =
		\frac{
			(1+T_\gamma)^{
				-\left[
				\beta\alpha_\gamma
				+(p-2)(\alpha_\gamma-1)
				-\gamma-1
				\right]
			}
		}{
			\beta\alpha_\gamma
			+(p-2)(\alpha_\gamma-1)
			-\gamma-1
		}.
\end{equation}
Moreover, by \eqref{E-4} and \eqref{E-5-1} we obtain
\[
\begin{aligned}
\mathcal L(T_\gamma+\tau)
&\leq
C_p\bigl(1+M_{p,v}^0\bigr)^{2/p}
+
\tau
\left[
4(\tau+\sigma^\infty)
+
\frac{6\tau(1+\sigma^\infty)}
{2-\sigma^\infty}
\right]
C_p\bigl(1+M_{p,v}^0\bigr)^{2/p} \\
&=
C_p
\left[
1
+
4\tau(\tau+\sigma^\infty)
+
\frac{6\tau^2(1+\sigma^\infty)}
{2-\sigma^\infty}
\right]
\bigl(1+M_{p,v}^0\bigr)^{2/p}.
\end{aligned}
\]
We combine the above estimate with \eqref{E-6} and \eqref{E-7-0} to obtain
\[
(1+t)^\gamma\mathcal L(t+\tau)
\leq C_\gamma,
\qquad t\geq T_\gamma.
\]
where $C_\gamma$ is defined by
\[
\begin{aligned}
	C_\gamma
	:=\;&
	C_p(1+T_\gamma)^\gamma
	\left[
	1
	+
	4\tau(\tau+\sigma^\infty)
	+
	\frac{6\tau^2(1+\sigma^\infty)}
	{2-\sigma^\infty}
	\right]
	\bigl(1+M_{p,v}^0\bigr)^{2/p}
	\\
	&+
	2^{(\gamma-1)_+}\tau C_p
	\left[
	4(\tau+\sigma^\infty)
	+
	\frac{6\tau(1+\sigma^\infty)}
	{2-\sigma^\infty}
	\right]
	\\
	&\times
	\Bigl[
	(1+T_\gamma)^\gamma
	-
	(1+T_\gamma-\tau)^\gamma
	\Bigr]
	\bigl(1+M_{p,v}^0\bigr)^{2/p}
	\\
	&+
	\frac{(2-\sigma^\infty)2^{-\beta}\tilde C_p}{24}
	\frac{
		(1+T_\gamma)^{
			-\left[
			\beta\alpha_\gamma
			+(p-2)(\alpha_\gamma-1)
			-\gamma-1
			\right]
		}
	}{
		\beta\alpha_\gamma
		+(p-2)(\alpha_\gamma-1)
		-\gamma-1
	}.
\end{aligned}
\]
With \eqref{C-11-0}, we have
\[
\mathcal E(t+\tau)\leq\mathcal L(t+\tau)
\leq
C_\gamma(1+t)^{-\gamma},
\qquad t\geq T_\gamma.
\]
Since $\tau\le 1$, for $t\geq T_\gamma+\tau$, we get\[
\mathcal E(t)
\leq
C_\gamma(1+t-\tau)^{-\gamma}
\leq
2^\gamma C_\gamma(1+t)^{-\gamma}.
\]
On the remaining bounded interval $[0,T_\gamma+\tau]$,
Proposition~\ref{P2.1} and Jensen's inequality yield a uniform bound.
Hence, enlarging $C_\gamma$ once more, we obtain
\[
\mathcal E(t)
\leq
C_\gamma(1+t)^{-\gamma},
\qquad t\geq0.
\]
This completes the proof of \eqref{E-0}.

\subsection{Spatial cohesion}\label{sec:6.2} 
Suppose that the same setting as in Section~\ref{sec:6.1} holds.
Recall from \(({\mathcal F}_A1)\) that
\[
p>\frac{2}{1-\beta}.
\]
Now, we claim that there exists a constant
\(C>0\), independent of \(N\), such that
\begin{equation}   \label{E-8}
	\sup_{0 \leq t <\infty}
	\E\left[
	\frac1N\sum_{i=1}^N
	\norm{\bx^i_t- \bar{\bx}_t}^2
	\right]
	\leq C.
\end{equation}
{\it Proof of \eqref{E-8}}:~If \(0<\beta<1\), then
\[
p>\frac{2}{1-\beta}\Longrightarrow
\frac{(p-2)(1-\beta)}{\beta}>2.
\]
Hence, we choose
\begin{equation*}
	2<\gamma_0<
	\frac{(p-2)(1-\beta)}{\beta}.
\end{equation*}
When \(\beta=0\), we choose any \(\gamma_0>2\).
Then, we use \eqref{E-0} to find 
\begin{equation}
	\mathcal E(t)
	\leq
	C_{\gamma_0}(1+t)^{-\gamma_0},
	\qquad t\geq0.
	\label{E-9}
\end{equation}
Since
\[
\bx^i_t- \bar{\bx}_t
= \bx_0^i -\bar{\bx}_0
+
\int_0^t
\bigl( \bv^i_s- \bar{\bv}_s \bigr)\,\dd s,
\]
we use Minkowski's inequality in
\(L^2(\Omega;(\R^d)^N,\norm{\cdot}_N)\) to get 
\begin{align}
	\begin{aligned} \label{E-9-1}
		\left(
		\E\left[
		\frac1N\sum_{i=1}^N
		\norm{\bx^i_t- \bar{\bx}_t}^2
		\right]
		\right)^{1/2} &\leq
		\left(
		\E\left[
		\frac1N\sum_{i=1}^N
		\norm{\bx^i_0- \bar{\bx}_0}^2
		\right]
		\right)^{1/2}
		+
		\int_0^t
		\mathcal E(s)^{1/2}\,\dd s \\
		&=: I_{51} + I_{52}.
	\end{aligned}
\end{align}
\noindent $\bullet$~Case 1 (Estimate of $I_{51}$): Note that the term $I_{51}$ is bounded by the estimates in \eqref{B-4-1}:
\begin{equation} \label{E-9-2}
		I_{51}
		\leq
		\left(
		\E\left[
		\sup_{s\in[-\tau,0]}
		\frac1N\sum_{i=1}^N
		\|\bphi^i_s-\bar{\bphi}_s\|^p
		\right]
		\right)^{1/p}
		\leq
		\bigl(M^0_{p,x}\bigr)^{1/p}.
\end{equation}
\noindent $\bullet$~Case 2 (Estimate of $I_{52}$):~We use $\gamma_0>2$ and \eqref{E-9} to get 
\begin{equation} \label{E-9-3}
	I_{52} \leq
	C_{\gamma_0}^{1/2}\int_0^\infty
	(1+s)^{-\gamma_0/2}\,\dd s
	<\infty.
\end{equation}
Finally, we combine \eqref{E-9-1}, \eqref{E-9-2} and \eqref{E-9-3} to get the desired spatial cohesion estimate. Note that the resulting bound is uniform in \(t\) and \(N\).

\vspace{0.3cm}

\section{Almost-sure stochastic flocking}\label{sec:7}
\setcounter{equation}{0}
In this section, we study almost-sure flocking for \eqref{B-1}. \newline

Suppose the framework $({\mathcal F}_A1) - ({\mathcal F}_A3)$ holds, and let  $(\bX_t, \bV_t)$ be a strong solution process to \eqref{B-1}. Then, there exist
$ \sigma_* = \sigma_*(p)>0,~~\tau_* = \tau_*(p)>0,$ 
such that whenever
\[
0\leq\sigma^\infty\leq\sigma_*,
\qquad
0\leq\tau\leq\tau_*,
\]
for every fixed \(N\), we have
\begin{equation}\label{F-0}
	\lim_{t\to\infty}
	\max_{1\leq i,j\leq N}
	\|\bv^i_t-\bv^j_t\|
	=0,
	\qquad
	\sup_{t\geq0}
	\frac1N\sum_{i=1}^N
	\|\bx^i_t-\bar{\bx}_t\|^2
	<\infty,
	\qquad
	\text{a.s.}
\end{equation}
In the following two subsections, we derive estimates in \eqref{F-0} one by one. \newline

\subsection{Spatial cohesion} \label{sec:7.1}
Choose \(\gamma_0>2\) as in Section~\ref{sec:6.2}. Then, by \eqref{E-9}, we have
\[
\mathcal E(t)
\leq
C_{\gamma_0}(1+t)^{-\gamma_0},
\qquad t\geq0.
\]
For the spatial bound, we use Jensen's inequality and Tonelli's
theorem to get 
\[
\begin{aligned}
&\E\left[
\int_0^\infty
\left(
\frac1N\sum_{i=1}^N
\norm{\bv^i_s- \bar{\bv}_s}^2
\right)^{1/2}
\dd s
\right]\leq
\int_0^\infty
\mathcal E(s)^{1/2}\,\dd s\leq
C_{\gamma_0}^{1/2}\int_0^\infty
(1+s)^{-\gamma_0/2}\,\dd s
<\infty.
\end{aligned}
\]
Hence, we have
\[
\int_0^\infty
\left(
\frac1N\sum_{i=1}^N
\norm{\bv^i_s- \bar{\bv}_s}^2
\right)^{1/2}
\dd s
<\infty
\qquad\text{a.s.}
\]
Now, we use
\[
\bx^i_t- \bar{\bx}_t
=
\bx^i_0 - \bar{\bx}_0
+
\int_0^t
\bigl( \bv^i_s -\bar{\bv}_s \bigr)\,\dd s
\]
and Minkowski's inequality in the empirical norm to get 
\begin{align*} \sup_{t\geq0}
&\left(
\frac1N\sum_{i=1}^N
\norm{\bx^i_t- \bar{\bx}_t}^2
\right)^{1/2}\\
&\qquad\leq
\left(
\frac1N\sum_{i=1}^N
\norm{\bx^i_0 - {\bar \bx}_0}^2
\right)^{1/2}+
\int_0^\infty
\left(
\frac1N\sum_{i=1}^N
\norm{\bv^i_s- \bar{\bv}_s}^2
\right)^{1/2}
\dd s<\infty\qquad\mbox{a.s.}
\end{align*}
This verifies the second estimate in \eqref{F-0}. 

\subsection{Velocity alignment} \label{sec:7.2}
It follows from 
\eqref{C-13} that 
\begin{equation*}
\int_0^\infty \Lambda^\omega(s)\,\dd s<\infty.
\end{equation*}
Now, we use \eqref{E-9} and \(\gamma_0>1\) to get 
\[
\begin{aligned}
\sum_{n=\lceil\tau\rceil}^{\infty}
\E\left[
\frac1N\sum_{i=1}^N
\| \bv^i_n- \bar{\bv}_n \|^2
\right]
&\leq
C_{\gamma_0}\sum_{n=\lceil\tau\rceil}^{\infty}
(1+n)^{-\gamma_0}<\infty.
\end{aligned}
\]
Then, Tonelli's theorem yields
\[
\sum_{n=\lceil\tau\rceil}^{\infty}
\frac1N\sum_{i=1}^N
\| \bv^i_n -\bar{\bv}_n \|^2
<\infty
\qquad\text{a.s.}
\]
In particular, we have
\begin{equation}
\frac1N\sum_{i=1}^N
\norm{\bv^i_n -\bar{\bv}_n}^2
\longrightarrow0 \quad \text{as} \quad n\to \infty 
\qquad\text{a.s.}
\label{E-13}
\end{equation}
For
\(n\geq\lceil\tau\rceil\), we define
\[
Z_n
:=
\sup_{n\leq t\leq n+1}
\left(
\frac1N\sum_{i=1}^N
\norm{
\bigl(\bv^i_t- \bar{\bv}_t \bigr)
-
\bigl(\bv^i_n- \bar{\bv}_n \bigr)
}^2
\right)^{1/2}.
\]
For \(t\in[n,n+1]\), we set
\[
M^{n,i}_t
:=
\int_n^t
\sqrt{\sigma^i}\,
\Diag(\bF^i_{s-\tau})\,\dd {\bf B}^i_s, \qquad
\overline{M}^{n}_t
:=
\frac1N\sum_{\ell=1}^N M_t^{n,\ell}.
\]
We integrate the centered velocity equation over \([n,t]\) to obtain
\[
\begin{aligned}
&\bigl(\bv^i_t- \bar{\bv}_t \bigr)
-
\bigl(\bv^i_n- \bar{\bv}_n \bigr)=  \int_n^t \bF^i_{s-\tau} \,\dd s
+ M^{n,i}_t  -
\overline M^{n}_t.
\end{aligned}
\]
Consequently, we have
\begin{align}
\begin{aligned} \label{E-14}
	Z_n^2
	&\leq
	2\sup_{n\leq t\leq n+1}
	\frac1N\sum_{i=1}^N
	\left\|
	\int_n^t \bF^i_{s-\tau}\,\dd s
	\right\|^2+
	2\sup_{n\leq t\leq n+1}
	\frac1N\sum_{i=1}^N
	\norm{
		M_t^{n,i} -\overline M^{n}_t
	}^2 \\
	&=: 2(I_{61} + I_{62}).
\end{aligned}
\end{align}
\noindent $\bullet$~Case 1 (Estimate of $ I_{61}$):~We use the Cauchy--Schwarz inequality and \eqref{C-4-2} to get 
\begin{equation*}
\begin{aligned}
	&\E\left[
	\sup_{n\leq t\leq n+1}
	\frac1N\sum_{i=1}^N
	\left\|
	\int_n^t \bF^i_{s-\tau} \,\dd s
	\right\|^2
	\right]\leq
	\int_n^{n+1}
	\E\left[
	\frac1N\sum_{i=1}^N
	\norm{\bF^i_{s-\tau}}^2
	\right]\dd s\leq
	\frac12
	\int_{n-\tau}^{n+1-\tau}\!\!
	\Lambda^\omega(r)\,\dd r.
\end{aligned}
\end{equation*}

\vspace{0.2cm}

\noindent $\bullet$~Case 2 (Estimate of $ I_{62}$):~Note that the subtraction of the empirical mean is an
orthogonal contraction in the empirical \(\ell^2\)-norm. Hence, we have
\[
\begin{aligned}
&\sup_{n\leq t\leq n+1}
\frac1N\sum_{i=1}^N
\norm{
	M^{n,i}_t- {\overline M}^{n}_t
}^2\leq
\frac1N\sum_{i=1}^N
\sup_{n\leq t\leq n+1}
\norm{M^{n,i}_t}^2.
\end{aligned}
\]
Doob's \(L^2\) inequality, It\^o's isometry, and \eqref{C-4-2} yield
\begin{equation}
\begin{aligned}
	&\E\left[
	\sup_{n\leq t\leq n+1}
	\frac1N\sum_{i=1}^N
	\norm{
		M^{n,i}_t -{\overline M}^{n}_t
	}^2
	\right]\\
	&\qquad\leq
	\frac4N\sum_{i=1}^N
	\E\left[
	\int_n^{n+1}
	\sigma^i
	\norm{\Diag(\bF^i_{s-\tau})}_{\mathrm{F}}^2
	\,\dd s
	\right]\leq
	2\sigma^{\infty}
	\int_{n-\tau}^{n+1-\tau}
	\Lambda^\omega(s)\,\dd s.
\end{aligned}
\label{E-15}
\end{equation}
We combine \eqref{E-14}--\eqref{E-15} to find
\[
\E[Z_n^2]
\leq
C(1+\sigma^{\infty})
\int_{n-\tau}^{n+1-\tau}
\Lambda^\omega(s)\,\dd s.
\]
Therefore, it follows from \eqref{C-13} that 
\[
\E\left[
\sum_{n=\lceil\tau\rceil}^{\infty}Z_n^2
\right]
<\infty.
\]
It follows that
\[
\sum_{n=\lceil\tau\rceil}^{\infty}Z_n^2
<\infty
\qquad\text{a.s.}
\]
Hence, we have
\[  Z_n~\longrightarrow~0
\qquad\text{a.s.} \]
For \(t\in[n,n+1]\), it holds
\[
\begin{aligned}
&\left(
\frac1N\sum_{i=1}^N
\norm{\bv^i_t- \bar{\bv}_t}^2
\right)^{1/2}\leq
\left(
\frac1N\sum_{i=1}^N
\norm{\bv^i_n - \bar{\bv}_n}^2
\right)^{1/2}
+
Z_n.
\end{aligned}
\]
This and \eqref{E-13} yield
\[
\sup_{n\leq t\leq n+1}
\left(
\frac1N\sum_{i=1}^N
\norm{\bv^i_t- \bar{\bv}_t}^2
\right)^{1/2}
\longrightarrow0
\qquad\text{a.s.}
\]
Finally, for fixed \(N\), we have
\[
\max_{1\leq i,j\leq N}
\norm{\bv^i_t-\bv^j_t}
\leq
2\sqrt N
\left(
\frac1N\sum_{i=1}^N
\norm{\bv^i_t- \bar{\bv}_t}^2
\right)^{1/2},
\]
which proves $\eqref{F-0}_1$. Similarly, we have
\[
\max_{1\leq i,j\leq N}
\norm{\bx^i_t-\bx^j_t}
\leq
2\sqrt N
\left(
\frac1N\sum_{i=1}^N
\norm{\bx^i_t- \bar{\bx}_t}^2
\right)^{1/2}.
\]
Together with the spatial bound obtained in Section~\ref{sec:7.1},
we have
\[
\sup_{t\geq0}
\max_{1\leq i,j\leq N}
\norm{\bx^i_t-\bx^j_t}
<\infty
\qquad\text{a.s.}
\]
This completes the proof of Theorem~\ref{T2.4}.

\section{Numerical experiments}
\label{sec:8}
In this section, we provide several numerical simulations to illustrate our theoretical results. More precisely, we illustrate the deterioration of flocking behavior
when the time delay or noise intensity becomes large.\newline 
\begin{figure}[H]
\centering
\includegraphics[width=0.68\textwidth]{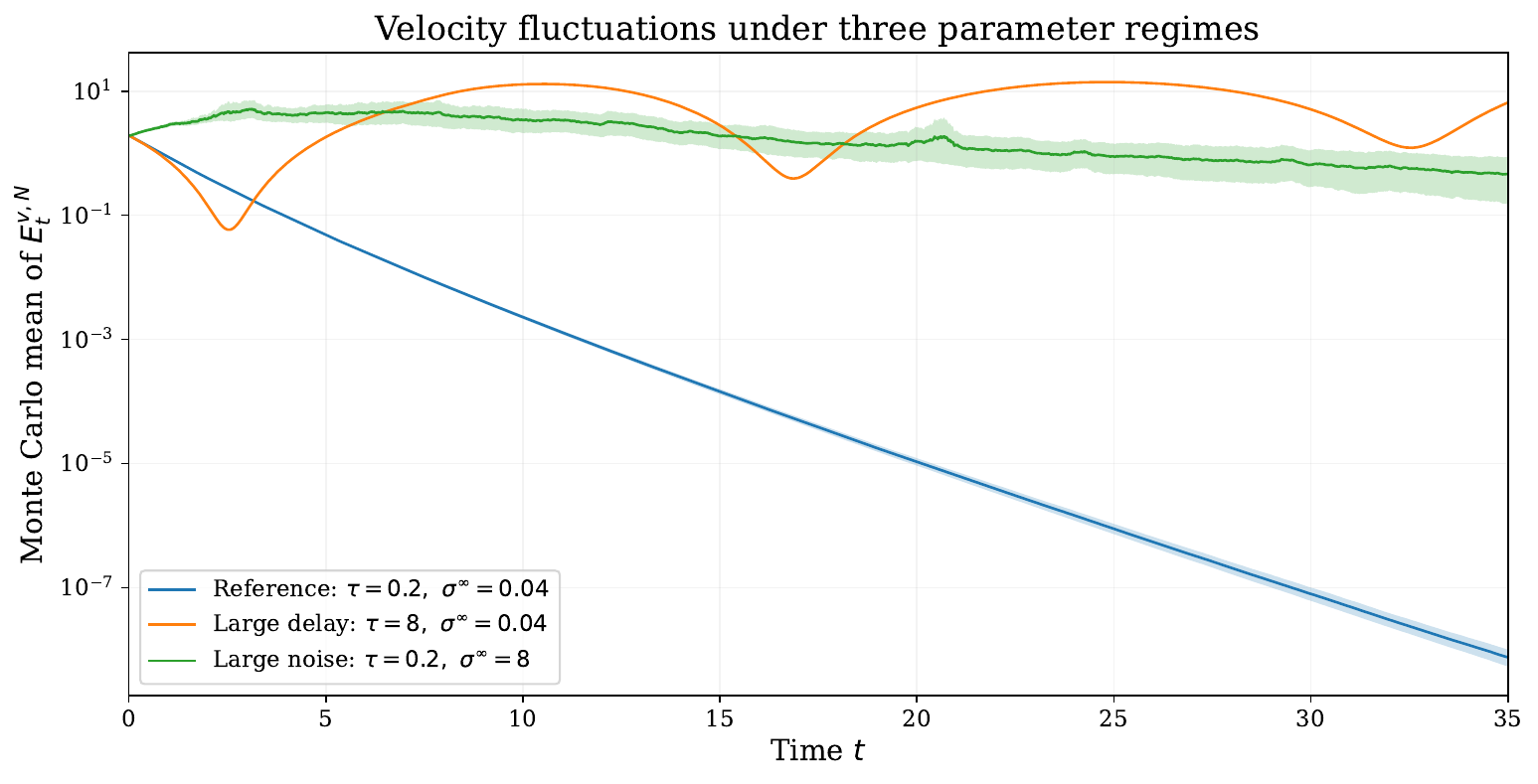}
\caption{
Velocity fluctuation under the three parameter regimes.
The solid curves represent the Monte Carlo means of
\(E_t^{v,N}\) computed from \(M=300\) sample paths, while the
shaded regions denote pointwise percentile-bootstrap
\(95\%\) confidence intervals based on \(2000\) bootstrap
resamples. 
}
\label{Fig-1}
\end{figure}

\begin{figure}[t]
\centering
\includegraphics[width=0.68\textwidth]{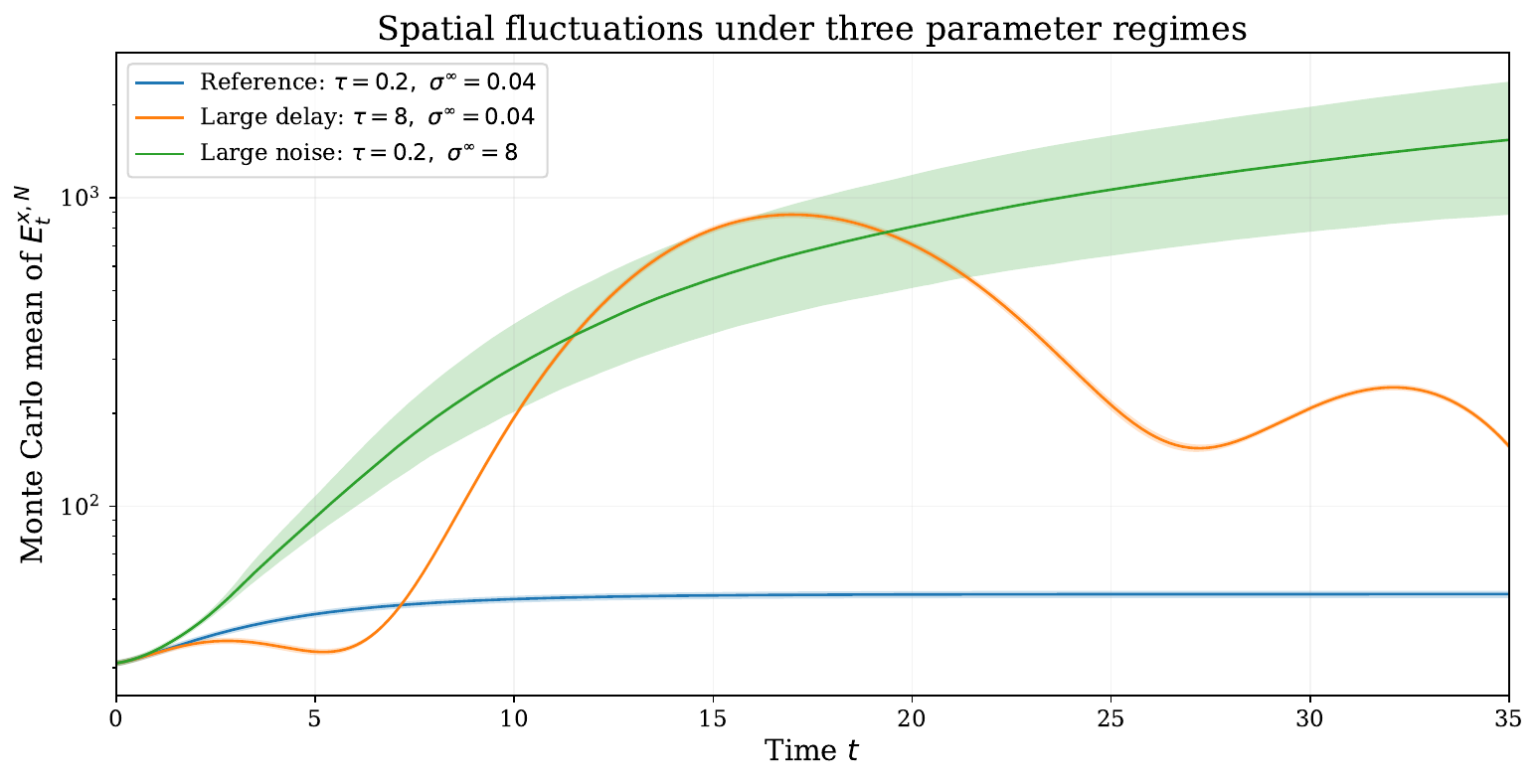}
\caption{
Spatial fluctuation under the three parameter regimes.
The solid curves represent the Monte Carlo means of
\(E_t^{x,N}\) computed from \(M=300\) sample paths, while the
shaded regions denote pointwise percentile-bootstrap
\(95\%\) confidence intervals based on \(2000\) bootstrap
resamples. 
}
\label{Fig-2}
\end{figure}
For \(s\in[-\tau,0]\), we prescribe the constant history
\[
\bphi^i_s \equiv \bphi^i_0,
\quad
\bpsi^i_s \equiv \bpsi^i_0.
\]
For homogeneous noise intensities
\[
\sigma^i\equiv\sigma^\infty,\qquad i\in[N],
\]
the Euler--Maruyama scheme reads as
\[
\begin{cases}
\displaystyle \bx^{i}_{n+1} =
\bx^{i}_{n} + \bv^{i}_{n} \Delta t, \\
\displaystyle \bv^{i}_{n+1} =
\bv^{i}_{n} +\bF^{i}_{n} \Delta t
+ \sqrt{\sigma^{\infty}}\,\Diag(\bF^{i}_{n})\Delta {\bf B}^{i}_{n},
\end{cases}	
\]
with
\[
m_\tau:=\frac{\tau}{\Delta t}\in\mathbb N,
\qquad
\Delta\mathbf{B}_n^i\sim\mathcal N(0,\Delta t I_d),
\]
\[
\bF_n^i
=
\frac1N\sum_{j=1}^N
\omega\left(
\norm{\bx_{n-m_\tau}^i-\bx_{n-m_\tau}^j}
\right)
\left(
\bv_{n-m_\tau}^j-\bv_{n-m_\tau}^i
\right),
\]
where the Brownian increments
\(\{\Delta\mathbf{B}_n^i\}_{i\in[N],\,n\geq0}\)
are independent. For the discrete history, we set
\[
\bx_n^i=\bx_0^i,
\qquad
\bv_n^i=\bv_0^i,
\qquad
-m_\tau\leq n\leq0,\quad i\in[N].
\]
In accordance with Remark~\ref{R2.2}(4), we take
\[
N=30,\qquad d=2,\qquad
\omega(r)=(1+r)^{-1/2},\qquad
\Delta t=10^{-2},\qquad T=35,
\]
and simulate \(M=300\) paths. For each path,
the initial data are chosen as follows:
\[
\bx^i_0=4\bx^i,
\quad
\bv^i_0
=
\eta^i
+
0.15\frac{\bx^i_0}{\sqrt{1+\norm{\bx^i_0}^2}},
\]
where \(\bx^i\) and \(\eta^i\) are independent standard Gaussian
vectors. The same initial realization and underlying standard Gaussian
increments are reused across the three regimes:
\[
\begin{array}{c|cc}
\text{regime} & \tau & \sigma^\infty \\ \hline
\text{reference}   & 0.2 & 0.04 \\
\text{large delay} & 8   & 0.04 \\
\text{large noise} & 0.2 & 8
\end{array}
\]
We consider the empirical fluctuations
\[
E^{v, N}_t
:=
\frac1N\sum_{i=1}^N
\norm{\bv^i_t- \bar{\bv}_t}^2,
\quad
E^{x, N}_t
:=
\frac1N\sum_{i=1}^N
\norm{\bx^i_t- \bar{\bx}_t}^2.
\]
The displayed error bands are pointwise percentile-bootstrap
\(95\%\) confidence intervals for the arithmetic mean, based on
\(2000\) resamples. \newline

In Figure~\ref{Fig-1}, we can see that the reference mean
decays to approximately \(10^{-8}\). In contrast, the large-delay
curve displays pronounced oscillations over the simulated interval
and remains of order one at \(T=35\). In the large-noise regime, the velocity fluctuation does not exhibit
sustained decay over the simulated interval. Figure~\ref{Fig-2} shows that the reference spatial
moment levels off over the simulated interval. The large-delay mean reaches a transient peak near \(882.40\)
and then decreases to about \(156.61\), whereas the large-noise regime exhibits substantially
stronger spatial spreading. These numerical simulations illustrate that large time-delay or noise
intensity may significantly deteriorate flocking behavior.

\vspace{0.5cm}

\section{Conclusion}
\setcounter{equation}{0} \label{sec:9}
In this paper, we established polynomial mean-square velocity alignment and uniform-in-time centered spatial confinement for a delayed stochastic CS system with an algebraically decaying communication kernel. Both estimates are uniform in the particle number \(N\) under suitable higher-order empirical moment assumptions on the initial history. We further proved almost-sure flocking for every fixed \(N\). The proof combines a memory-corrected Lyapunov functional, uniform
high-order empirical velocity moments and a method of time-varying effective region. This allows us to recover sufficient alignment dissipation without assuming an a priori bound on the spatial diameter or a uniform spectral gap. Our analysis does not cover the critical exponent $\beta=1$.  Moreover, determining sharp joint
delay--noise thresholds and extending the analysis to the corresponding
mean-field and kinetic equations remain open questions. We leave these interesting problems for future study.

\vspace{0.2cm}
\section*{Conflict of interest statement}
The authors declare no conflicts of interest.

\section*{Data availability statement}
The data supporting the findings of this study are available from the corresponding author upon reasonable request.

\section*{Declaration of generative AI use}
During the preparation of this work, the last author used ChatGPT (OpenAI) to assist with language editing and consistency checks. The authors reviewed and verified all content and take full
responsibility for the content of the published article.

\section*{Acknowledgements}

The work of S.-Y. Ha is supported by National Research Foundation (NRF) grant funded by the Korea government (MSIT) (RS-2025-00514472). The work of X. Wang is supported by the Natural Science Foundation of China (grants 123B2003), the China Postdoctoral Science Foundation (grants 2025M774290), and Heilongjiang Province Postdoctoral Funding (grants LBH-Z24167). The work of C. Wu was supported by the National Natural Science Foundation of China (Grant No.~12526501) and the China Postdoctoral Science Foundation
(Grant No. ~GZB20260744 and No.~2025M783063). This work was carried out while X. Wang and F. Zeng were visiting
the Department of Mathematical Sciences, Seoul National University,
as visiting professors during the 2026--2027 academic year.

\newpage

\appendix

\renewcommand{\theequation}{\thesection.\arabic{equation}}

\section{Auxiliary technical estimates} \label{App-A}
\setcounter{equation}{0}
In this appendix, we provide the proof of Lemma \ref{L2.3} and proofs
for several lemmas in Section \ref{sec:4} which have been used in the
proof of our main results.

\subsection{Proof of Lemma~\ref{L2.3}}
\label{App-A-0}
	For \(p=2\), all the assertions are obvious. Below, we consider
	\(p>2\). We set
	\[
	A_p(\bz):=\norm{\bz}^{p-2}\bz.
	\]
	Then \(A_p\in C^1(\R^d;\R^d)\). For \(\bz\neq0\) and
	\(\bh\in\R^d\), by direct calculation, one has
	\[
	DA_p(\bz)\bh
	=
	\norm{\bz}^{p-2}\bh
	+
	(p-2)\norm{\bz}^{p-4}(\bz\cdot\bh)\bz,
	\]
	and \(DA_p(0)=0\). Hence,
	\[
	\begin{aligned}
		DA_p(\bz)\bh\cdot\bh
		&=
		\norm{\bz}^{p-2}\norm{\bh}^2
		+
		(p-2)\norm{\bz}^{p-4}(\bz\cdot\bh)^2
\geq
		\norm{\bz}^{p-2}\norm{\bh}^2,
	\end{aligned}
	\]
	and
	\[
	\norm{DA_p(\bz)\bh}
	\leq
	(p-1)\norm{\bz}^{p-2}\norm{\bh}.
	\]
	
	\medskip
	
	\noindent
	(i)~By the fundamental theorem of calculus, we have
	\[
	\begin{aligned}
		&\bigl(A_p(\bx)-A_p(\by)\bigr)\cdot(\bx-\by)
		\\
		&\quad=
		\int_0^1
		DA_p\bigl(\by+\vartheta(\bx-\by)\bigr)
		(\bx-\by)\cdot(\bx-\by)\,\dd\vartheta
		\\
		&\quad\geq
		\norm{\bx-\by}^2
		\int_0^1
		\norm{\by+\vartheta(\bx-\by)}^{p-2}
		\,\dd\vartheta.
	\end{aligned}
	\]
	We set
	\[
	M:=\max\{\norm{\bx},\norm{\by}\}.
	\]
	If \(M=0\), then the desired estimate is trivial. For \(M>0\),
	without loss of generality, we assume
	\[
	M=\norm{\by}.
	\]
	Since
	\[
	\norm{\bx-\by}
	\leq
	\norm{\bx}+\norm{\by}
	\leq
	2M,
	\]
	for \(0\leq\vartheta\leq1/4\), we obtain
	\[
	\begin{aligned}
		\norm{\by+\vartheta(\bx-\by)}
		&\geq
		\norm{\by}
		-
		\vartheta\norm{\bx-\by}
		\geq
		M-\frac14(2M)
		=
		\frac{M}{2}.
	\end{aligned}
	\]
	It follows that
	\[
	\begin{aligned}
		\int_0^1
		\norm{\by+\vartheta(\bx-\by)}^{p-2}
		\,\dd\vartheta
		&\geq
		\int_0^{1/4}
		\norm{\by+\vartheta(\bx-\by)}^{p-2}
		\,\dd\vartheta
	\geq
		\frac14
		\left(\frac{M}{2}\right)^{p-2}.
	\end{aligned}
	\]
	On the other hand, since
	\[
	M
	\geq
	\frac{\norm{\bx}+\norm{\by}}{2},
	\]
we have
	\[
	\int_0^1
	\norm{\by+\vartheta(\bx-\by)}^{p-2}
	\,\dd\vartheta
	\geq
	c_p
	\bigl(\norm{\bx}+\norm{\by}\bigr)^{p-2},
	\]
	where \(c_p>0\) depends only on \(p\). Combining the above estimates,
	we obtain
	\[
	\bigl(A_p(\bx)-A_p(\by)\bigr)\cdot(\bx-\by)
	\geq
	c_p
	\bigl(\norm{\bx}+\norm{\by}\bigr)^{p-2}
	\norm{\bx-\by}^2.
	\]
	This proves (i).
	
	\vspace{0.2cm}
	
	\noindent
	(ii)~It follows from (i) and the triangle inequality
	\(
	\norm{\bx}+\norm{\by}
	\geq
	\norm{\bx-\by}
	\)
	that
	\[
	\bigl(A_p(\bx)-A_p(\by)\bigr)\cdot(\bx-\by)
	\geq
	c_p\norm{\bx-\by}^p.
	\]
	This proves (ii).
	
	\vspace{0.2cm}
	
	\noindent
	(iii)~By the fundamental theorem of calculus, we have
	\[
	A_p(\bx)-A_p(\by)
	=
	\int_0^1
	DA_p\bigl(\by+\vartheta(\bx-\by)\bigr)
	(\bx-\by)\,\dd\vartheta.
	\]
	Hence, by the estimate for \(DA_p\), one has
	\[
	\begin{aligned}
		\norm{A_p(\bx)-A_p(\by)}
		&\leq
		(p-1)\norm{\bx-\by}
		\int_0^1
		\norm{\by+\vartheta(\bx-\by)}^{p-2}
		\,\dd\vartheta.
	\end{aligned}
	\]
	Since 
	\[
	\norm{\by+\vartheta(\bx-\by)}
	\leq
	\norm{\by}+\norm{\bx-\by},\qquad \forall~0\leq\vartheta\leq1,
	\]
we get
	\[
	\begin{aligned}
		\int_0^1
		\norm{\by+\vartheta(\bx-\by)}^{p-2}
		\,\dd\vartheta
		&\leq
		\bigl(\norm{\by}+\norm{\bx-\by}\bigr)^{p-2}
		\leq
		C_p\left(
		\norm{\by}^{p-2}
		+
		\norm{\bx-\by}^{p-2}
		\right).
	\end{aligned}
	\]
	Combining the above estimates, we obtain
	\[
	\norm{A_p(\bx)-A_p(\by)}
	\leq
	C_p\left(
	\norm{\by}^{p-2}\norm{\bx-\by}
	+
	\norm{\bx-\by}^{p-1}
	\right).
	\]
	This verifies (iii).

\subsection{Proof of Lemma~\ref{L4.2}}
\label{App-A.1} 
In this subsection, we derive four estimates in \eqref{D-5} one by one. Recall that
\[
\bF^i
=
\frac{1}{N}\sum_{j=1}^N\omega^{ij}(\bz^j-\bz^i).
\]
By Jensen's inequality, with the missing mass assigned to the zero vector,
for any \(q\geq1\), we have
\begin{equation}\label{Ap-7}
	\norm{\bF^i}^q
	\leq
	\frac1N\sum_{j=1}^N
	\omega^{ij}\norm{\bz^i-\bz^j}^q.
\end{equation}

\noindent (i)~\textbf{(Derivation of $\eqref{D-5}_1$)}: Since\[
\bF^i
=
\frac{1}{N}\sum_{j=1}^N\omega^{ij}\bz^j
-
\left(
\frac{1}{N}\sum_{j=1}^N\omega^{ij}
\right)\bz^i,
\]
 we use \(\omega^{ij}\leq1\), the triangle inequality and
	Jensen's inequality, with the missing mass assigned to the zero vector,
	to find
	\[
	\begin{aligned}
		\norm{\bF^i}^p
		&\leq
		2^{p-1}
		\left(
		\norm{\bz^i}^p
		+
		\left\|
		\frac1N\sum_{j=1}^N\omega^{ij}\bz^j
		\right\|^p
		\right)
	\leq
		2^{p-1}
		\left(
		\norm{\bz^i}^p
		+
		\frac1N\sum_{j=1}^N\norm{\bz^j}^p
		\right).
	\end{aligned}
	\]
	We sum the above inequality over \(i\in[N]\) and divide it by \(N\)
	to obtain \eqref{D-5}$_1$ with \(C_p=2^p\).
	
	\vspace{0.2cm}
	
	\noindent (ii)~\textbf{(Derivation of $\eqref{D-5}_2$)}: For \eqref{Ap-7} with \(q=p\), we use Lemma~\ref{L2.3} (ii)
	and index exchange symmetry to find
	\[
	\begin{aligned}
		\frac1N\sum_{i=1}^N\norm{\bF^i}^p
		&\leq
		\frac{C_p}{N^2}\sum_{i,j=1}^N
		\omega^{ij}
		\left(
		\norm{\bz^i}^{p-2}\bz^i
		-
		\norm{\bz^j}^{p-2}\bz^j
		\right)
		\cdot(\bz^i-\bz^j)
		\\
		&\le
		\frac{C_p}{N^2}\sum_{i,j=1}^N
		\omega^{ij}
		\norm{\bz^i}^{p-2}\bz^i
		\cdot(\bz^i-\bz^j).
	\end{aligned}
	\]
	
	\vspace{0.2cm}
	
	 \noindent (iii) \textbf{(Derivation of $\eqref{D-5}_3$)}:~
	 For \eqref{Ap-7} with $q=2$, we use index exchange symmetry to see
	 \begin{align*}
	 	&\frac1N\sum_{i=1}^N
	 	\norm{\bz^i}^{p-2}\norm{\bF^i}^2
	 	\\
	 	&\hspace{1cm}\leq
	 	\frac1{N^2}\sum_{i,j=1}^N
	 	\omega^{ij}
	 	\norm{\bz^i}^{p-2}\norm{\bz^i-\bz^j}^2
	 	\\
	 	&\hspace{1cm}=
	 	\frac1{2N^2}\sum_{i,j=1}^N
	 	\omega^{ij}
	 	\left(
	 	\norm{\bz^i}^{p-2}
	 	+
	 	\norm{\bz^j}^{p-2}
	 	\right)
	 	\norm{\bz^i-\bz^j}^2
	 	\\
	 	&\hspace{1cm}\leq
	 	\frac{C_p}{N^2}\sum_{i,j=1}^N
	 	\omega^{ij}
	 	\left(
	 	\norm{\bz^i}^{p-2}\bz^i
	 	-
	 	\norm{\bz^j}^{p-2}\bz^j
	 	\right)
	 	\cdot(\bz^i-\bz^j)
	 	\\
	 	&\hspace{1cm}=
	 	\frac{2C_p}{N^2}\sum_{i,j=1}^N
	 	\omega^{ij}
	 	\norm{\bz^i}^{p-2}\bz^i
	 	\cdot(\bz^i-\bz^j),
	 \end{align*}
	where the second inequality follows from
	Lemma~\ref{L2.3}(i) and
	\[
	\norm{\bz^{i}}^{p-2}+\norm{\bz^j}^{p-2}
	\leq
	C_p(\norm{\bz^i}+\norm{\bz^j})^{p-2}.
	\]
	
	\vspace{0.2cm}
	
\noindent (iv)~\textbf{(Derivation of $\eqref{D-5}_4$)}:
For the nonnegativity, we use index exchange symmetry and
Lemma~\ref{L2.3}(ii) to find
\[
\begin{aligned}
	&\frac1{N^2}\!\!\sum_{i,j=1}^N
	\omega^{ij}\!
	\norm{\bz^i}^{p-2}\bz^i\cdot(\bz^i-\bz^j)
=
	\!\frac1{2N^2}\sum_{i,j=1}^N
	\omega^{ij}\!
	\left(
	\norm{\bz^i}^{p-2}\bz^i
	-
	\norm{\bz^j}^{p-2}\bz^j
	\right)\!
	\cdot\!(\bz^i-\bz^j)
	\geq0.
\end{aligned}
\]
For the upper bound, by the definition of \(\bF^i\), we use
the Cauchy--Schwarz inequality, H\"older's inequality and
\eqref{D-5}$_1$ to find
\[
\begin{aligned}
	&\frac1{N^2}\sum_{i,j=1}^N
	\omega^{ij}
	\norm{\bz^i}^{p-2}\bz^i\cdot(\bz^i-\bz^j)
=
	-\frac1N\sum_{i=1}^N
	\norm{\bz^i}^{p-2}\bz^i\cdot\bF^i
	\\
	&\hspace{1cm}\leq
	\frac1N\sum_{i=1}^N
	\norm{\bz^i}^{p-1}\norm{\bF^i}
\leq
	\left(
	\frac1N\sum_{i=1}^N\norm{\bz^i}^p
	\right)^{\frac{p-1}{p}}
	\left(
	\frac1N\sum_{i=1}^N\norm{\bF^i}^p
	\right)^{\frac1p}
\leq
	\frac{C_p}{N}\sum_{i=1}^N\norm{\bz^i}^p.
\end{aligned}
\]
This completes the proof.
\subsection{Proof of Lemma~\ref{L4.4}}
\label{App-A.2}
In this subsection, we derive two estimates \eqref{D-6} and \eqref{D-7} one by one. \newline

\noindent (i)~(Derivation of \eqref{D-6}): By Lemma~\ref{L4.2}(iv), we have
\[
0\leq
{\mathcal D}_p(t)
\leq
C_p\E\left[
\frac1N\sum_{i=1}^N
\|\hat{\bv}_t^i\|^p
\right].
\]
Hence, by Lemma~\ref{L4.3}, we obtain
\begin{equation}\label{i-4}
{\mathcal D}_p\in L^1_{\mathrm{loc}}([0,\infty)).
\end{equation}
Next, we claim that for every \(r \geq2\tau\) and
\(0\leq h\leq\tau\), the following relation holds:
\begin{equation}\label{i-5}
	\begin{aligned}
		&\E\left[
		\frac1N\sum_{i=1}^N
		\norm{\hat{\bv}^i_r-\hat{\bv}^i_{r-h}}^p
		\right]
		\leq
		C_p\left(
		h^{p-1}
		+
		(\sigma^{\infty})^{p/2}h^{p/2-1}
		\right)
		\int_{r-h-\tau}^{r-\tau}{\mathcal D}_p(s)\,\dd s.
	\end{aligned}
\end{equation}
{\it Derivation of \eqref{i-5}}:~We integrate \eqref{D-2} over \([r-h,r]\) to obtain
\begin{equation}\label{i-6}
	\begin{aligned}
		\hat{\bv}^i_r -\hat{\bv}^i_{r-h}
		&=
		\int_{r-h}^{r} \bF^i_{s-\tau} \,\dd s
		+ \int_{r-h}^{r}
		\sqrt{\sigma^i}\,\Diag(\bF^i_{s-\tau})\dd {\bf B}^i_s.
	\end{aligned}
\end{equation}
This implies 
\begin{align*}
	\begin{aligned}
		& \E \left[ \frac1N\sum_{i=1}^N \| \hat{\bv}^i_r -\hat{\bv}^i_{r-h} \|^p \right]  \\
		& \hspace{0.2cm} \leq C_p \left ( 
		\E \left[   \frac1N\sum_{i=1}^N 
		\Big \| \int_{r-h}^{r} \bF^i_{s-\tau} \,\dd s \Big\|^p \right] 
		+ 
		\E \left[  \frac1N\sum_{i=1}^N 
		\Big\| \int_{r-h}^{r}
		\sqrt{\sigma^i}\,\Diag(\bF^i_{s-\tau})\dd {\bf B}^i_s  
		\Big \|^p    \right]             
		\right ) \\
		& \hspace{0.2cm} =: C_p  \Big (I_{71} + I_{72} \Big ).
	\end{aligned}
\end{align*}	
Below, we estimate the terms \(I_{7i}\) one by one. \newline

\noindent $\bullet$~(Estimate of $I_{71}$):~We use H\"older's inequality in time to get 
\[
\begin{aligned}
	I_{71} 
	&\leq
	h^{p-1}
	\int_{r-h}^{r}
	\E\left[
	\frac1N\sum_{i=1}^N\norm{\bF^i_{s-\tau}}^p
	\right]\dd s
\leq
	C_ph^{p-1}
	\int_{r-h-\tau}^{r-\tau}
	{\mathcal D}_p(s)\,\dd s,
\end{aligned}
\]
where the last inequality follows from Lemma~\ref{L4.2}(ii). \newline

\noindent $\bullet$~(Estimate of $I_{72}$):~We use the Burkholder--Davis--Gundy inequality
\cite[Theorem~7.3, p.~40]{M} and
\(
\norm{\Diag(\bF^i_{s-\tau})}_{\mathrm{F}}^2
=
\norm{\bF^i_{s-\tau}}^2
\)
to obtain
\begin{align*}
I_{72}& =
\E\left[
\frac1N\sum_{i=1}^N
\left\|
\int_{r-h}^{r}
\sqrt{\sigma^i}\,
\Diag(\bF^i_{s-\tau})\dd {\bf B}^i_s
\right\|^p
\right]\\
&\leq
C_p (\sigma^{\infty})^{p/2}
\E\left[
\frac1N\sum_{i=1}^N
\left(
\int_{r-h}^{r}
\norm{\bF^i_{s-\tau}}^2\,\dd s
\right)^{p/2}
\right].
\end{align*}
We apply H\"older's inequality in time once more and use
\eqref{D-5}$_2$ to obtain
\[
\begin{aligned}
	&\E\left[
	\frac1N\sum_{i=1}^N
	\left\|
	\int_{r-h}^{r}
	\sqrt{\sigma^i}\,
	\Diag(\bF^i_{s-\tau})\dd {\bf B}^i_s
	\right\|^p
	\right]
	\\
	&\hspace{1cm}\leq
	C_p (\sigma^{\infty})^{p/2}h^{p/2-1}
	\int_{r-h}^{r}{\mathcal D}_p(s-\tau)\,\dd s
	\\
	&\hspace{1cm}=
	C_p (\sigma^{\infty})^{p/2}h^{p/2-1}
	\int_{r-h-\tau}^{r-\tau}{\mathcal D}_p(s)\,\dd s.
\end{aligned}
\]
We combine two estimates to get \eqref{i-5}. \newline

\noindent Next, we take \(r=t\) and \(h=\tau\), and observe that for \(p\geq2\) and \(0\leq\tau,\sigma^{\infty} \leq1\), it holds
\[
\tau^{p-1}\leq\tau,
\qquad
(\sigma^{\infty})^{p/2}\tau^{p/2-1}\leq \sigma^{\infty},
\]
so we verify \eqref{D-6}. \newline

\noindent (ii)~(Derivation of \eqref{D-7}):~Next we estimate the weighted increment in \eqref{D-7}. Going back to \eqref{i-6}, with
\[
\|a+b\|^2\leq 2(\|a\|^2+\|b\|^2)
\]
and the Cauchy--Schwarz inequality in time,
\[
\left\|
\int_{t-\tau}^{t}\bF^i_{s-\tau}\,\dd s
\right\|^2
\leq
\tau\int_{t-\tau}^{t}
\|\bF^i_{s-\tau}\|^2\,\dd s,
\]
we have
\begin{align}
	\begin{aligned} \label{i-7}
		&\E\left[
		\frac1N\sum_{i=1}^N
		\norm{\hat{\bv}^i_{t-\tau}}^{p-2}
		\norm{\hat{\bv}^i_t-\hat{\bv}^i_{t-\tau}}^2
		\right]   \\
		& \hspace{1cm} \leq
		2\tau\E\left[
		\frac1N\sum_{i=1}^N
		\norm{\hat{\bv}^i_{t-\tau}}^{p-2}
		\int_{t-\tau}^{t}\|\bF^i_{s-\tau}\|^2\,\dd s
		\right]  \\
		& \hspace{1cm} +
		2\E\left[
		\frac1N\sum_{i=1}^N
		\norm{\hat{\bv}^i_{t-\tau}}^{p-2}
		\left\|
		\int_{t-\tau}^{t}
		\sqrt{\sigma^i}\,
		\Diag(\bF^i_{s-\tau})\dd {\bf B}^i_s
		\right\|^2
		\right] \\
		& \hspace{1cm} =:2\bigl(\tau I_{81}+I_{82}\bigr).	
	\end{aligned}
\end{align}
Below, we estimate the term \(I_{81}\) and \(I_{82}\) one by one. \newline

\noindent $\bullet$~Case 1 (Estimate of $I_{82}$):~Since the weight
\(\norm{\hat{\bv}^i_{t-\tau}}^{p-2}\) is unbounded, we use a truncation argument to estimate $I_{82}$. For \(R>0\), we set
\[
\hat{v}^{i,R}_{t-\tau} 
:=
\norm{\hat{\bv}^i_{t-\tau}}^{p-2}\wedge R.
\]
For \(s\in[t-\tau,t]\), both \(\hat{v}^{i,R}_{t-\tau}\) and $\bF^i_{s-\tau}$
are \(\mathcal F_{t-\tau}\)-measurable. Hence, we use the conditional It\^o
isometry to get 
\[
\begin{aligned}
	&\E\left[
	\hat{v}^{i,R}_{t-\tau}
	\left\|
	\int_{t-\tau}^{t}
	\sqrt{\sigma^i}\,
	\Diag(\bF^i_{s-\tau})\dd {\bf B}^i_s
	\right\|^2
	\right]=
	\E\left[
	\hat{v}^{i,R}_{t-\tau}
	\int_{t-\tau}^{t}
	\sigma^i \norm{\bF^i_{s-\tau}}^2\,\dd s
	\right].
\end{aligned}
\]
Letting \(R\to\infty\) and using the monotone convergence theorem, we
obtain
\begin{align}
	\begin{aligned}  \label{i-7-1}
		I_{82} &= \E\left[ \frac1N\sum_{i=1}^N
		\norm{\hat{\bv}^i_{t-\tau}}^{p-2}
		\left\|
		\int_{t-\tau}^{t}
		\sqrt{\sigma^i}\,
		\Diag(\bF^i_{s-\tau})\dd {\bf B}^i_s
		\right\|^2
		\right] \\
		&=
		\E\left[ \frac1N\sum_{i=1}^N
		\int_{t-\tau}^{t}
		\sigma^i
		\norm{\hat{\bv}^i_{t-\tau}}^{p-2}
		\norm{\bF^i_{s-\tau}}^2\,\dd s
		\right] \\
		&\leq \sigma^{\infty} \E\left[ \frac1N\sum_{i=1}^N
		\int_{t-\tau}^{t}
		\norm{\hat{\bv}^i_{t-\tau}}^{p-2}
		\norm{\bF^i_{s-\tau}}^2\,\dd s
		\right].
	\end{aligned}
\end{align}

\vspace{0.2cm}

\noindent $\bullet$~Case 2 (Estimate of $I_{81}$):~
\begin{equation} \label{i-7-2}
	I_{81}
	=
	\E\left[
	\frac1N\sum_{i=1}^N
	\norm{\hat{\bv}^i_{t-\tau}}^{p-2}
	\int_{t-\tau}^{t}
	\|\bF^i_{s-\tau}\|^2\,\dd s
	\right].
\end{equation}
In \eqref{i-7}, we combine estimates \eqref{i-7-1} and \eqref{i-7-2} to get 
\begin{align}
	\begin{aligned} \label{i-8}
		& \E\left[
		\frac1N\sum_{i=1}^N
		\norm{\hat{\bv}^i_{t-\tau}}^{p-2}
		\norm{\hat{\bv}^i_t-\hat{\bv}^i_{t-\tau}}^2
		\right]  \!\!\leq
		C_p(\tau+\sigma^{\infty})
		\int_{t-\tau}^{t}
		\E\left[
		\frac1N\sum_{i=1}^N
		\norm{\hat{\bv}^i_{t-\tau}}^{p-2}
		\norm{\bF^i_{s-\tau}}^2
		\right]\dd s.
	\end{aligned}
\end{align}
The remaining difficulty is that the weight in \eqref{i-8} is
evaluated at time \(t-\tau\), whereas \(\bF^i_{s-\tau}\) is generated by the
velocity at time \(s-\tau\). To match these two instants, we set 
\[
\hat{\bv}_{t-\tau}
=
\hat{\bv}_{s-\tau}+\bDelta_{t,s},
\qquad
\bDelta_{t,s}
:=
\hat{\bv}_{t-\tau}-\hat{\bv}_{s-\tau}.
\]
For \(p\geq2\), we have
\[
\norm{\hat{\bv}^i_{t-\tau}}^{p-2}
\leq
C_p\norm{\hat{\bv}^i_{s-\tau}}^{p-2}
+
C_p\norm{\bDelta^i_{t,s}}^{p-2}.
\]
Then, the integrand on the right-hand side of \eqref{i-8} becomes 
\begin{align*}
		& \E\left[ \frac1N\sum_{i=1}^N
		\norm{\hat{\bv}^i_{t-\tau}}^{p-2}
		\norm{\bF^i_{s-\tau}}^2
		\right] \\
		& \hspace{0.5cm} \leq C_p \Bigg(
		\E\left[ \frac1N\sum_{i=1}^N
		\norm{\hat{\bv}^i_{s-\tau}}^{p-2}
		\norm{\bF^i_{s-\tau}}^2
		\right]
		+
		\E\left[ \frac1N\sum_{i=1}^N
		\norm{\bDelta^i_{t,s}}^{p-2}
		\norm{\bF^i_{s-\tau}}^2
		\right]
		\Bigg) \\
		& \hspace{0.5cm} =: C_p \Big( I_{91} + I_{92} \Big).		
\end{align*}
Below, we estimate the terms $I_{9i},\,i=1,2$ one by one. \newline

\noindent $\diamond$~Case 2.1 (Estimate of $I_{91}$):~By \eqref{D-5}$_3$, one has 
\[
I_{91} = \E\left[ \frac1N\sum_{i=1}^N
\norm{\hat{\bv}^i_{s-\tau}}^{p-2}
\norm{\bF^i_{s-\tau}}^2
\right]  \leq
C_p{\mathcal D}_p(s-\tau).
\]

\noindent $\diamond$~Case 2.2 (Estimate of $I_{92}$):~We use Young's inequality to find 
\[
\norm{\bDelta^i_{t,s}}^{p-2}
\norm{\bF^i_{s-\tau}}^2
\leq
C_p \Big( \norm{\bF^i_{s-\tau}}^p
+ \norm{\bDelta^i_{t,s}}^p \Big).
\]
By \eqref{D-5}$_2$, we have
\[
\E\left[
\frac1N\sum_{i=1}^N\norm{\bF^i_{s-\tau}}^p
\right]
\leq
C_p{\mathcal D}_p(s-\tau).
\]
Now, we need to estimate \(\bDelta_{t,s}\). For
\(s\in[t-\tau,t]\), we apply \eqref{i-5} with
$r=t-\tau$ and
$h=t-s$ to get
\[
\begin{aligned}
	\E\left[
	\frac1N\sum_{i=1}^N
	\norm{\bDelta^i_{t,s}}^p
	\right]
	&\leq
	C_p\left(
	(t-s)^{p-1}
	+
	(\sigma^\infty)^{p/2}(t-s)^{p/2-1}
	\right)
	\int_{s-2\tau}^{t-2\tau}
	{\mathcal D}_p(r)\,\dd r
	\\
	&\leq
	C_p(\tau+\sigma^\infty)
	\int_{s-2\tau}^{t-2\tau}
	{\mathcal D}_p(r)\,\dd r.
\end{aligned}
\]
Therefore, we have
\begin{equation}
	\begin{aligned}
		&\E\left[
		\frac1N\sum_{i=1}^N
		\norm{\hat{\bv}^i_{t-\tau}}^{p-2}
		\norm{\bF^i_{s-\tau}}^2
		\right]
	\leq
		C_p{\mathcal D}_p(s-\tau)
		+
		C_p(\tau+\sigma^\infty)
		\int_{s-2\tau}^{t-2\tau}
		{\mathcal D}_p(r)\,\dd r.
	\end{aligned}
	\label{i-9}
\end{equation}
Next, we substitute \eqref{i-9} into \eqref{i-8} to find
\[
\begin{aligned}
	&\E\left[
	\frac1N\sum_{i=1}^N
	\norm{\hat{\bv}^i_{t-\tau}}^{p-2}
	\norm{\hat{\bv}^i_t-\hat{\bv}^i_{t-\tau}}^2
	\right]
	\\
	&\qquad\leq
	C_p(\tau+\sigma^\infty)
	\int_{t-2\tau}^{t-\tau}
	{\mathcal D}_p(s)\,\dd s
	+
	C_p(\tau+\sigma^\infty)^2
	\int_{t-\tau}^{t}
	\int_{s-2\tau}^{t-2\tau}
	{\mathcal D}_p(r)\,\dd r\,\dd s
	\\
	&\qquad\leq
	C_p(\tau+\sigma^\infty)
	\int_{t-2\tau}^{t-\tau}
	{\mathcal D}_p(s)\,\dd s
	+
	C_p(\tau+\sigma^\infty)^2\tau
	\int_{t-3\tau}^{t-2\tau}
	{\mathcal D}_p(r)\,\dd r\\
	&\qquad\leq C_p(\tau+\sigma^\infty)
	\int_{t-3\tau}^{t-\tau}
	{\mathcal D}_p(s)\,\dd s,
\end{aligned}
\]
where we used
\[
s-2\tau\geq t-3\tau,
\qquad
s\in[t-\tau,t],
\]
together with
\[
(\tau+\sigma^\infty)^2\tau
\leq
2(\tau+\sigma^\infty),
\]
valid for \(0\leq\tau,\,\sigma^\infty\leq1\).
This verifies \eqref{D-7}.

\subsection{Martingale verification for high-order moment estimate}
\label{App-A.3}
In this subsection, we show that the local martingale in It\^o's formula for $ \frac1N\sum_{i=1}^N\norm{\hat{\bv}^i_t}^p$ has zero mean. \newline

For fixed \(T>0\), we first apply the
Burkholder--Davis--Gundy inequality to the stopped local martingale,
and then use H\"older's inequality in probability and
H\"older's inequality in time to obtain
\begin{align}
\begin{aligned} \label{i-10}
	& \E\left[\sup_{0\le t\le T} \left |   \frac{p}{N}\sum_{i=1}^N\sqrt{\sigma^i}
	\int_0^t \Diag(\bF^i_{s-\tau}) \| \hat{\bv}^i_s \|^{p-2} \hat{\bv}^i_s \cdot\dd {\bf B}^i_s \right | \right] \\
	& \hspace{0.5cm} \le
	C_{p,N,\sigma}\sum_{i=1}^N
	\E\left[
	\sup_{0\le s\le T}\norm{\hat{\bv}^i_s}^{p-1}
	\left(\int_0^T\norm{\bF^i_{s-\tau}}^2\dd s\right)^{1/2}
	\right] \\
	&\hspace{0.5cm} \le
	C_{p,N,\sigma,T}\sum_{i=1}^N
	\left(\E\left[\sup_{0\le s\le T}\norm{\hat{\bv}^i_s}^p\right]\right)^{(p-1)/p}
	\left(\E \left[ \int_0^T\norm{\bF^i_{s-\tau}}^p\dd s \right] \right)^{1/p}.
\end{aligned}
\end{align}
The estimate follows by letting the stopping level tend to infinity
and using Fatou's lemma.
Next, we study each factor in the right-hand side of \eqref{i-10}. \newline

\noindent $\bullet$~Case 1 (Estimate of the first term): Since \(N\) is fixed, Lemma~\ref{L4.3} implies
\[ \E\left[\sup_{0\le s\le T}\norm{\hat{\bv}^i_s}^p\right] < \infty. \]
\noindent $\bullet$~Case 2 (Estimate of the second term): First, we assume that \(\tau>0\). For \(s\in[0,\tau]\),
\[
\bF_{s-\tau}^i
=
\frac1N\sum_{j=1}^N
\omega_{s-\tau}^{ij}
\left[
\bigl(\bv_{s-\tau}^j-\bar{\bv}_{s-\tau}\bigr)
-
\bigl(\bv_{s-\tau}^i-\bar{\bv}_{s-\tau}\bigr)
\right].
\]
Hence, Lemma~\ref{L4.2}(i) and the centered history bound imply
\[
\int_0^\tau
\E\left[
\frac1N\sum_{i=1}^N\norm{ \bF^i_{s-\tau}}^p
\right]\dd s<\infty.
\]
On each subsequent delay interval, Lemma~\ref{L4.2}(i) and
Lemma~\ref{L4.3} yield the same integrability. Finally, we sum the above estimates over the finitely
many intervals intersecting \([0,T]\) to get 
\[
\E \left[ \int_0^T
\frac1N\sum_{i=1}^N\norm{\bF^i_{s-\tau}}^p\dd s \right]<\infty.
\]
For \(\tau=0\), this follows directly from Lemma~\ref{L4.2}(i) and
Lemma~\ref{L4.3}. Combining the two estimates above with \eqref{i-10} and using that \(N\) is fixed, we have
\[
\E\left[\sup_{0\le t\le T}\left|\frac{p}{N}\sum_{i=1}^N\sqrt{\sigma^i}\int_0^t \Diag(\bF^i_{s-\tau}) \| \hat{\bv}^i_s \|^{p-2} \hat{\bv}^i_s \cdot\dd {\bf B}^i_s \right| \right]  <\infty.
\]
Thus the stochastic integral
\[
\frac{p}{N}\sum_{i=1}^N\sqrt{\sigma^i}
\int_0^t
\Diag(\bF^i_{s-\tau})
\|\hat{\bv}_s^i\|^{p-2}\hat{\bv}_s^i
\cdot\dd{\bf B}_s^i
\] is a true martingale and 
\[  \E \left[    \frac{p}{N}\sum_{i=1}^N\sqrt{\sigma^i}
\int_0^t \Diag(\bF^i_{s-\tau}) \| \hat{\bv}^i_s \|^{p-2} \hat{\bv}^i_s \cdot\dd {\bf B}^i_s \right] = 0. \]
The verification may depend on the
fixed \(N\), but no such constant enters the preceding \(N\)-uniform a priori estimates.

\end{document}